% 
% New Complex and Quaternion-Hyperbolic Reflection Groups
%
% Daniel Allcock

%					
%	dja.tex	
%	Daniel Allcock   (allcock@math.utah.edu)
%	17 December 1996			
%						
%	A simple set of macros for plain tex
%	documentation is in dja.doc, a text file
%
% suppresses duplicate inputs of this file:
\ifx\djatexLoaded\relax\endinput\else\let\djatexLoaded=\relax\fi

% FONTS
%
% title font (scaled, not a larger font, for portability)
%\font\titlefont=cmb10 scaled \magstep3
\def\titlefont{\bf}
% calligraphic font
\newfam\calfont
\font\tencal=eusm10 \font\sevencal=eusm7 \font\fivecal=eusm5
\textfont\calfont=\tencal 
\scriptfont\calfont=\sevencal 
\scriptscriptfont\calfont=\fivecal
\def\cal{\fam=\calfont}	
% fraktur
\newfam\frakturfont
\font\tenfrak=eufm10 \font\sevenfrak=eufm7 \font\fivefrak=eufm5
\textfont\frakturfont=\tenfrak 
\scriptfont\frakturfont=\sevenfrak 
\scriptscriptfont\frakturfont=\fivefrak

% blackboard bold (includes many math symbols)
\newfam\bbbfont
\font\tenbbb=msbm10 \font\sevenbbb=msbm7 \font\fivebbb=msbm5	
\textfont\bbbfont=\tenbbb 
\scriptfont\bbbfont=\sevenbbb 
\scriptscriptfont\bbbfont=\fivebbb
\def\bbb{\fam=\bbbfont}
% script and script-script fonts for typewriter font (all in
% 10pt for portability)
\scriptfont\ttfam=\tentt
\scriptscriptfont\ttfam=\tentt

% PROGRAMMING STUFF
% 
% `if not defined' :
\long\def\ifndef#1{\expandafter\ifx\csname#1\endcsname\relax}
% sometimes more convenient, but you can't say
% ``\ifdef...\else...\fi''
\long\def\ifdef#1{\ifndef{#1}\else}
% My warning messages
\def\mywarning#1{\ifnum\warningsoff=0
	\immediate\write16{l.\the\inputlineno: #1}\fi}
\ifndef{warningsoff}\def\warningsoff{0}\fi

% FRONT MATTER
%
% note: \author, \title and \note get redefined when
% \bibliography is called.
\long\def\author#1{\def\Zauthor{#1}}
\long\def\title#1{\def\Ztitle{#1}}
\long\def\address#1{\def\Zaddress{#1}}
\long\def\email#1{\def\Zemail{#1}}
\long\def\homepage#1{\def\Zhomepage{#1}}
\long\def\date#1{\def\Zdate{#1}}
\long\def\subject#1{\def\Zsubject{#1}}
\long\def\note#1{\def\Znote{#1}}
\def\plaintitlepage{{\parindent=0pt
	\ifdef{Ztitle}\par{\titlefont\Ztitle}\medskip\fi
	\ifdef{Zauthor}\par{\Zauthor}\fi
	\ifdef{Zdate}\par{\Zdate}\fi
	\smallskip
	\ifdef{Zemail}\par{\it \Zemail}\fi
	\ifdef{Zhomepage}\par{web page: \it\Zhomepage}\fi
	\ifdef{Zaddress}\par{\Zaddress}\fi
	\ifdef{Zsubject}\smallskip\par{1991 mathematics subject
		classification: \Zsubject}\fi
	\ifndef{Znote}\else\smallskip\par{\Znote}\fi
	\vskip 0pt plus 10pt}}
\def\abstract{\bigbreak\noindent{\bf Abstract}\smallskip\noindent} 

% SECTIONS, THEOREMS, ETC.
%
%\let\section=\beginsection
\outer\def\section#1\par{\bigbreak\noindent{\bf #1}\nobreak\medskip\noindent}
% Proclaimed Things (theorems, lemmas, etc.)
%   (modeled on plain tex's ``\proclaim'', but theorems can be more
%   than one paragraph long.)
\outer\def\beginproclaim#1. {\medbreak
	\noindent{\bf #1.\enspace}\begingroup\sl}
\def\endproclaim{\endgroup\par\ifdim\lastskip<\medskipamount 
	\removelastskip\penalty55\medskip\fi}
% Proofs (n.b. \drawbox is defined in `OTHER')
\def\beginproof#1{\smallskip{\it#1\/}}
\def\endproof{\leavevmode\vrule height0pt width0pt
	depth0pt\nobreak\hfill\proofbox\smallskip}
\def\QED{\endproof}
\def\proofbox{\drawbox{1.2ex}{1.2ex}{.1ex}}
% Remarks (n.b. remarks are only one par long)
\def\remark#1#2\par{\ifdim\lastskip<\smallskipamount
	\removelastskip\penalty55\smallskip\fi
	{\it #1\/}#2\smallbreak}

% CROSS-REFERENCES AND BIBLIOGRAPHY
%
% If the tag has not been defined (say by \deftag), then
% \tag warns and inserts question
% marks. If the tag is defined, it
% just produces the tag (symbolic if \userawtags=0, 
% or symbolic:raw, otherwise.) Default is \userawtags=1.

\def\symbolictags{\def\userawtags{0}}
\ifndef{userawtags}\def\userawtags{1}\fi
\def\tag#1{\ifndef{#1}\mywarning{tag `#1' 
	undefined.}{\hbox{\bf?`}\tt #1\hbox{\bf ?}}\else
	\ifnum\userawtags=0 \csname #1\endcsname
	\else\csname #1\endcsname{ \tt [#1]}\fi\fi}
\def\Tag#1{\tag{#1}}
\def\eqtag#1{(\tag{#1})}
\def\eqTag#1{\ifnum\userawtags=0
	\eqtag{#1}\else
	\lower12pt\hbox{\eqtag{#1}}\fi}
\def\cite#1{[\tag{#1}]}
\def\ecite#1#2{[\tag{#1}, #2]}
\def\nocite#1{}
% syntax: \deftag{label to insert in text}{label in file}
\def\deftag#1#2{\ifndef{#2}\else\mywarning{tag `#2' defined
	more than once.}\fi 
	\expandafter\def\csname #2\endcsname{#1}}
\def\defcite#1#2{\deftag{#1}{#2}}
% bibliography entries
\def\bibitem#1{\ifnum\userawtags=0
	\item{[{\tag{#1}}]}\else
	\noindent\hangindent=\parindent\hangafter=1
			\cite{#1}\enskip\fi}

% FONT SHORTCUTS 
%
\def\rom#1{({\it\romannumeral#1\/})}
\def\Rom#1{{\rm\uppercase\expandafter{\romannumeral#1}}}
\def\a{\alpha}		
  	
  \def\cals{{\cal S}}

\def\cale{{\cal E}}  
  
\def\calg{{\cal G}}  
\def\calh{{\cal H}}  
 \def\calr{{\cal R}}

% SHORTCUTS AND MNEMONICS FOR STANDARD MATH SYMBOLS
%
\def\Z{{\bbb Z}} % integers
\def\Q{{\bbb Q}} % rational numbers 
\def\R{{\bbb R}} % real numbers
\def\C{{\bbb C}} % complex numbers
\def\F{{\bbb F}} % A finite field
\let\sset=\subseteq		\let\dimension=\dim%
\let\congruent=\equiv%
\def\setminus{\mathop{\bbb \char"72}\nolimits}

\let\tensor=\otimes	

\let\isomorphism=\cong		
\let\to=\rightarrow
\def\aut{\mathop{\rm Aut}\nolimits}

\def\re{\mathop{\rm Re}\nolimits}
\def\im{\mathop{\rm Im}\nolimits}

\def\mod{\mathop{\rm mod}\nolimits}
% quick parentheses and brackets
\def\({\left(} \def\){\right)}
\def\[{\left[} \def\]{\right]}

% OTHER STUFF	
%
\def\hideboxes{\overfullrule=0pt}
\def\frac#1/#2{\leavevmode\kern.1em\raise.5ex\hbox{\the
	\scriptfont0 #1}\kern-.1em/\kern-.15em\lower.25ex
	\hbox{\the\scriptfont0 #2}} 
% 2x2 matrices for use in text:
\def\smallmatrix#1#2#3#4{%
	\left({#1\atop #3}\;{#2\atop #4}\right)}
% syntax: ``\drawbox{internal width}{internal height}{line
% 	thickness}''; arguments should be dimens.
\def\drawbox#1#2#3{\leavevmode\vbox{\hrule height #3%
	\hbox{\vrule width #3 height #2\kern #1%
	\vrule width #3}\hrule height #3}}

% does an \halign but centers it with a little vertical space
% around it. Uses diplay math mode.

% an \llap that you can use inside formulas without
% tex changing to textstyle inside the \hbox of \llap
\def\mathllap#1{\mathchoice
{\llap{$\displaystyle #1$}}%
{\llap{$\textstyle #1$}}%
{\llap{$\scriptstyle #1$}}%
{\llap{$\scriptscriptstyle #1$}}}
% the set of all #1 such that #2; makes the braces and the
% middle bar big enough to accomodate everything.
\def\set#1#2{\left\{\,#1\mathllap{\phantom{#2}}\mathrel{}\right|\left.#2\mathllap{\phantom{#1}}\,\right\}}

\parskip=0pt

\catcode`!=11 %  ***** THIS MUST NEVER BE OMITTED
% *******************************
% *** HACKS  (Utility macros) ***
% *******************************
 
% ** User commands
% **   \PiC{P\kern-.12em\lower.5ex\hbox{I}\kern-.075emC}
% **   \PiCTeX{\PiC\kern-.11em\TeX}
% **   \placevalueinpts of <DIMENSION REGISTER> in {CONTROL SEQUENCE}
  
% ** Internal commands
% **   \!ifnextchar{CHARACTER}{TRUE ACTION}{FALSE ACTION}
% **   \!tfor NAME := LIST \do {BODY}
% **   \!etfor NAME:= LIST \do {BODY}
% **   \!cfor NAME := LIST \do {BODY}
% **   \!ecfor NAME:= LIST \do {BODY}
% **   \!ifempty{MACRO}{TRUE ACTION}{FALSE ACTION}
% **   \!getnext\\ITEMfrom\LIST
% **   \!getnextvalueof\DIMEN\from\LIST
% **   \!copylist\LISTMACRO_A\to\LISTMACRO_B
% **   \!wlet\CONTROL_SEQUENCE_A=\CONTROL_SEQUENCE_B
% **   \!listaddon ITEM LIST
% **   \!rightappendITEM\withCS\to\LISTMACRO
% **   \!leftappendITEM\withCS\to\LISTMACRO
% **   \!lop\LISTMACRO\to\ITEM
% **   \!loop ... repeat
% **   \!!loop ... repeat
% **   \!removept{DIMENSION REGISTER}{CONTROL SEQUENCE}
% **   \!mlap{...}  
% **   \!vmlap{...}
% **   \!not{TEK if-CONDITION}

% ** First, here are the the PiCTeX logo, and the syllable PiC:
\def\PiC{P\kern-.12em\lower.5ex\hbox{I}\kern-.075emC}
\def\PiCTeX{\PiC\kern-.11em\TeX}

% ** The following macro expands to parameter #2 or parameter #3 according to
% ** whether the next non-blank character following the macro is or is not #1. 
% ** Blanks following the macro are gobbled.
\def\!ifnextchar#1#2#3{%
  \let\!testchar=#1%
  \def\!first{#2}%
  \def\!second{#3}%
  \futurelet\!nextchar\!testnext}
\def\!testnext{%
  \ifx \!nextchar \!spacetoken 
    \let\!next=\!skipspacetestagain
  \else
    \ifx \!nextchar \!testchar
      \let\!next=\!first
    \else 
      \let\!next=\!second 
    \fi 
  \fi
  \!next}
\def\\{\!skipspacetestagain} 
  \expandafter\def\\ {\futurelet\!nextchar\!testnext} 
\def\\{\let\!spacetoken= } \\  %  ** set \spacetoken to a space token

% ** Borrow the "tfor" macro from Latex:
% **   \!tfor NAME := LIST \do {BODY}
% **   if, before expansion, LIST = T1 ... Tn,  where each  Ti  is a token
% **   or  {...},  then executes  BODY  n  times, with  NAME = Ti  on the
% **   i-th iteration.  Works for  n=0.
\def\!tfor#1:=#2\do#3{%
  \edef\!fortemp{#2}%
  \ifx\!fortemp\!empty 
    \else
    \!tforloop#2\!nil\!nil\!!#1{#3}%
  \fi}
\def\!tforloop#1#2\!!#3#4{%
  \def#3{#1}%
  \ifx #3\!nnil
    \let\!nextwhile=\!fornoop
  \else
    #4\relax
    \let\!nextwhile=\!tforloop
  \fi 
  \!nextwhile#2\!!#3{#4}}

% **   \!etfor NAME:= LIST\do {BODY}
% **   This is like \!cfor, but LIST is any balanced token list whose complete
% **     expansion has the form  T1 ... Tn
\def\!etfor#1:=#2\do#3{%
  \def\!!tfor{\!tfor#1:=}%
  \edef\!!!tfor{#2}%
  \expandafter\!!tfor\!!!tfor\do{#3}}

% **   modify the Latex \tfor (token-for) loop to a \cfor (comma-for) loop.
% **   \!cfor NAME := LIST \do {BODY}
% **     if, before expansion, LIST = a1,a2,...an, then executes  BODY n times,
% **     with  NAME = ai  on the i-th iteration.  Works for  n=0.
\def\!cfor#1:=#2\do#3{%
  \edef\!fortemp{#2}%
  \ifx\!fortemp\!empty 
  \else
    \!cforloop#2,\!nil,\!nil\!!#1{#3}%
  \fi}
\def\!cforloop#1,#2\!!#3#4{%
  \def#3{#1}%
  \ifx #3\!nnil
    \let\!nextwhile=\!fornoop 
  \else
    #4\relax
    \let\!nextwhile=\!cforloop
  \fi
  \!nextwhile#2\!!#3{#4}}

% **   \!ecfor NAME:= LIST\do {BODY}
% **   This is like \!cfor, but LIST is any balanced token list whose complete
% **     expansion has the form  a1,a2,...,an.
\def\!ecfor#1:=#2\do#3{%
  \def\!!cfor{\!cfor#1:=}%
  \edef\!!!cfor{#2}%
  \expandafter\!!cfor\!!!cfor\do{#3}}

\def\!empty{}
\def\!nnil{\!nil}
\def\!fornoop#1\!!#2#3{}

% **  \!ifempty{ARG}{TRUE ACTION}{FALSE ACTION}
\def\!ifempty#1#2#3{%
  \edef\!emptyarg{#1}%
  \ifx\!emptyarg\!empty
    #2%
  \else
    #3%
  \fi}
 
% **  \!getnext\\ITEMfrom\LIST
% **    \LIST has the form \\{item1}\\{item2}\\{item3}...\\{itemk}
% **    This routine sets \ITEM to item1, and cycles \LIST to
% **    \\{item2}\\{item3}...\\{itemk}\\{item1}
\def\!getnext#1\from#2{%
  \expandafter\!gnext#2\!#1#2}%
\def\!gnext\\#1#2\!#3#4{%
  \def#3{#1}%
  \def#4{#2\\{#1}}%
  \ignorespaces}

% ** \!getnextvalueof\DIMEN\from\LIST
% **   Similar to !getnext.  
% **   \LIST has the form \\{dimen1}\\{dimen2}\\{dimen3} ... 
% **   \DIMEN is a dimension register
% **   Works also for counts
%
\def\!getnextvalueof#1\from#2{%
  \expandafter\!gnextv#2\!#1#2}%
\def\!gnextv\\#1#2\!#3#4{%
  #3=#1%
  \def#4{#2\\{#1}}%
  \ignorespaces}

% ** \!copylist\LISTMACROA\to\LISTMACROB
% **   makes the replacement text of LISTMACRO B identical to that of
% **   list macro A.
\def\!copylist#1\to#2{%
  \expandafter\!!copylist#1\!#2}
\def\!!copylist#1\!#2{%
  \def#2{#1}\ignorespaces}

% **  \!wlet\CSA=\CSB
% **  lets control sequence \CSB = control sequence \CSA, and writes a
% **    message to that effect in the log file using plain TEK's \wlog
\def\!wlet#1=#2{%
  \let#1=#2 
  \wlog{\string#1=\string#2}}
 
% ** \!listaddon ITEM LIST
% ** LIST <-- LIST \\ ITEM
\def\!listaddon#1#2{%
  \expandafter\!!listaddon#2\!{#1}#2}
\def\!!listaddon#1\!#2#3{%
  \def#3{#1\\#2}}
 
% **  \!rightappendITEM\to\LISTMACRO
% **    \LISTMACRO --> \LISTMACRO\\{ITEM}
%\def\!rightappend#1\to#2{\expandafter\!!rightappend#2\!{#1}#2}
%\def\!!rightappend#1\!#2#3{\def#3{#1\\{#2}}}

% **  \!rightappendITEM\withCS\to\LISTMACRO
% **    \LISTMACRO --> \LISTMACRO||CS||{ITEM}
\def\!rightappend#1\withCS#2\to#3{\expandafter\!!rightappend#3\!#2{#1}#3}
\def\!!rightappend#1\!#2#3#4{\def#4{#1#2{#3}}}

% **  \!leftappendITEM\withCS\to\LISTMACRO
% **    \LISTMACRO --> CS||{ITEM}||\LISTMACRO
\def\!leftappend#1\withCS#2\to#3{\expandafter\!!leftappend#3\!#2{#1}#3}
\def\!!leftappend#1\!#2#3#4{\def#4{#2{#3}#1}}

% **  \!lop\LISTMACRO\to\ITEM
% **    \\{item1}\\{item2}\\{item3} ... --> \\{item2}\\{item3} ...
% **    item1 --> \ITEM
\def\!lop#1\to#2{\expandafter\!!lop#1\!#1#2}
\def\!!lop\\#1#2\!#3#4{\def#4{#1}\def#3{#2}}

% **  \!placeNUMBER\of\LISTMACRO\in\ITEM
% **    the NUMBERth item of \LISTMACRO --> replacement text of \ITEM
%\def\!place#1\of#2\in#3{\def#3{\outofrange}%
%{\count0=#1\def\\##1{\advance\count0-1 \ifnum\count0=0 \gdef#3{##1}\fi}#2}}

% **  Following code converts a commalist to a list macro, with all items 
% **    fully expanded.
%\!ecfor\item:=\commalist\do{\expandafter\!rightappend\item\to\list}

% ** \!loop ... repeat
% ** This is exactly like TEX's \loop ... repeat.  It can be used in nesting
% ** two loops, without puting the inner one inside a group.
\def\!loop#1\repeat{\def\!body{#1}\!iterate}
\def\!iterate{\!body\let\!next=\!iterate\else\let\!next=\relax\fi\!next}
 
% ** \!!loop ... repeat
% ** This is exactly like TEX's \loop ... repeat.  It can be used in nesting
% ** two loops, without puting the inner one inside a group.
\def\!!loop#1\repeat{\def\!!body{#1}\!!iterate}
\def\!!iterate{\!!body\let\!!next=\!!iterate\else\let\!!next=\relax\fi\!!next}
%  (\multiput uses \!!loop)
 
% ** \!removept{DIMENREG}{\CS}
% ** Defines the control sequence CS to be the value (in points) in the
% ** dimension register DIMENREG (but without the "pt" TEK usually adds)
% ** E.g., after  \dimen0=12.3pt \!removept\dimen0\A, \A expands to 12.3
\def\!removept#1#2{\edef#2{\expandafter\!!removePT\the#1}}
{\catcode`p=12 \catcode`t=12 \gdef\!!removePT#1pt{#1}}

% ** \pladevalueinpts of <DIMENSION REGISTER> in {CONTROL SEQUENCE}
\def\placevalueinpts of <#1> in #2 {%
  \!removept{#1}{#2}}
 
% ** \!mlap{...}  \!vmlap{...}
% ** Center  ...  in a box of width 0.
\def\!mlap#1{\hbox to 0pt{\hss#1\hss}}
\def\!vmlap#1{\vbox to 0pt{\vss#1\vss}}
 
% ** \!not{TEK if-CONDITION}
% ** By a TEK if-CONDITION is meant something like 
% **     \ifnum\N<0,   or   \ifdim\A>\B
% ** \!not produces an if-condition which is false if the original condition
% ** is true, and true if the original condition is false.
\def\!not#1{%
  #1\relax
    \!switchfalse
  \else
    \!switchtrue
  \fi
  \if!switch
  \ignorespaces}

% *******************
% *** ALLOCATIONS ***
% *******************

% This section allocates all the registers PiCTeX uses. Following
% each allocation is a string of the form  ....N.D...L......... ;
% the various letters show which sections of PiCTeX make explicit
% reference to that register, according to the following code:
 
%   H Hacks
%   A Areas
%   W arroWs
%   B Bars
%   X boXes
%   C Curves
%   D Dashpattterns
%   V diVision
%   E Ellipses
%   U rUles
%   L Linear arc
%   G loGten
%   P Pictures
%   O plOtting
%   Y pYthagoras
%   Q Quadratic arc
%   R Rotations
%   S Shading
%   T Ticks

% Turn off messages from TeX's allocation macros
\let\!!!wlog=\wlog              % "\wlog" is defined in plain TeX
\def\wlog#1{}    

\newdimen\headingtoplotskip     %.A.................
\newdimen\linethickness         %.A..X....U........T
\newdimen\longticklength        %.A................T
\newdimen\plotsymbolspacing     %......D...L....Q...
\newdimen\shortticklength       %.A................T
\newdimen\stackleading          %.A..........P......
\newdimen\tickstovaluesleading  %.A................T
\newdimen\totalarclength        %......D...L....Q...
\newdimen\valuestolabelleading  %.A.................

\newbox\!boxA                   %.AW...............T
\newbox\!boxB                   %..W................
\newbox\!picbox                 %............P......
\newbox\!plotsymbol             %..........L..O.....
\newbox\!putobject              %............PO...S.
\newbox\!shadesymbol            %.................S.

\newcount\!countA               %.A....D..UL....Q.ST
\newcount\!countB               %......D..U.....Q.ST
\newcount\!countC               %...............Q..T
\newcount\!countD               %...................
\newcount\!countE               %.............O....T
\newcount\!countF               %.............O....T
\newcount\!countG               %..................T
\newcount\!fiftypt              %.........U.........
\newcount\!intervalno           %..........L....Q...
\newcount\!npoints              %..........L........
\newcount\!nsegments            %.........U.........
\newcount\!ntemp                %............P......
\newcount\!parity               %.................S.
\newcount\!scalefactor          %..................T
\newcount\!tfs                  %.......V...........
\newcount\!tickcase             %..................T

\newdimen\!Xleft                %............P......
\newdimen\!Xright               %............P......
\newdimen\!Xsave                %.A................T
\newdimen\!Ybot                 %............P......
\newdimen\!Ysave                %.A................T
\newdimen\!Ytop                 %............P......
\newdimen\!angle                %........E..........
\newdimen\!arclength            %..W......UL....Q...
\newdimen\!areabloc             %.A........L........
\newdimen\!arealloc             %.A........L........
\newdimen\!arearloc             %.A........L........
\newdimen\!areatloc             %.A........L........
\newdimen\!bshrinkage           %.................S.
\newdimen\!checkbot             %..........L........
\newdimen\!checkleft            %..........L........
\newdimen\!checkright           %..........L........
\newdimen\!checktop             %..........L........
\newdimen\!dimenA               %.AW.X.DVEUL..OYQRST
\newdimen\!dimenB               %....X.DVEU...O.QRS.
\newdimen\!dimenC               %..W.X.DVEU......RS.
\newdimen\!dimenD               %..W.X.DVEU....Y.RS.
\newdimen\!dimenE               %..W........G..YQ.S.
\newdimen\!dimenF               %...........G..YQ.S.
\newdimen\!dimenG               %...........G..YQ.S.
\newdimen\!dimenH               %...........G..Y..S.
\newdimen\!dimenI               %...BX.........Y....
\newdimen\!distacross           %..........L....Q...
\newdimen\!downlength           %..........L........
\newdimen\!dp                   %.A..X.......P....S.
\newdimen\!dshade               %.................S.
\newdimen\!dxpos                %..W......U..P....S.
\newdimen\!dxprime              %...............Q...
\newdimen\!dypos                %..WB.....U..P......
\newdimen\!dyprime              %...............Q...
\newdimen\!ht                   %.A..X.......P....S.
\newdimen\!leaderlength         %......D..U.........
\newdimen\!lshrinkage           %.................S.
\newdimen\!midarclength         %...............Q...
\newdimen\!offset               %.A................T
\newdimen\!plotheadingoffset    %.A.................
\newdimen\!plotsymbolxshift     %..........L..O.....
\newdimen\!plotsymbolyshift     %..........L..O.....
\newdimen\!plotxorigin          %..........L..O.....
\newdimen\!plotyorigin          %..........L..O.....
\newdimen\!rootten              %...........G.......
\newdimen\!rshrinkage           %.................S.
\newdimen\!shadesymbolxshift    %.................S.
\newdimen\!shadesymbolyshift    %.................S.
\newdimen\!tenAa                %...........G.......
\newdimen\!tenAc                %...........G.......
\newdimen\!tenAe                %...........G.......
\newdimen\!tshrinkage           %.................S.
\newdimen\!uplength             %..........L........
\newdimen\!wd                   %....X.......P....S.
\newdimen\!wmax                 %...............Q...
\newdimen\!wmin                 %...............Q...
\newdimen\!xB                   %...............Q...
\newdimen\!xC                   %...............Q...
\newdimen\!xE                   %..W.....E.L....Q.S.
\newdimen\!xM                   %..W.....E......Q.S.
\newdimen\!xS                   %..W.....E.L....Q.S.
\newdimen\!xaxislength          %.A................T
\newdimen\!xdiff                %..........L........
\newdimen\!xleft                %............P......
\newdimen\!xloc                 %..WB.....U.......S.
\newdimen\!xorigin              %.A........L.P....S.
\newdimen\!xpivot               %................R..
\newdimen\!xpos                 %..........L.P..Q.ST
\newdimen\!xprime               %...............Q...
\newdimen\!xright               %............P......
\newdimen\!xshade               %.................S.
\newdimen\!xshift               %..W.........PO...S.
\newdimen\!xtemp                %............P......
\newdimen\!xunit                %.AWBX...EUL.P..QRS.
\newdimen\!xxE                  %........E..........
\newdimen\!xxM                  %........E..........
\newdimen\!xxS                  %........E..........
\newdimen\!xxloc                %..WB....EU.........
\newdimen\!yB                   %...............Q...
\newdimen\!yC                   %...............Q...
\newdimen\!yE                   %..W.....E.L....Q...
\newdimen\!yM                   %..W.....E......Q...
\newdimen\!yS                   %..W.....E.L....Q...
\newdimen\!yaxislength          %.A................T
\newdimen\!ybot                 %............P......
\newdimen\!ydiff                %..........L........
\newdimen\!yloc                 %..WB.....U.......S.
\newdimen\!yorigin              %.A........L.P....S.
\newdimen\!ypivot               %................R..
\newdimen\!ypos                 %..........L.P..Q.ST
\newdimen\!yprime               %...............Q...
\newdimen\!yshade               %.................S.
\newdimen\!yshift               %..W.........PO...S.
\newdimen\!ytemp                %............P......
\newdimen\!ytop                 %............P......
\newdimen\!yunit                %.AWBX...EUL.P..QRS.
\newdimen\!yyE                  %........E..........
\newdimen\!yyM                  %........E..........
\newdimen\!yyS                  %........E..........
\newdimen\!yyloc                %..WB....EU.........
\newdimen\!zpt                  %.AWBX.DVEULGP.YQ.ST

\newif\if!axisvisible           %.A.................
\newif\if!gridlinestoo          %..................T
\newif\if!keepPO                %...................
\newif\if!placeaxislabel        %.A.................
\newif\if!switch                %H..................
\newif\if!xswitch               %.A................T

\newtoks\!axisLaBeL             %.A.................
\newtoks\!keywordtoks           %.A.................

\newwrite\!replotfile           %.............O.....

\newhelp\!keywordhelp{The keyword mentioned in the error message in unknown. 
Replace NEW KEYWORD in the indicated response by the keyword that 
should have been specified.}    %.A.................

% The following commands assign alternate names to some of the 
% above registers.  "\!wlet"  is defined in  Hacks.
\!wlet\!!origin=\!xM                   %.A................T
\!wlet\!!unit=\!uplength               %.A................T
\!wlet\!Lresiduallength=\!dimenG       %.........U.........
\!wlet\!Rresiduallength=\!dimenF       %.........U.........
\!wlet\!axisLength=\!distacross        %.A................T
\!wlet\!axisend=\!ydiff                %.A................T
\!wlet\!axisstart=\!xdiff              %.A................T
\!wlet\!axisxlevel=\!arclength         %.A................T
\!wlet\!axisylevel=\!downlength        %.A................T
\!wlet\!beta=\!dimenE                  %...............Q...
\!wlet\!gamma=\!dimenF                 %...............Q...
\!wlet\!shadexorigin=\!plotxorigin     %.................S.
\!wlet\!shadeyorigin=\!plotyorigin     %.................S.
\!wlet\!ticklength=\!xS                %..................T
\!wlet\!ticklocation=\!xE              %..................T
\!wlet\!ticklocationincr=\!yE          %..................T
\!wlet\!tickwidth=\!yS                 %..................T
\!wlet\!totalleaderlength=\!dimenE     %.........U.........
\!wlet\!xone=\!xprime                  %....X..............
\!wlet\!xtwo=\!dxprime                 %....X..............
\!wlet\!ySsave=\!yM                    %...................
\!wlet\!ybB=\!yB                       %.................S.
\!wlet\!ybC=\!yC                       %.................S.
\!wlet\!ybE=\!yE                       %.................S.
\!wlet\!ybM=\!yM                       %.................S.
\!wlet\!ybS=\!yS                       %.................S.
\!wlet\!ybpos=\!yyloc                  %.................S.
\!wlet\!yone=\!yprime                  %....X..............
\!wlet\!ytB=\!xB                       %.................S.
\!wlet\!ytC=\!xC                       %.................S.
\!wlet\!ytE=\!downlength               %.................S.
\!wlet\!ytM=\!arclength                %.................S.
\!wlet\!ytS=\!distacross               %.................S.
\!wlet\!ytpos=\!xxloc                  %.................S.
\!wlet\!ytwo=\!dyprime                 %....X..............

% Initial values for registers
\!zpt=0pt                              % static
\!xunit=1pt
\!yunit=1pt
\!arearloc=\!xunit
\!areatloc=\!yunit
\!dshade=5pt
\!leaderlength=24in
\!tfs=256                              % static
\!wmax=5.3pt                           % static
\!wmin=2.7pt                           % static
\!xaxislength=\!xunit
\!xpivot=\!zpt
\!yaxislength=\!yunit 
\!ypivot=\!zpt
\plotsymbolspacing=.4pt
  \!dimenA=50pt \!fiftypt=\!dimenA     % static

\!rootten=3.162278pt                   % static
\!tenAa=8.690286pt                     % static  (A5)
\!tenAc=2.773839pt                     % static  (A3)
\!tenAe=2.543275pt                     % static  (A1)

% Initial values for control sequences
\def\!cosrotationangle{1}      %................R..
\def\!sinrotationangle{0}      %................R..
\def\!xpivotcoord{0}           %................R..
\def\!xref{0}                  %............P......
\def\!xshadesave{0}            %.................S.
\def\!ypivotcoord{0}           %................R..
\def\!yref{0}                  %............P......
\def\!yshadesave{0}            %.................S.
\def\!zero{0}                  %..................T

% Reset TeX to report allocations
\let\wlog=\!!!wlog
%  *************************************
%  ***  AREAS: Deals with plot areas ***
%  *************************************
%
%  ** User commands
%  **   \setplotarea x from LEFT XCOORD to RIGTH XCOORD, y from BOTTOM YCOORD
%  **      to TOP YCOORD
%  **   \axis BOTTOM-LEFT-TOP-RIGHT  [SHIFTEDTO xy=COORD] [VISIBLE-INVISIBLE]
%  **      [LABEL {label}] [TICKS] /
%  **   \visibleaxes
%  **   \invisibleaxes
%  **   \plotheading {HEADING}
%  **   \grid {# of columns} {# of rows}
%  **   \normalgraphs 
  
%  **  \normalgraphs
%  **    Sets defaults for graph setup. See Subsection 3.4 of manual.
\def\normalgraphs{%
  \longticklength=.4\baselineskip
  \shortticklength=.25\baselineskip
  \tickstovaluesleading=.25\baselineskip
  \valuestolabelleading=.8\baselineskip
  \linethickness=.4pt
  \stackleading=.17\baselineskip
  \headingtoplotskip=1.5\baselineskip
  \visibleaxes
  \ticksout
  \nogridlines
  \unloggedticks}
%
% **  \setplotarea x from LEFT XCOORD to RIGTH XCOORD, y from BOTTOM YCOORD
% **    to TOP YCOORD
% **  Reserves space in PICBOX for a rectangular box with the indicated
% **   coordinates.  Must be specified before calls to  \axis, 
% **   \grid, \plotheading.
% **  See Subsection 3.1 of the manual.
\def\setplotarea x from #1 to #2, y from #3 to #4 {%
  \!arealloc=\!M{#1}\!xunit \advance \!arealloc -\!xorigin
  \!areabloc=\!M{#3}\!yunit \advance \!areabloc -\!yorigin
  \!arearloc=\!M{#2}\!xunit \advance \!arearloc -\!xorigin
  \!areatloc=\!M{#4}\!yunit \advance \!areatloc -\!yorigin
  \!initinboundscheck
  \!xaxislength=\!arearloc  \advance\!xaxislength -\!arealloc
  \!yaxislength=\!areatloc  \advance\!yaxislength -\!areabloc
  \!plotheadingoffset=\!zpt
  \!dimenput {{\setbox0=\hbox{}\wd0=\!xaxislength\ht0=\!yaxislength\box0}}
     [bl] (\!arealloc,\!areabloc)}
%
% ** \visibleaxes, \invisibleaxes 
% ** Switches for setting visibility of subsequent axes.
% ** See Subsection 3.2 of the manual.
\def\visibleaxes{%
  \def\!axisvisibility{\!axisvisibletrue}}

%
% ** The next few macros enable the user to fix up an erroneous keyword
% **   in the \axis command.
%  \newhelp is in ALLOCATIONS
%  \newhelp\!keywordhelp{The keyword mentioned in the error message in unknown. 
%  Replace NEW KEYWORD in the indicated response by the keyword that 
%  should have been specified.}

\def\!fixkeyword#1{%
  \errhelp=\!keywordhelp
  \errmessage{Unrecognized keyword `#1': \the\!keywordtoks{NEW KEYWORD}'}}

%  \newtoks\!keywordtoks    In ALLOCATIONS.
\!keywordtoks={enter `i\fixkeyword}

\def\fixkeyword#1{%
  \!nextkeyword#1 }

% ** \axis BOTTOM-LEFT-TOP-RIGHT  [SHIFTEDTO xy=COORD] [VISIBLE-INVISIBLE]
% **   [LABEL {label}] [TICKS] /
% ** Exactly one of the keywords BOTTOM, LEFT, TOP, RIGHT must be
% ** specified. Axis is drawn along the indicated edge of the current
% ** plot area, shifted if the SHIFTEDTO option is used, visible or
% ** invisible according the selected option, with an optional LABEL,
% ** and optional TICKS (see ticks.tex for the options avialabel with
% ** TICKS). The TICKS option must be the last one specified. The \axis
% ** MUST be terminated with a / followed by a space.
% ** See Subsection 3.2 of the manual for more information.

% ** The various options of the \axis command are processed by the
% ** \!nextkeyword macro defined below. For example, 
% ** `\!nextkeyword shiftedto ' expands to `\!axisshiftedto'.
\def\axis {%
  \def\!nextkeyword##1 {%
    \expandafter\ifx\csname !axis##1\endcsname \relax
      \def\!next{\!fixkeyword{##1}}%
    \else
      \def\!next{\csname !axis##1\endcsname}%
    \fi
    \!next}%
  \!offset=\!zpt
  \!axisvisibility
  \!placeaxislabelfalse
  \!nextkeyword}

% ** This and the various macros that follow handle the keyword
% ** specifications on the \axis command
% ** See Subsection 3.2 of the manual.
\def\!axisbottom{%
  \!axisylevel=\!areabloc
  \def\!tickxsign{0}%
  \def\!tickysign{-}%
  \def\!axissetup{\!axisxsetup}%
  \def\!axislabeltbrl{t}%
  \!nextkeyword}

\def\!axistop{%
  \!axisylevel=\!areatloc
  \def\!tickxsign{0}%
  \def\!tickysign{+}%
  \def\!axissetup{\!axisxsetup}%
  \def\!axislabeltbrl{b}%
  \!nextkeyword}

\def\!axisleft{%
  \!axisxlevel=\!arealloc
  \def\!tickxsign{-}%
  \def\!tickysign{0}%
  \def\!axissetup{\!axisysetup}%
  \def\!axislabeltbrl{r}%
  \!nextkeyword}

\def\!axisright{%
  \!axisxlevel=\!arearloc
  \def\!tickxsign{+}%
  \def\!tickysign{0}%
  \def\!axissetup{\!axisysetup}%
  \def\!axislabeltbrl{l}%
  \!nextkeyword}

\def\!axisshiftedto#1=#2 {%
  \if 0\!tickxsign
    \!axisylevel=\!M{#2}\!yunit
    \advance\!axisylevel -\!yorigin
  \else
    \!axisxlevel=\!M{#2}\!xunit
    \advance\!axisxlevel -\!xorigin
  \fi
  \!nextkeyword}

\def\!axisvisible{%
  \!axisvisibletrue  
  \!nextkeyword}

\def\!axisinvisible{%
  \!axisvisiblefalse
  \!nextkeyword}

\def\!axislabel#1 {%
  \!axisLaBeL={#1}%
  \!placeaxislabeltrue
  \!nextkeyword}

\expandafter\def\csname !axis/\endcsname{%
  \!axissetup % This could done already by "ticks"; if so, now \relax
  \if!placeaxislabel
    \!placeaxislabel
  \fi
  \if +\!tickysign %                 ** (A "top" axis)
    \!dimenA=\!axisylevel
    \advance\!dimenA \!offset %      ** dimA = top of the axis structure
    \advance\!dimenA -\!areatloc %   ** dimA = excess over the plot area
    \ifdim \!dimenA>\!plotheadingoffset
      \!plotheadingoffset=\!dimenA % ** Greatest excess over the plot area
    \fi
  \fi}

% ** \grid {c} {r} 
% ** Partitions the plot area into c columns and r rows; see Subsection 3.3
% ** of the manual.
% ** (Other grid patterns can be drawn with the TICKS option of the \axis 
% ** command.
\def\grid #1 #2 {%
  \!countA=#1\advance\!countA 1
  \axis bottom invisible ticks length <\!zpt> andacross quantity {\!countA} /
  \!countA=#2\advance\!countA 1
  \axis left   invisible ticks length <\!zpt> andacross quantity {\!countA} / }

% ** \plotheading{HEADING}
% ** Places HEADING centered above the top of the plotarea (and above
% ** any top axis ticks marks, tick labels, and axis label); see
% ** Subsection 3.3 of the manual.
\def\plotheading#1 {%
  \advance\!plotheadingoffset \headingtoplotskip
  \!dimenput {#1} [B] <.5\!xaxislength,\!plotheadingoffset>
    (\!arealloc,\!areatloc)}

% ** From here on, the routines are internal.
\def\!axisxsetup{%
  \!axisxlevel=\!arealloc
  \!axisstart=\!arealloc
  \!axisend=\!arearloc
  \!axisLength=\!xaxislength
  \!!origin=\!xorigin
  \!!unit=\!xunit
  \!xswitchtrue
  \if!axisvisible 
    \!makeaxis
  \fi}

\def\!axisysetup{%
  \!axisylevel=\!areabloc
  \!axisstart=\!areabloc
  \!axisend=\!areatloc
  \!axisLength=\!yaxislength
  \!!origin=\!yorigin
  \!!unit=\!yunit
  \!xswitchfalse
  \if!axisvisible
    \!makeaxis
  \fi}

\def\!makeaxis{%
  \setbox\!boxA=\hbox{% (Make a pseudo-y[x] tick for an x[y]-axis)
    \beginpicture
      \!setdimenmode
      \setcoordinatesystem point at {\!zpt} {\!zpt}   
      \putrule from {\!zpt} {\!zpt} to
        {\!tickysign\!tickysign\!axisLength} 
        {\!tickxsign\!tickxsign\!axisLength}
    \endpicturesave <\!Xsave,\!Ysave>}%
    \wd\!boxA=\!zpt
    \!placetick\!axisstart}

\def\!placeaxislabel{%
  \advance\!offset \valuestolabelleading
  \if!xswitch
    \!dimenput {\the\!axisLaBeL} [\!axislabeltbrl]
      <.5\!axisLength,\!tickysign\!offset> (\!axisxlevel,\!axisylevel)
    \advance\!offset \!dp  % ** advance offset by the "tallness"
    \advance\!offset \!ht  % ** of the label
  \else
    \!dimenput {\the\!axisLaBeL} [\!axislabeltbrl]
      <\!tickxsign\!offset,.5\!axisLength> (\!axisxlevel,\!axisylevel)
  \fi
  \!axisLaBeL={}}

% *******************************
% *** ARROWS  (Draws arrows)  ***
% *******************************
%
% ** User commands
% **  \arrow <ARROW HEAD LENGTH> [MID FRACTION, BASE FRACTION]
% **    [<XSHIFT,YSHIFT>] from XFROM YFROM to XTO YTO
% **  \betweenarrows {TEXT} [orientation & shift] from XFROM YFROM to XTO YTO

% ** \arrow <ARROW HEAD LENGTH> [MID FRACTION, BASE FRACTION]
% **    [<XSHIFT,YSHIFT>] from XFROM YFROM to XTO YTO
% ** Draws an arrow from (XFROM,YFROM) to (XTO,YTO).  The arrow head
% ** is constructed two quadratic arcs, which extend back a distance
% ** ARROW HEAD LENGTH (a dimension) on both sides of the arrow shaft.
% ** All the way back the arcs are a distance BASE FRACTION*ARROW HEAD
% ** LENGTH apart, while half-way back they are a distance MID FRACTION*
% ** ARROW HEAD LENGTH apart. <XSHIFT,YSHIFT> is optional, and has
% ** its usual interpreation. See Subsection 5.4 of the manual.

\def\arrow <#1> [#2,#3]{%
  \!ifnextchar<{\!arrow{#1}{#2}{#3}}{\!arrow{#1}{#2}{#3}<\!zpt,\!zpt> }}

\def\!arrow#1#2#3<#4,#5> from #6 #7 to #8 #9 {%
%
% ** convert to dimensions
  \!xloc=\!M{#8}\!xunit   
  \!yloc=\!M{#9}\!yunit
  \!dxpos=\!xloc  \!dimenA=\!M{#6}\!xunit  \advance \!dxpos -\!dimenA
  \!dypos=\!yloc  \!dimenA=\!M{#7}\!yunit  \advance \!dypos -\!dimenA
  \let\!MAH=\!M%                         ** save current c/d mode
  \!setdimenmode%                        ** go into dimension mode
  \!xshift=#4\relax  \!yshift=#5\relax%  ** pick up shift
  \!reverserotateonly\!xshift\!yshift%   ** back rotate shift
  \advance\!xshift\!xloc  \advance\!yshift\!yloc
%
% **  draw shaft of arrow
  \!xS=-\!dxpos  \advance\!xS\!xshift
  \!yS=-\!dypos  \advance\!yS\!yshift
  \!start (\!xS,\!yS)
  \!ljoin (\!xshift,\!yshift)
%
% ** find 32*cosine and 32*sine of angle of rotation
  \!Pythag\!dxpos\!dypos\!arclength
  \!divide\!dxpos\!arclength\!dxpos  
  \!dxpos=32\!dxpos  \!removept\!dxpos\!!cos
  \!divide\!dypos\!arclength\!dypos  
  \!dypos=32\!dypos  \!removept\!dypos\!!sin
% 
% ** construct arrowhead
  \!halfhead{#1}{#2}{#3}%                ** draw half of arrow head
  \!halfhead{#1}{-#2}{-#3}%              ** draw other half
  \let\!M=\!MAH%                         ** restore old c/d mode
  \ignorespaces}
%
% ** draw half of arrow head
  \def\!halfhead#1#2#3{%
    \!dimenC=-#1%                
    \divide \!dimenC 2 %                 ** half way back
    \!dimenD=#2\!dimenC%                 ** half the mid width
    \!rotate(\!dimenC,\!dimenD)by(\!!cos,\!!sin)to(\!xM,\!yM)
    \!dimenC=-#1%                        ** all the way back
    \!dimenD=#3\!dimenC
    \!dimenD=.5\!dimenD%                 ** half the full width
    \!rotate(\!dimenC,\!dimenD)by(\!!cos,\!!sin)to(\!xE,\!yE)
    \!start (\!xshift,\!yshift)
    \advance\!xM\!xshift  \advance\!yM\!yshift
    \advance\!xE\!xshift  \advance\!yE\!yshift
    \!qjoin (\!xM,\!yM) (\!xE,\!yE) 
    \ignorespaces}

% ** \betweenarrows {TEXT} [orientation & shift] from XFROM YFROM to XTO YTO
% **   Makes things like <--- text --->, using arrow heads from TeX's fonts.
% **   See Subsection 5.4 of the manual.
\def\betweenarrows #1#2 from #3 #4 to #5 #6 {%
  \!xloc=\!M{#3}\!xunit  \!xxloc=\!M{#5}\!xunit%   
  \!yloc=\!M{#4}\!yunit  \!yyloc=\!M{#6}\!yunit%           
  \!dxpos=\!xxloc  \advance\!dxpos by -\!xloc
  \!dypos=\!yyloc  \advance\!dypos by -\!yloc
  \advance\!xloc .5\!dxpos
  \advance\!yloc .5\!dypos
  \let\!MBA=\!M%           ** save current coord\dimen mode
  \!setdimenmode%          ** express locations in dimens
  \ifdim\!dypos=\!zpt
    \ifdim\!dxpos<\!zpt \!dxpos=-\!dxpos \fi
    \put {\!lrarrows{\!dxpos}{#1}}#2{} at {\!xloc} {\!yloc}
  \else
    \ifdim\!dxpos=\!zpt
      \ifdim\!dypos<\!zpt \!dypos=-\!dypos \fi
      \put {\!udarrows{\!dypos}{#1}}#2{} at {\!xloc} {\!yloc}
    \fi
  \fi
  \let\!M=\!MBA%           ** restore previous c/d mode
  \ignorespaces}

% ** Subroutine for left-right between arrows 
\def\!lrarrows#1#2{% #1=width, #2=text
  {\setbox\!boxA=\hbox{$\mkern-2mu\mathord-\mkern-2mu$}%
   \setbox\!boxB=\hbox{$\leftarrow$}\!dimenE=\ht\!boxB
   \setbox\!boxB=\hbox{}\ht\!boxB=2\!dimenE
   \hbox to #1{$\mathord\leftarrow\mkern-6mu
     \cleaders\copy\!boxA\hfil
     \mkern-6mu\mathord-$%
     \kern.4em $\vcenter{\box\!boxB}$$\vcenter{\hbox{#2}}$\kern.4em
     $\mathord-\mkern-6mu
     \cleaders\copy\!boxA\hfil
     \mkern-6mu\mathord\rightarrow$}}}

% ** Subroutine for up-down between arrows 
\def\!udarrows#1#2{% #1=width, #2=text
  {\setbox\!boxB=\hbox{#2}%
   \setbox\!boxA=\hbox to \wd\!boxB{\hss$\vert$\hss}%
   \!dimenE=\ht\!boxA \advance\!dimenE \dp\!boxA \divide\!dimenE 2
   \vbox to #1{\offinterlineskip
      \vskip .05556\!dimenE
      \hbox to \wd\!boxB{\hss$\mkern.4mu\uparrow$\hss}\vskip-\!dimenE
      \cleaders\copy\!boxA\vfil
      \vskip-\!dimenE\copy\!boxA
      \vskip\!dimenE\copy\!boxB\vskip.4em
      \copy\!boxA\vskip-\!dimenE
      \cleaders\copy\!boxA\vfil
      \vskip-\!dimenE \hbox to \wd\!boxB{\hss$\mkern.4mu\downarrow$\hss}
      \vskip .05556\!dimenE}}}

% ***************************
% *** BARS  (Draws bars)  ***
% ***************************
%
% ** User commands:
% ** \putbar [<XSHIFT,YSHIFT>] breadth <BREADTH> from XSTART YSTART
% **   to XEND YEND
% ** \setbars [<XSHIFT,YSHIFT>] breadth <BREADTH> baseline at XY = COORD
% **   [baselabels ([B_ORIENTATION_x,B_ORIENTATION_y] <B_XSHIFT,B_YSHIFT>)]
% **   [endlabels  ([E_ORIENTATION_x,E_ORIENTATION_y] <E_XSHIFT,E_YSHIFT>)]

% ** \putbar [<XSHIFT,YSHIFT>] breadth <BREADTH> from XSTART YSTART
% **   to XEND YEND
% ** Either XSTART=XEND or YSTART=YEND. Draws a rectangle between
% **   (XSTART,YSTART) & (XEND,YEND). The "depth" of the rectangle
% **   is determined by those two plot positions; its other
% **   dimension "breadth" is specified by the dimension BREADTH.
% ** See Subsection 4.2 of the manual.
\def\putbar#1breadth <#2> from #3 #4 to #5 #6 {%
  \!xloc=\!M{#3}\!xunit  \!xxloc=\!M{#5}\!xunit%   
  \!yloc=\!M{#4}\!yunit  \!yyloc=\!M{#6}\!yunit%           
  \!dypos=\!yyloc  \advance\!dypos by -\!yloc
  \!dimenI=#2  
  \ifdim \!dimenI=\!zpt %            ** If 0 breadth
    \putrule#1from {#3} {#4} to {#5} {#6} % ** Then draw line
  \else %                            ** Else, put in a rectangle
    \let\!MBar=\!M%                  ** save current c/d mode
    \!setdimenmode %                 ** go into dimension mode
    \divide\!dimenI 2
    \ifdim \!dypos=\!zpt             
      \advance \!yloc -\!dimenI %    ** Equal y coordinates
      \advance \!yyloc \!dimenI
    \else
      \advance \!xloc -\!dimenI %    ** Equal x coordinates
      \advance \!xxloc \!dimenI
    \fi
    \putrectangle#1corners at {\!xloc} {\!yloc} and {\!xxloc} {\!yyloc}
    \let\!M=\!MBar %                 ** restore c/d mode
  \fi
  \ignorespaces}

% ** \setbars [<XSHIFT,YSHIFT>] breadth <BREADTH> baseline at XY = COORD
% **   [baselabels ([B_ORIENTATION_x,B_ORIENTATION_y] <B_XSHIFT,B_YSHIFT>)]
% **   [endlabels  ([E_ORIENTATION_x,E_ORIENTATION_y] <E_XSHIFT,E_YSHIFT>)]
% ** This command puts PiCTeX into the bar graph drawing mode described
% **   in Subsection 4.4 of the manual.
\def\setbars#1breadth <#2> baseline at #3 = #4 {%
  \edef\!barshift{#1}%
  \edef\!barbreadth{#2}%
  \edef\!barorientation{#3}%
  \edef\!barbaseline{#4}%
  \def\!bardobaselabel{\!bardoendlabel}%
  \def\!bardoendlabel{\!barfinish}%
  \let\!drawcurve=\!barcurve
  \!setbars}
\def\!setbars{%
  \futurelet\!nextchar\!!setbars}
\def\!!setbars{%
  \if b\!nextchar
    \def\!!!setbars{\!setbarsbget}%
  \else 
    \if e\!nextchar
      \def\!!!setbars{\!setbarseget}%
    \else
      \def\!!!setbars{\relax}%
    \fi
  \fi
  \!!!setbars}
\def\!setbarsbget baselabels (#1) {%
  \def\!barbaselabelorientation{#1}%
  \def\!bardobaselabel{\!!bardobaselabel}%
  \!setbars}
\def\!setbarseget endlabels (#1) {%
  \edef\!barendlabelorientation{#1}%
  \def\!bardoendlabel{\!!bardoendlabel}%
  \!setbars}

% ** \!barcurve
% ** Draws a bargraph with preset values of barshift, barbreadth,
% ** barorientation (x or y) and barbaseline (coordinate)
\def\!barcurve #1 #2 {%
  \if y\!barorientation
    \def\!basexarg{#1}%
    \def\!baseyarg{\!barbaseline}%
  \else
    \def\!basexarg{\!barbaseline}%
    \def\!baseyarg{#2}%
  \fi
  \expandafter\putbar\!barshift breadth <\!barbreadth> from {\!basexarg}
    {\!baseyarg} to {#1} {#2}
  \def\!endxarg{#1}%
  \def\!endyarg{#2}%
  \!bardobaselabel}

\def\!!bardobaselabel "#1" {%
  \put {#1}\!barbaselabelorientation{} at {\!basexarg} {\!baseyarg}
  \!bardoendlabel}
 
\def\!!bardoendlabel "#1" {%
  \put {#1}\!barendlabelorientation{} at {\!endxarg} {\!endyarg}
  \!barfinish}

\def\!barfinish{%
  \!ifnextchar/{\!finish}{\!barcurve}}

% ********************************
% *** BOXES (Draws rectangles) ***
% ********************************
%
% ** User commands:
% **   \putrectangle [<XSHIFT,YSHIFT>] corners at  XCOORD1 YCOORD1
% **     and  XCOORD2 YCOORD2 
% **   \shaderectangleson
% **   \shaderectanglesoff
% **   \frame [<SEPARATION>] {TEXT}
% **   \rectangle <WIDTH> <HEIGHT>
%
%
% **  \putrectangle [<XSHIFT,YSHIFT>] corners at XCOORD1 YCOORD1 
% **    and  XCOORD2 YCOORD2 
% **  Draws a rectangle with corners at (X1,Y1), (X2,Y1), (X1,Y2), (X2,Y2)
% **  Lines have thickness \linethickness, and overlap at the corners.
% **  The optional field  <XSHIFT,YSHIFT>  functions as with a \put command.
% **  See Subsection 4.2 of the manual.
\def\putrectangle{%
  \!ifnextchar<{\!putrectangle}{\!putrectangle<\!zpt,\!zpt> }}
\def\!putrectangle<#1,#2> corners at #3 #4 and #5 #6 {%
%
% ** get locations
  \!xone=\!M{#3}\!xunit  \!xtwo=\!M{#5}\!xunit%   
  \!yone=\!M{#4}\!yunit  \!ytwo=\!M{#6}\!yunit%           
  \ifdim \!xtwo<\!xone
    \!dimenI=\!xone  \!xone=\!xtwo  \!xtwo=\!dimenI
  \fi
  \ifdim \!ytwo<\!yone
    \!dimenI=\!yone  \!yone=\!ytwo  \!ytwo=\!dimenI
  \fi
  \!dimenI=#1\relax  \advance\!xone\!dimenI  \advance\!xtwo\!dimenI
  \!dimenI=#2\relax  \advance\!yone\!dimenI  \advance\!ytwo\!dimenI
  \let\!MRect=\!M%                  ** save current coord/dimen mode
  \!setdimenmode
%
% ** shade rectangle if appropriate
  \!shaderectangle
%
% ** draw horizontal edges
  \!dimenI=.5\linethickness
  \advance \!xone  -\!dimenI%       ** adjust x-location to overlap corners
  \advance \!xtwo   \!dimenI%       ** ditto
  \putrule from {\!xone} {\!yone} to {\!xtwo} {\!yone} 
  \putrule from {\!xone} {\!ytwo} to {\!xtwo} {\!ytwo} 
%
% ** draw vertical edges
  \advance \!xone   \!dimenI%       ** restore original x-values
  \advance \!xtwo  -\!dimenI% 
  \advance \!yone  -\!dimenI%       ** adjust y-location to overlap corners
  \advance \!ytwo   \!dimenI%       ** ditto
  \putrule from {\!xone} {\!yone} to {\!xone} {\!ytwo} 
  \putrule from {\!xtwo} {\!yone} to {\!xtwo} {\!ytwo} 
  \let\!M=\!MRect%                  ** restore coord/dimen mode
  \ignorespaces}
 
% ** \shaderectangleson 
% **   Subsequent rectangles will be shaded according to 
% **   the current shading pattern.  Affects \putrectangle, \putbar,
% **   \frame, \sethistograms, and \setbars. See Subsection 7.5 of the manual.

% ** \shaderectanglesoff 
% **    Suppresses  \shaderectangleson.  The default.
\def\shaderectanglesoff{%
  \def\!shaderectangle{}%
  \ignorespaces}

\shaderectanglesoff
 
% ** The following internal routine shades the current rectangle, when
% **   \!shaderectangle = \!!shaderectangle . 
\def\!!shaderectangle{%
  \!dimenA=\!xtwo  \advance \!dimenA -\!xone
  \!dimenB=\!ytwo  \advance \!dimenB -\!yone
  \ifdim \!dimenA<\!dimenB
    \!startvshade (\!xone,\!yone,\!ytwo)
    \!lshade      (\!xtwo,\!yone,\!ytwo)
  \else
    \!starthshade (\!yone,\!xone,\!xtwo)
    \!lshade      (\!ytwo,\!xone,\!xtwo)
  \fi
  \ignorespaces}
  
% ** \frame [<SEPARATION>] {TEXT}
% ** Draws a frame of thickness linethickness about the box enclosing
% **   TEXT; the frame is separated from the box by a distance of
% **   SEPARATION.  The result is an hbox with the same baseline as TEXT.
% **   If <SEPARATION> is omitted, you get the effect of <0pt>.
% ** See Subsection 4.2 of the manual.
\def\frame{%
  \!ifnextchar<{\!frame}{\!frame<\!zpt> }}
\long\def\!frame<#1> #2{%
  \beginpicture
    \setcoordinatesystem units <1pt,1pt> point at 0 0 
    \put {#2} [Bl] at 0 0 
    \!dimenA=#1\relax
    \!dimenB=\!wd \advance \!dimenB \!dimenA
    \!dimenC=\!ht \advance \!dimenC \!dimenA
    \!dimenD=\!dp \advance \!dimenD \!dimenA
    \let\!MFr=\!M
    \!setdimenmode
    \putrectangle corners at {-\!dimenA} {-\!dimenD} and {\!dimenB} {\!dimenC}
    \!setcoordmode
    \let\!M=\!MFr
  \endpicture
  \ignorespaces}
 
% ** \rectangle <WIDTH> <HEIGHT>
% ** Constructs a rectangle of width WIDTH and heigth HEIGHT. 
% ** See Subsection 4.2 of the manual.
\def\rectangle <#1> <#2> {%
  \setbox0=\hbox{}\wd0=#1\ht0=#2\frame {\box0}}

% *********************************************
% ***  CURVES  (Upper level \plot commands) ***
% *********************************************
%
% ** User commands
% **   \plot  DATA  /
% **   \plot  "FILE NAME"
% **   \setquadratic
% **   \setlinear
% **   \sethistograms
% **   \vshade  ...
% **   \hshade  ...

% \plot: multi-purpose command. Draws histograms, bar graphs, piecewise-linear
% or piecewise quadratic curves, depending on the setting of \!drawcurve.
% See Subsections 4.3-4.5, 5.1, 5.2 of the manual.
\def\plot{%
  \!ifnextchar"{\!plotfromfile}{\!drawcurve}}
\def\!plotfromfile"#1"{%
  \expandafter\!drawcurve \input #1 /}

% Command to set piecewise quadratic mode
% See Subsections 5.1, 7.3, and 7.4 of the manual.
\def\setquadratic{%
  \let\!drawcurve=\!qcurve
  \let\!!Shade=\!!qShade
  \let\!!!Shade=\!!!qShade}

% Command to set piecewise linear mode
% See Subsections 5.1, 7.3, and 7.4 of the manual.
\def\setlinear{%
  \let\!drawcurve=\!lcurve
  \let\!!Shade=\!!lShade
  \let\!!!Shade=\!!!lShade}

% Command to set histogram mode
% See Subsection 4.3 of the manual.
\def\sethistograms{%
  \let\!drawcurve=\!hcurve}

% Commands to cycle through list of coordinates in piecewise quadratic 
% interpolation mode
\def\!qcurve #1 #2 {%
  \!start (#1,#2)
  \!Qjoin}
\def\!Qjoin#1 #2 #3 #4 {%
  \!qjoin (#1,#2) (#3,#4)             % \!qjoin  is defined in QUADRATIC
  \!ifnextchar/{\!finish}{\!Qjoin}}

% Commands to cycle through list of coordinates in piecewise linear 
% interpolation mode
\def\!lcurve #1 #2 {%
  \!start (#1,#2)
  \!Ljoin}
\def\!Ljoin#1 #2 {%
  \!ljoin (#1,#2)                    % \!ljoin  is defined in LINEAR
  \!ifnextchar/{\!finish}{\!Ljoin}}

\def\!finish/{\ignorespaces}

% Command to cycle through list of coordinates in histogram mode
\def\!hcurve #1 #2 {%
  \edef\!hxS{#1}%
  \edef\!hyS{#2}%
  \!hjoin}
\def\!hjoin#1 #2 {%
  \putrectangle corners at {\!hxS} {\!hyS} and {#1} {#2}
  \edef\!hxS{#1}%
  \!ifnextchar/{\!finish}{\!hjoin}}

% \vshade: See Subsection 7.3 of the manual.
\def\vshade #1 #2 #3 {%
  \!startvshade (#1,#2,#3)
  \!Shadewhat}

% \hshade: See Subsection 7.4 of the manual.
\def\hshade #1 #2 #3 {%
  \!starthshade (#1,#2,#3)
  \!Shadewhat}

% Commands to cycle through coordinates and optional "edge effect"
% fields while shading.
\def\!Shadewhat{%
  \futurelet\!nextchar\!Shade}
\def\!Shade{%
  \if <\!nextchar
    \def\!nextShade{\!!Shade}%
  \else
    \if /\!nextchar
      \def\!nextShade{\!finish}%
    \else
      \def\!nextShade{\!!!Shade}%
    \fi
  \fi
  \!nextShade}
\def\!!lShade<#1> #2 #3 #4 {%
  \!lshade <#1> (#2,#3,#4)                 % \!lshade is defined in SHADING
  \!Shadewhat}
\def\!!!lShade#1 #2 #3 {%
  \!lshade (#1,#2,#3)
  \!Shadewhat} 
\def\!!qShade<#1> #2 #3 #4 #5 #6 #7 {%
  \!qshade <#1> (#2,#3,#4) (#5,#6,#7)      % \!qshade is defined in SHADING
  \!Shadewhat}
\def\!!!qShade#1 #2 #3 #4 #5 #6 {%
  \!qshade (#1,#2,#3) (#4,#5,#6)
  \!Shadewhat} 

% ** Set default interpolation mode
\setlinear

\def\setdashpattern <#1>{%
  \def\!Flist{}\def\!Blist{}\def\!UDlist{}%
  \!countA=0
  \!ecfor\!item:=#1\do{%
    \!dimenA=\!item\relax
    \expandafter\!rightappend\the\!dimenA\withCS{\\}\to\!UDlist%
    \advance\!countA  1
    \ifodd\!countA
      \expandafter\!rightappend\the\!dimenA\withCS{\!Rule}\to\!Flist%
      \expandafter\!leftappend\the\!dimenA\withCS{\!Rule}\to\!Blist%
    \else 
      \expandafter\!rightappend\the\!dimenA\withCS{\!Skip}\to\!Flist%
      \expandafter\!leftappend\the\!dimenA\withCS{\!Skip}\to\!Blist%
    \fi}%
  \!leaderlength=\!zpt
  \def\!Rule##1{\advance\!leaderlength  ##1}%
  \def\!Skip##1{\advance\!leaderlength  ##1}%
  \!Flist%
  \ifdim\!leaderlength>\!zpt 
  \else
    \def\!Flist{\!Skip{24in}}\def\!Blist{\!Skip{24in}}\ignorespaces
    \def\!UDlist{\\{\!zpt}\\{24in}}\ignorespaces
    \!leaderlength=24in
  \fi
  \!dashingon}

%  **  \!dashingon  -- puts the curve drawing routines into dash mode
%  **  \!dashingoff -- puts the curve drawing routines into solid mode
%  **  These are internal commands, invoked by \setdashpattern and \setsolid
\def\!dashingon{%
  \def\!advancedashing{\!!advancedashing}%
  \def\!drawlinearsegment{\!lineardashed}%
  \def\!puthline{\!putdashedhline}%
  \def\!putvline{\!putdashedvline}%
  \ignorespaces}% 
\def\!dashingoff{%
  \def\!advancedashing{\relax}%
  \def\!drawlinearsegment{\!linearsolid}%
  \def\!puthline{\!putsolidhline}%
  \def\!putvline{\!putsolidvline}%
  \ignorespaces}

%  **  \setdots <LENGTH>  --  sets up a dot/skip pattern where dot (actually
%  **    the current plotsymbol) is plunked down once for every LENGTH 
%  **    traveled along the curve.  LENGTH defaults to 5pt.
%  **    See Subsection 6.1 of the manual.
\def\setdots{%
  \!ifnextchar<{\!setdots}{\!setdots<5pt>}}
\def\!setdots<#1>{%
  \!dimenB=#1\advance\!dimenB -\plotsymbolspacing
  \ifdim\!dimenB<\!zpt
    \!dimenB=\!zpt
  \fi
\setdashpattern <\plotsymbolspacing,\!dimenB>}
 
% ** \setdotsnear <LENGTH> for <ARC LENGTH>
% ** sets up a dot pattern where the dots are approximately LENGTH apart,
% ** the total length of the pattern is ARC LENGTH, and the pattern
% ** begins and ends with a dot. See Subsection 6.3 of the manual.
\def\setdotsnear <#1> for <#2>{%
  \!dimenB=#2\relax  \advance\!dimenB -.05pt  
  \!dimenC=#1\relax  \!countA=\!dimenC 
  \!dimenD=\!dimenB  \advance\!dimenD .5\!dimenC  \!countB=\!dimenD
  \divide \!countB  \!countA
  \ifnum 1>\!countB 
    \!countB=1
  \fi
  \divide\!dimenB  \!countB
  \setdots <\!dimenB>}
 
%  **  \setdashes <LENGTH>  --  sets up a dash/skip pattern where the dash
%  **    and the skip are each of length LENGTH (the dash is formed by
%  **    plunking down the current plotsymbol over an arc of length LENGTH
%  **    and so may actually be longer than LENGTH.  LENGTH defaults to 5pt.
%  **    See Subsection 6.1 of the manual.
\def\setdashes{%
  \!ifnextchar<{\!setdashes}{\!setdashes<5pt>}}
\def\!setdashes<#1>{\setdashpattern <#1,#1>}
 
% ** \setdashesnear ...
% ** Like \setdotsnear; the pattern begins and ends with a dash.
% ** See Subsection 6.3 of the manual.
\def\setdashesnear <#1> for <#2>{%
  \!dimenB=#2\relax  
  \!dimenC=#1\relax  \!countA=\!dimenC 
  \!dimenD=\!dimenB  \advance\!dimenD .5\!dimenC  \!countB=\!dimenD
  \divide \!countB  \!countA
  \ifodd \!countB 
  \else 
    \advance \!countB  1
  \fi
  \divide\!dimenB  \!countB
  \setdashes <\!dimenB>}
 
%  **  \setsolid  --  puts the curve drawing routines in "solid line" mode,
%  **    the default mode.  See Subsection 6.1 of the manual.
\def\setsolid{%
  \def\!Flist{\!Rule{24in}}\def\!Blist{\!Rule{24in}}%  
  \def\!UDlist{\\{24in}\\{\!zpt}}%
  \!dashingoff}  
\setsolid

%  **  \findlength {CURVE CMDS}
%  **  PiCTeX executes the \start, \ljoin, and \qjoin cmds comprising
%  **  CURVE CMDS without plotting anything, but stashes the length
%  **  of the phantom curve away in \totalarclength.
%  **  See Subsection 6.3 of the manual.

% *************************************************************
% *** DIVISION  (Does long division of dimension registers) ***
% *************************************************************
 
% ** User command:
% **   \Divide {DIVIDEND} by {DIVISOR} forming {RESULT}
  
% ** Internal command
% **   \!divide{DIVIDEND}{DIVISOR}{RESULT}
 
% **  \!divide DIVIDEND [by] DIVISOR [to get] ANSWER
% **  Divides the dimension DIVIDEND by the dimension DIVISOR, placing the 
% **  quotient in the dimension register ANSWER.  Values are understood to 
% **  be in points.  E.g.  12.5pt/1.4pt=8.92857pt.
% **  Quotient is accurate to 1/65536pt=2**[-16]pt
% **  |DIVISOR| should be < 2048pt (about 28 inches).
\def\!divide#1#2#3{%
  \!dimenB=#1%                      **  dimB  holds current remainder (r)
  \!dimenC=#2%                      **  dimC  holds divisor (d)
  \!dimenD=\!dimenB%                **  dimD  holds quotient q=r/d for this 
  \divide \!dimenD \!dimenC%        **    step, in units of scaled pts
  \!dimenA=\!dimenD%                **  dimA  eventually holds answer (a)
  \multiply\!dimenD \!dimenC%       **  r <-- r - dq
  \advance\!dimenB -\!dimenD%       **  First step complete. Have integer part
%                                   **  of a, and corresponding remainder.
  \!dimenD=\!dimenC%                **  Temporarily use dimD to hold |d|
    \ifdim\!dimenD<\!zpt \!dimenD=-\!dimenD 
  \fi
  \ifdim\!dimenD<64pt%              **  Branch on the magnitude of |d|
    \!divstep[\!tfs]\!divstep[\!tfs]%
  \else 
    \!!divide
  \fi
  #3=\!dimenA\ignorespaces}

% **  The following code handles divisors  d  with 
% **    (1)  .88in =  64pt <= d <  256pt =  3.54in
% **    (2) 3.54in = 256pt <= d < 2048pt = 28.34in
% **  Anything bigger than that may result in an overflow condition.
% **  For our purposes, we should never even see case (2).
\def\!!divide{%
  \ifdim\!dimenD<256pt
    \!divstep[64]\!divstep[32]\!divstep[32]%
  \else 
    \!divstep[8]\!divstep[8]\!divstep[8]\!divstep[8]\!divstep[8]%
    \!dimenA=2\!dimenA
  \fi}

% **  The following macro does the real long division work.
\def\!divstep[#1]{%                 **  #1 = "B"
  \!dimenB=#1\!dimenB%              **  r <-- B*r
  \!dimenD=\!dimenB%                **  dimD  holds quotient q=r/d for this 
    \divide \!dimenD by \!dimenC%   **    step, in units of scaled pts
  \!dimenA=#1\!dimenA%              **  a <-- B*a + q
    \advance\!dimenA by \!dimenD%
  \multiply\!dimenD by \!dimenC%    **  r <-- r - dq
    \advance\!dimenB by -\!dimenD}
 
% **  \Divide:  See Subsection 9.3 of the manual.
\def\Divide <#1> by <#2> forming <#3> {%
  \!divide{#1}{#2}{#3}}

% *********************************************
% *** ELLIPSES (Draws ellipses and circles) ***
% *********************************************
 
% ** User commands
% **   \ellipticalarc  axes ratio A:B  DEGREES degrees from XSTART YSTART 
% **      center at XCENTER YCENTER 
% **   \circulararc DEGREES degrees from XSTART YSTART 
% **      center at XCENTER YCENTER 
 
% ** Internal command
% **   \!sinandcos{32*ANGLE in radians}{32*SIN}{32*COS}

% **   \ellipticalarc  axes ratio A:B  DEGREES degrees from XSTART YSTART 
% **      center at XCENTER YCENTER 
% **    Draws a elliptical arc starting at the coordinate point (XSTART,YSTART).
% **    The center of the ellipse of which the arc is a segment is at 
% **      (XCENTER,YCENTER).
% **    The arc extends through an angle of DEGREES degrees (may be + or -).
% **    A:B is the ratio of the length of the xaxis to the length of
% **      the yaxis of the ellipse
% **    Sqrt{[(XSTART-XCENTER)/A]**2 + [(YSTART-YCENTER)/B]**2}
% **      must be < 512pt (about 7in).
% **    Doesn't modify the dimensions (ht, dp, wd) of the PiCture under
% **      construction.
 
% ** \circulararc  --  See Subsection 5.3 of the manual.

% ** \ellipticalarc  --  See Subsection 5.3 of the manual.
\def\ellipticalarc axes ratio #1:#2 #3 degrees from #4 #5 center at #6 #7 {%
  \!angle=#3pt\relax%                    ** get angle
  \ifdim\!angle>\!zpt 
    \def\!sign{}%                        ** counterclockwise
  \else 
    \def\!sign{-}\!angle=-\!angle%       ** clockwise
  \fi
  \!xxloc=\!M{#6}\!xunit%                ** convert CENTER to dimension
  \!yyloc=\!M{#7}\!yunit     
  \!xxS=\!M{#4}\!xunit%                  ** get STARTing point on rim of ellipse
  \!yyS=\!M{#5}\!yunit
  \advance\!xxS -\!xxloc%                ** make center of ellipse (0,0)
  \advance\!yyS -\!yyloc
  \!divide\!xxS{#1pt}\!xxS %             ** scale point on ellipse to point on 
  \!divide\!yyS{#2pt}\!yyS %                 corresponding circle
  \let\!MC=\!M%                          ** save current c/d mode
  \!setdimenmode%                        ** go into dimension mode
  \!xS=#1\!xxS  \advance\!xS\!xxloc
  \!yS=#2\!yyS  \advance\!yS\!yyloc
  \!start (\!xS,\!yS)%
  \!loop\ifdim\!angle>14.9999pt%         ** draw in major portion of ellipse 
    \!rotate(\!xxS,\!yyS)by(\!cos,\!sign\!sin)to(\!xxM,\!yyM) 
    \!rotate(\!xxM,\!yyM)by(\!cos,\!sign\!sin)to(\!xxE,\!yyE)
    \!xM=#1\!xxM  \advance\!xM\!xxloc  \!yM=#2\!yyM  \advance\!yM\!yyloc
    \!xE=#1\!xxE  \advance\!xE\!xxloc  \!yE=#2\!yyE  \advance\!yE\!yyloc
    \!qjoin (\!xM,\!yM) (\!xE,\!yE)
    \!xxS=\!xxE  \!yyS=\!yyE 
    \advance \!angle -15pt
  \repeat
  \ifdim\!angle>\!zpt%                   ** complete remaining arc, if any
    \!angle=100.53096\!angle%            ** convert angle to radians, divide
    \divide \!angle 360 %                **   by 2, and multiply by 32
    \!sinandcos\!angle\!!sin\!!cos%      ** get 32*sin & 32*cos
    \!rotate(\!xxS,\!yyS)by(\!!cos,\!sign\!!sin)to(\!xxM,\!yyM) 
    \!rotate(\!xxM,\!yyM)by(\!!cos,\!sign\!!sin)to(\!xxE,\!yyE)
    \!xM=#1\!xxM  \advance\!xM\!xxloc  \!yM=#2\!yyM  \advance\!yM\!yyloc
    \!xE=#1\!xxE  \advance\!xE\!xxloc  \!yE=#2\!yyE  \advance\!yE\!yyloc
    \!qjoin (\!xM,\!yM) (\!xE,\!yE)
  \fi
  \let\!M=\!MC%                          ** restore c/d mode
  \ignorespaces}%                        **   if appropriate

%  ** \!rotate(XREG,YREG)by(32cos,32sin)to(XXREG,YYREG)
%  ** rotates (XREG,YREG) by angle with specfied scaled cos & sin to
%  ** (XXREG,YYREG).  Uses \!dimenA & \!dimenB as scratch registers.
\def\!rotate(#1,#2)by(#3,#4)to(#5,#6){% 
  \!dimenA=#3#1\advance \!dimenA -#4#2%   ** Rcos(x+t)=Rcosx*cost - Rsinx*sint
  \!dimenB=#3#2\advance \!dimenB  #4#1%   ** Rsin(x+t)=Rsinx*cost + Rcosx*sint
  \divide \!dimenA 32  \divide \!dimenB 32 
  #5=\!dimenA  #6=\!dimenB
  \ignorespaces}
\def\!sin{4.17684}%                       ** 32*sin(pi/24) (pi/24=7.5deg)
\def\!cos{31.72624}%                      ** 32*cos(pi/24)

%  ** \!sinandcos{32*ANGLE in radians}{\SINCS}{\COSCS}
%  **   Computes the 32*sine and 32*cosine of a small ANGLE expressed in 
%  **   radians/32 and puts these values in the replacement texts of 
%  **   \SINCS and \COSCS
\def\!sinandcos#1#2#3{%
 \!dimenD=#1%                **  angle is expressed in radians/32: 1pt = 1/32rad
 \!dimenA=\!dimenD%          **  dimA will eventually contain 32sin(angle)in pts
 \!dimenB=32pt%              **  dimB will eventually contain 32cos(angle)in pts
 \!removept\!dimenD\!value%  **  get value of 32*angle, without "pt"
 \!dimenC=\!dimenD%          **  holds 32*angle**i/i! in pts
 \!dimenC=\!value\!dimenC \divide\!dimenC by 64 %   ** now 32*angle**2/2
 \advance\!dimenB by -\!dimenC%                     ** 32-32*angle**2/2
 \!dimenC=\!value\!dimenC \divide\!dimenC by 96 %   ** now 32*angle**3/3!
 \advance\!dimenA by -\!dimenC%                     ** now 32*(angle-angle**3/6)
 \!dimenC=\!value\!dimenC \divide\!dimenC by 128 %  ** now 32*angle**4/4!
 \advance\!dimenB by \!dimenC%
 \!removept\!dimenA#2%                              ** set 32*sin(angle)
 \!removept\!dimenB#3%                              ** set 32*cos(angle)
 \ignorespaces}

% *****************************************************************
% ***  RULES  (Draws rules, i.e., horizontal & vertical lines)  ***
% *****************************************************************

% **  User command:
% **    \putrule [<XDIMEN,YDIMEN>] from  XCOORD1 YCOORD1 
% **      to  XCOORD2 YCOORD2 

% **  Internal commands:
% **    \!puthline [<XDIMEN,YDIMEN>]    (h = horizontal)
% **      Set by dashpat to either: \!putsolidhline  or \!putdashedhline
% **    \!putvline [<XDIMEN,YDIMEN>]    (v = vertical)
% **      Either:  \!putsolidvline  or  \!putdashedvline

% **  \putrule [<XDIMEN,YDIMEN>] from XCOORD1 YCOORD1
% **    to XCOORD2 YCOORD2
% **  Draws a rule -- dashed or solid depending on the current dash pattern --
% **    from (X1,Y1) to (X2,Y2).  Uses TEK's  \hrule & \vrule & \leaders  
% **    constructions to handle horizontal & vertical lines efficiently both
% **    in terms of execution time and space in the DVI file.  
% **  See Subsection 4.1 of the manual.
\def\putrule#1from #2 #3 to #4 #5 {%
  \!xloc=\!M{#2}\!xunit  \!xxloc=\!M{#4}\!xunit%   
  \!yloc=\!M{#3}\!yunit  \!yyloc=\!M{#5}\!yunit%           
  \!dxpos=\!xxloc  \advance\!dxpos by -\!xloc
  \!dypos=\!yyloc  \advance\!dypos by -\!yloc
  \ifdim\!dypos=\!zpt
    \def\!!Line{\!puthline{#1}}\ignorespaces
  \else
    \ifdim\!dxpos=\!zpt
      \def\!!Line{\!putvline{#1}}\ignorespaces
    \else 
       \def\!!Line{}
    \fi
  \fi
  \let\!ML=\!M%           ** save current coord\dimen mode
  \!setdimenmode%         ** express locations in dimens
  \!!Line%
  \let\!M=\!ML%           ** restore previous c/d mode
  \ignorespaces}

% **  \!putsolidhline [<XDIMEN,YDIMEN>]
% **  Place horizontal solid line
\def\!putsolidhline#1{%
  \ifdim\!dxpos>\!zpt 
    \put{\!hline\!dxpos}#1[l] at {\!xloc} {\!yloc}
  \else 
    \put{\!hline{-\!dxpos}}#1[l] at {\!xxloc} {\!yyloc}
  \fi
  \ignorespaces}
 
% **  \!putsolidvline [shifted <XDIMEN,YDIMEN>]
% **  Place vertical solid line
\def\!putsolidvline#1{%
  \ifdim\!dypos>\!zpt 
    \put{\!vline\!dypos}#1[b] at {\!xloc} {\!yloc}
  \else 
    \put{\!vline{-\!dypos}}#1[b] at {\!xxloc} {\!yyloc}
  \fi
  \ignorespaces}
 
\def\!hline#1{\hbox to #1{\leaders \hrule height\linethickness\hfill}}
\def\!vline#1{\vbox to #1{\leaders \vrule width\linethickness\vfill}}

% **  \!putdashedhline [<XDIMEN,YDIMEN>]
% **  Place dashed horizontal line
\def\!putdashedhline#1{%
  \ifdim\!dxpos>\!zpt 
    \!DLsetup\!Flist\!dxpos
    \put{\hbox to \!totalleaderlength{\!hleaders}\!hpartialpattern\!Rtrunc}
      #1[l] at {\!xloc} {\!yloc} 
  \else 
    \!DLsetup\!Blist{-\!dxpos}
    \put{\!hpartialpattern\!Ltrunc\hbox to \!totalleaderlength{\!hleaders}}
      #1[r] at {\!xloc} {\!yloc} 
  \fi
  \ignorespaces}
 
% **  \!putdashedhline [<XDIMEN,YDIMEN>]
% **  Place dashed vertical line
\def\!putdashedvline#1{%
  \!dypos=-\!dypos%            ** vertical leaders go from top to bottom
  \ifdim\!dypos>\!zpt 
    \!DLsetup\!Flist\!dypos 
    \put{\vbox{\vbox to \!totalleaderlength{\!vleaders}
      \!vpartialpattern\!Rtrunc}}#1[t] at {\!xloc} {\!yloc} 
  \else 
    \!DLsetup\!Blist{-\!dypos}
    \put{\vbox{\!vpartialpattern\!Ltrunc
      \vbox to \!totalleaderlength{\!vleaders}}}#1[b] at {\!xloc} {\!yloc} 
  \fi
  \ignorespaces}

% **  The rest of the macros in this section are subroutines used by 
% **  \!putdashedhline and \!putdashedvline.
\def\!DLsetup#1#2{%            ** Dashed-Line set up
  \let\!RSlist=#1%             ** set !Rule-Skip list
  \!countB=#2%                 ** convert rule length to integer (number of sps)
  \!countA=\!leaderlength%     ** ditto, leaderlength
  \divide\!countB by \!countA% ** number of complete leader units
  \!totalleaderlength=\!countB\!leaderlength
  \!Rresiduallength=#2%
  \advance \!Rresiduallength by -\!totalleaderlength%  \** excess length
  \!Lresiduallength=\!leaderlength
  \advance \!Lresiduallength by -\!Rresiduallength
  \ignorespaces}
 
\def\!hleaders{%
  \def\!Rule##1{\vrule height\linethickness width##1}%
  \def\!Skip##1{\hskip##1}%
  \leaders\hbox{\!RSlist}\hfill}
 
\def\!hpartialpattern#1{%
  \!dimenA=\!zpt \!dimenB=\!zpt 
  \def\!Rule##1{#1{##1}\vrule height\linethickness width\!dimenD}%
  \def\!Skip##1{#1{##1}\hskip\!dimenD}%
  \!RSlist}
 
\def\!vleaders{%
  \def\!Rule##1{\hrule width\linethickness height##1}%
  \def\!Skip##1{\vskip##1}%
  \leaders\vbox{\!RSlist}\vfill}
 
\def\!vpartialpattern#1{%
  \!dimenA=\!zpt \!dimenB=\!zpt 
  \def\!Rule##1{#1{##1}\hrule width\linethickness height\!dimenD}%
  \def\!Skip##1{#1{##1}\vskip\!dimenD}%
  \!RSlist}
 
\def\!Rtrunc#1{\!trunc{#1}>\!Rresiduallength}
\def\!Ltrunc#1{\!trunc{#1}<\!Lresiduallength}
 
\def\!trunc#1#2#3{%          
  \!dimenA=\!dimenB         
  \advance\!dimenB by #1%
  \!dimenD=\!dimenB  \ifdim\!dimenD#2#3\!dimenD=#3\fi
  \!dimenC=\!dimenA  \ifdim\!dimenC#2#3\!dimenC=#3\fi
  \advance \!dimenD by -\!dimenC}

\def\!start (#1,#2){%
  \!plotxorigin=\!xorigin  \advance \!plotxorigin by \!plotsymbolxshift
  \!plotyorigin=\!yorigin  \advance \!plotyorigin by \!plotsymbolyshift
  \!xS=\!M{#1}\!xunit \!yS=\!M{#2}\!yunit
  \!rotateaboutpivot\!xS\!yS
  \!copylist\!UDlist\to\!!UDlist% **\!UDlist has the form \\{dimen1}\\{dimen2}..
%                                 ** Routine will draw dashed line with pen
%                                 ** down for dimen1, up for dimen2, ...
  \!getnextvalueof\!downlength\from\!!UDlist
  \!distacross=\!zpt%             ** 1st point goes at start of curve
  \!intervalno=0 %                ** initialize interval counter
  \global\totalarclength=\!zpt%   ** initialize distance traveled along curve
  \ignorespaces}

%  **  \!ljoin (XCOORD,YCOORD) 
%  **  Draws a straight line starting at the last point specified
%  **    by the most recent \!start, \!ljoin, or \!qjoin, and
%  **    ending at (XCOORD,YCOORD).
\def\!ljoin (#1,#2){%
  \advance\!intervalno by 1
  \!xE=\!M{#1}\!xunit \!yE=\!M{#2}\!yunit
  \!rotateaboutpivot\!xE\!yE
  \!xdiff=\!xE \advance \!xdiff by -\!xS%**  xdiff = xE - xS
  \!ydiff=\!yE \advance \!ydiff by -\!yS%**  ydiff = yE - yS
  \!Pythag\!xdiff\!ydiff\!arclength%     **  arclength = sqrt(xdiff**2+ydiff**2) 
  \global\advance \totalarclength by \!arclength%
  \!drawlinearsegment%   ** set by dashpat to \!linearsolid or \!lineardashed
  \!xS=\!xE \!yS=\!yE%   ** shift ending points to starting points
  \ignorespaces}

% **  The following routine is used to draw a "solid" line between (xS,yS)
% **  and (xE,yE).  Points are spaced nearly every  \plotsymbolspacing length
% **  along the line.  
\def\!linearsolid{%
  \!npoints=\!arclength
  \!countA=\plotsymbolspacing
  \divide\!npoints by \!countA%      ** now #pts =. arclength/plotsymbolspacing
  \ifnum \!npoints<1 
    \!npoints=1 
  \fi
  \divide\!xdiff by \!npoints
  \divide\!ydiff by \!npoints
  \!xpos=\!xS \!ypos=\!yS
  \loop\ifnum\!npoints>-1
    \!plotifinbounds
    \advance \!xpos by \!xdiff
    \advance \!ypos by \!ydiff
    \advance \!npoints by -1
  \repeat
  \ignorespaces}

% ** The following routine is used to draw a dashed line between (xS,yS)
% ** and (xE,yE). The dash pattern continues from the previous segment.
\def\!lineardashed{%
% **
  \ifdim\!distacross>\!arclength
    \advance \!distacross by -\!arclength  %nothing to plot in this interval
  \else
    \loop\ifdim\!distacross<\!arclength
%     ** plot point, interpolating linearly in x and y
      \!divide\!distacross\!arclength\!dimenA%  ** dimA = across/arclength
      \!removept\!dimenA\!t%  ** \!t holds value in dimA, without the "pt"
      \!xpos=\!t\!xdiff \advance \!xpos by \!xS
      \!ypos=\!t\!ydiff \advance \!ypos by \!yS
      \!plotifinbounds
      \advance\!distacross by \plotsymbolspacing
      \!advancedashing
    \repeat  
    \advance \!distacross by -\!arclength%    ** prepare for next interval 
  \fi
  \ignorespaces}

\def\!!advancedashing{%
  \advance\!downlength by -\plotsymbolspacing
  \ifdim \!downlength>\!zpt
  \else
    \advance\!distacross by \!downlength
    \!getnextvalueof\!uplength\from\!!UDlist
    \advance\!distacross by \!uplength
    \!getnextvalueof\!downlength\from\!!UDlist
  \fi}

% ** \inboundscheckoff & \inboundscheckon: See Subsection 5.5 of the manual.
\def\inboundscheckoff{%
  \def\!plotifinbounds{\!plot(\!xpos,\!ypos)}%
  \def\!initinboundscheck{\relax}\ignorespaces}
 
\inboundscheckoff
 
% ** The following code plots the current point only if it falls in the
% ** current plotarea.  It doesn't matter if the coordinate system has
% ** changed since the plotarea was set up.  However, shifts of the plot
% ** are ignored (how the plotsymbol stands relative to its plot position is
% ** unknown anyway).
\def\!!plotifinbounds{%
  \ifdim \!xpos<\!checkleft
  \else
    \ifdim \!xpos>\!checkright
    \else
      \ifdim \!ypos<\!checkbot
      \else
         \ifdim \!ypos>\!checktop
         \else
           \!plot(\!xpos,\!ypos)
         \fi 
      \fi
    \fi
  \fi}

\def\!!initinboundscheck{%
  \!checkleft=\!arealloc     \advance\!checkleft by \!xorigin
  \!checkright=\!arearloc    \advance\!checkright by \!xorigin
  \!checkbot=\!areabloc      \advance\!checkbot by \!yorigin
  \!checktop=\!areatloc      \advance\!checktop by \!yorigin}

% *********************************
% *** LOGTEN  (Log_10 function) ***
% *********************************
%
% ** \!logten{X}
% ** Calculates log_10 of X.  X and LOG10(X) are in fixed point notation.
% **  X must be positive; it may have an optional `+' sign; any number
% **  of digits may be specified for X.  The absolute error in LOG10(X) is
% **  less than .0001 (probably < .00006).  That's about as good as you
% **  hope for, since TEX only operates to 5 figures after the decimal
% **  point anyway.

%  \!rootten=3.162278pt       **** These are values are set in ALLOCATIONS
%  \!tenAe=2.543275pt  (=A5)
%  \!tenAc=2.773839pt  (=A3)
%  \!tenAa=8.690286pt  (=A1)

\def\!logten#1#2{%
  \expandafter\!!logten#1\!nil
  \!removept\!dimenF#2%
  \ignorespaces}

\def\!!logten#1#2\!nil{%
  \if -#1%
    \!dimenF=\!zpt
    \def\!next{\ignorespaces}%
  \else
    \if +#1%
      \def\!next{\!!logten#2\!nil}%
    \else
      \if .#1%
        \def\!next{\!!logten0.#2\!nil}%
      \else
        \def\!next{\!!!logten#1#2..\!nil}%
      \fi
    \fi
  \fi
  \!next}

\def\!!!logten#1#2.#3.#4\!nil{%
  \!dimenF=1pt %                 ** DimF holds log10 original argument
  \if 0#1%                      
    \!!logshift#3pt %            ** Argument < 1
  \else %                        ** Argument >= 1
    \!logshift#2/%               ** Shift decimal pt as many places
    \!dimenE=#1.#2#3pt %         **   as there are figures in #2
  \fi %                          ** Now dimE holds revised X want log10 of
  \ifdim \!dimenE<\!rootten%          ** Transform X to XX between sqrt(10) 
    \multiply \!dimenE 10 %           **   and 10*sqrt(10)
    \advance  \!dimenF -1pt
  \fi
  \!dimenG=\!dimenE%                  ** dimG <- (XX + 10)
    \advance\!dimenG 10pt
  \advance\!dimenE -10pt %            ** dimE <- (XX - 10)
  \multiply\!dimenE 10 %              ** dimE = 10*(XX-10)
  \!divide\!dimenE\!dimenG\!dimenE%   ** Now dimE=10t==10*(XX-10)/(XX+10)
  \!removept\!dimenE\!t%              ** !t=10t, with "pt" removed
  \!dimenG=\!t\!dimenE%               ** dimG=100t**2
  \!removept\!dimenG\!tt%             ** !tt=100t**2, with "pt" removed
  \!dimenH=\!tt\!tenAe%               ** dimH=10*a5*(10t)**2 /100
    \divide\!dimenH 100
  \advance\!dimenH \!tenAc%           ** ditto + 10*a3
  \!dimenH=\!tt\!dimenH%              ** ditto * (10t)**2 /100
    \divide\!dimenH 100   
  \advance\!dimenH \!tenAa%           ** ditto + 10*a1
  \!dimenH=\!t\!dimenH%               ** ditto * 10t / 100
    \divide\!dimenH 100 %             ** Now dimH = log10(XX) - 1
  \advance\!dimenF \!dimenH}%         ** dimF = log10(X)

\def\!logshift#1{%
  \if #1/%
    \def\!next{\ignorespaces}%
  \else
    \advance\!dimenF 1pt 
    \def\!next{\!logshift}%
  \fi 
  \!next}
 
 \def\!!logshift#1{%
   \advance\!dimenF -1pt
   \if 0#1%
     \def\!next{\!!logshift}%
   \else
     \if p#1%
       \!dimenF=1pt
       \def\!next{\!dimenE=1p}%
     \else
       \def\!next{\!dimenE=#1.}%
     \fi
   \fi
   \!next}

\def\beginpicture{%
  \setbox\!picbox=\hbox\bgroup%
  \!xleft=\maxdimen  
  \!xright=-\maxdimen
  \!ybot=\maxdimen
  \!ytop=-\maxdimen}
 
% **  \endpicture : See Subsection 1.1 of the manual.
\def\endpicture{%
  \ifdim\!xleft=\maxdimen%  ** check if nothing was put in picbox
    \!xleft=\!zpt \!xright=\!zpt \!ybot=\!zpt \!ytop=\!zpt 
  \fi
  \global\!Xleft=\!xleft \global\!Xright=\!xright
  \global\!Ybot=\!ybot \global\!Ytop=\!ytop
  \egroup%
  \ht\!picbox=\!Ytop  \dp\!picbox=-\!Ybot
  \ifdim\!Ybot>\!zpt
  \else 
    \ifdim\!Ytop<\!zpt
      \!Ybot=\!Ytop
    \else
      \!Ybot=\!zpt
    \fi
  \fi
  \hbox{\kern-\!Xleft\lower\!Ybot\box\!picbox\kern\!Xright}}
 
% **  \endpicturesave : See Subsection 8.4 of the manual.
\def\endpicturesave <#1,#2>{%
  \endpicture \global #1=\!Xleft \global #2=\!Ybot \ignorespaces}

% **   \setcoordinatesystem units <XUNIT,YUNIT> 
% **     point at XREF YREF  
% **   Each of `units <XUNIT,YUNIT>' and `point at XREF YREF' 
% **     are optional.
% **   Unit lengths must be given in dimensions (e.g., <10pt,1in>).
% **     Default unit lengths are 1pt, 1pt, or previous unit lengths.
% **   Reference point is specified in current units (e.g., 3 5 ). 
% **     Default reference point is 0 0 , or previous reference point.
% **   Unit lengths and reference points obey TEX's scoping rules.
% **   See Subsection 1.2 of the manual.
\def\setcoordinatesystem{%
  \!ifnextchar{u}{\!getlengths }
    {\!getlengths units <\!xunit,\!yunit>}}
\def\!getlengths units <#1,#2>{%
  \!xunit=#1\relax
  \!yunit=#2\relax
  \!ifcoordmode 
    \let\!SCnext=\!SCccheckforRP
  \else
    \let\!SCnext=\!SCdcheckforRP
  \fi
  \!SCnext}
\def\!SCccheckforRP{%
  \!ifnextchar{p}{\!cgetreference }
    {\!cgetreference point at {\!xref} {\!yref} }}
\def\!cgetreference point at #1 #2 {%
  \edef\!xref{#1}\edef\!yref{#2}%
  \!xorigin=\!xref\!xunit  \!yorigin=\!yref\!yunit  
  \!initinboundscheck % ** See linear.tex
  \ignorespaces}
\def\!SCdcheckforRP{%
  \!ifnextchar{p}{\!dgetreference}%
    {\ignorespaces}}
\def\!dgetreference point at #1 #2 {%
  \!xorigin=#1\relax  \!yorigin=#2\relax
  \ignorespaces}

%  ** \put {OBJECT} [XY] <XDIMEN,YDIMEN> at (XCOORD,YCOORD)
%  **   `[XY]' and `<XDIMEN,YDIMEN>' are optional.
%  **   First OBJECT is placed in an hbox (the "objectbox") and then a
%  **     "reference point" is assigned to the objectbox as follows:
%  **     [1] first, the reference point is taken to be the center of the box;
%  **     [2] next, centering is overridden by the specifications
%  **           X=l -- reference point along the left edge of the objectbox
%  **           X=r -- reference point along the right edge of the objectbox
%  **           Y=b -- reference point along the bottom edge of the objectbox
%  **           Y=B -- reference point along the Baseline of the objectbox
%  **           Y=t -- reference point along the top edge of the objectbox;
%  **     [3] finally the reference point is shifted left by XDIMEN, down
%  **           by YDIMEN  (both default to 0pt).
%  **   The objectbox is placed within PICBOX with its reference point at  
%  **     (XCOORD,YCOORD). 
%  **   If OBJECT is a saved box, say  box0, you have to write
%  **     \put{\box0}...   or  \put{\copy0}...
%  **   The objectbox is void after the put.
%  **   See Subsection 2.1 of the manual.
\long\def\put#1#2 at #3 #4 {%
  \!setputobject{#1}{#2}%
  \!xpos=\!M{#3}\!xunit  \!ypos=\!M{#4}\!yunit  
  \!rotateaboutpivot\!xpos\!ypos%
  \advance\!xpos -\!xorigin  \advance\!xpos -\!xshift
  \advance\!ypos -\!yorigin  \advance\!ypos -\!yshift
  \kern\!xpos\raise\!ypos\box\!putobject\kern-\!xpos%
  \!doaccounting\ignorespaces}
 
%  **   \multiput etc.  Like  \put.  The objectbox is not voided until the
%  **     termininating /, and is placed repeatedly with:
%  **     XCOORD YCOORD -- the objectbox is put down with its reference point
%  **       at (XCOORD,YCOORD);
%  **     *N DXCOORD DYCOORD -- each of N times the current
%  **       (xcoord,ycoord) is incremented by (DXCOORD,DYCOORD), and the
%  **       objectbox is put down with its reference point at (xcoord,ycoord)
%  **       (This specification has to follow an XCOORD YCOORD pair)
%  **     See Subsection 2.2 of the manual.
\long\def\multiput #1#2 at {%
  \!setputobject{#1}{#2}%
  \!ifnextchar"{\!putfromfile}{\!multiput}}
\def\!putfromfile"#1"{%
  \expandafter\!multiput \input #1 /}
\def\!multiput{%
  \futurelet\!nextchar\!!multiput}
\def\!!multiput{%
  \if *\!nextchar
    \def\!nextput{\!alsoby}%
  \else
    \if /\!nextchar
      \def\!nextput{\!finishmultiput}%
    \else
      \def\!nextput{\!alsoat}%
    \fi
  \fi
  \!nextput}
\def\!finishmultiput/{%
  \setbox\!putobject=\hbox{}%
  \ignorespaces}
 
%  **   \!alsoat XCOORD YCOORD 
%  **     The objectbox is put down with reference point at XCOORD,YCOORD
\def\!alsoat#1 #2 {%
  \!xpos=\!M{#1}\!xunit  \!ypos=\!M{#2}\!yunit  
  \!rotateaboutpivot\!xpos\!ypos%
  \advance\!xpos -\!xorigin  \advance\!xpos -\!xshift
  \advance\!ypos -\!yorigin  \advance\!ypos -\!yshift
  \kern\!xpos\raise\!ypos\copy\!putobject\kern-\!xpos%
  \!doaccounting
  \!multiput}
 
% **   \!alsoby*N DXCOORD DYCOORD
% **     N times, the current (XCOORD,YCOORD) is advanced by (DXCOORD,DYCOORD),
% **     and the current (shifted, oriented) OBJECT is put down.
\def\!alsoby*#1 #2 #3 {%
  \!dxpos=\!M{#2}\!xunit \!dypos=\!M{#3}\!yunit 
  \!rotateonly\!dxpos\!dypos
  \!ntemp=#1%
  \!!loop\ifnum\!ntemp>0
    \advance\!xpos by \!dxpos  \advance\!ypos by \!dypos
    \kern\!xpos\raise\!ypos\copy\!putobject\kern-\!xpos%
    \advance\!ntemp by -1
  \repeat
  \!doaccounting 
  \!multiput}
 
% **  \accountingoff : Suspends PiCTeX's accounting of the aggregate
% **    size of the picture box.
% **  \accounting on : Reinstates accounting.
% **  See Subsection 8.2 of the manual.
\def\accountingon{\def\!doaccounting{\!!doaccounting}\ignorespaces}

\accountingon
\def\!!doaccounting{%
  \!xtemp=\!xpos  
  \!ytemp=\!ypos
  \ifdim\!xtemp<\!xleft 
     \!xleft=\!xtemp 
  \fi
  \advance\!xtemp by  \!wd 
  \ifdim\!xright<\!xtemp 
    \!xright=\!xtemp
  \fi
  \advance\!ytemp by -\!dp
  \ifdim\!ytemp<\!ybot  
    \!ybot=\!ytemp
  \fi
  \advance\!ytemp by  \!dp
  \advance\!ytemp by  \!ht 
  \ifdim\!ytemp>\!ytop  
    \!ytop=\!ytemp  
  \fi}
 
\long\def\!setputobject#1#2{%
  \setbox\!putobject=\hbox{#1}%
  \!ht=\ht\!putobject  \!dp=\dp\!putobject  \!wd=\wd\!putobject
  \wd\!putobject=\!zpt
  \!xshift=.5\!wd   \!yshift=.5\!ht   \advance\!yshift by -.5\!dp
  \edef\!putorientation{#2}%
  \expandafter\!SPOreadA\!putorientation[]\!nil%
  \expandafter\!SPOreadB\!putorientation<\!zpt,\!zpt>\!nil\ignorespaces}
 
\def\!SPOreadA#1[#2]#3\!nil{\!etfor\!orientation:=#2\do\!SPOreviseshift}
 
\def\!SPOreadB#1<#2,#3>#4\!nil{\advance\!xshift by -#2\advance\!yshift by -#3}
 
\def\!SPOreviseshift{%
  \if l\!orientation 
    \!xshift=\!zpt
  \else 
    \if r\!orientation 
      \!xshift=\!wd
    \else 
      \if b\!orientation
        \!yshift=-\!dp
      \else 
        \if B\!orientation 
          \!yshift=\!zpt
        \else 
          \if t\!orientation 
            \!yshift=\!ht
          \fi 
        \fi
      \fi
    \fi
  \fi}

%  **  \!dimenput{OBJECT} <XDIMEN,YDIMEN> [XY] (XLOC,YLOC)
%  **    This is an internal put routine, similar to \put, except that
%  **    XLOC=distance right from reference point, YLOC=distance up from
%  **    reference point. XLOC and YLOC are dimensions, so this routine
%  **    is completely independent of the current coordinate system. 
%  **    This routine does NOT do ROTATIONS.
\long\def\!dimenput#1#2(#3,#4){%
  \!setputobject{#1}{#2}%
  \!xpos=#3\advance\!xpos by -\!xshift
  \!ypos=#4\advance\!ypos by -\!yshift
  \kern\!xpos\raise\!ypos\box\!putobject\kern-\!xpos%
  \!doaccounting\ignorespaces}

%  ** The following macros permit the picture drawing routines to be used 
%  ** either in the default "coordinate mode", or in "dimension mode".
%  **   In coordinate mode  \!M(1.5,\!xunit)    expands to  1.5\!xunit
%  **   In dimension  mode  \!M(1.5pt,\!xunit)  expands to  1.5pt
%  ** Dimension mode is useful in coding macros.
%  ** Any special purpose picture macro that sets dimension mode should 
%  ** reset coordinate mode before completion.
%  ** See Subsection 9.2 of the manual.
\def\!setdimenmode{%
  \let\!M=\!M!!\ignorespaces}
\def\!setcoordmode{%
  \let\!M=\!M!\ignorespaces}
\def\!ifcoordmode{%
  \ifx \!M \!M!}
\def\!ifdimenmode{%
  \ifx \!M \!M!!}
\def\!M!#1#2{#1#2} 
\def\!M!!#1#2{#1}
\!setcoordmode
\let\setdimensionmode=\!setdimenmode
\let\setcoordinatemode=\!setcoordmode

%  ** \Xdistance{XCOORD}, \Ydistance{YCOORD}  are the horizontal and
%  **   vertical distances from the origin (0,0) to the point
%  **   (XCOORD,YCOORD)  in the current coordinate system.
%  ** See Subsection 9.2 of the manual.

% ** The following macros -- \stack, \line, and \Lines -- are useful for
% **   annotating PiCtures. They can be used outside the \beginpicture ...
% **   \endpicture environment.

% ** \stack [POSITIONING] <LEADING> {VALUESLIST}
% ** Builds a vertical stack of the values in VALUESLIST. Values in
% ** VALUESLIST are separated by commas.  In the resulting stack, values are
% ** centered by default, and positioned flush left (right) if 
% ** POSITIONING = l (r).  Values are separated vertically by LEADING,
% ** which defaults to \stackleading.
% ** See Subsection 2.3 of the manual.

\def\!stack[#1]{%
  \let\!lglue=\hfill \let\!rglue=\hfill
  \expandafter\let\csname !#1glue\endcsname=\relax
  \!ifnextchar<{\!!stack}{\!!stack<\stackleading>}}
\def\!!stack<#1>#2{%
  \vbox{\def\!valueslist{}\!ecfor\!value:=#2\do{%
    \expandafter\!rightappend\!value\withCS{\\}\to\!valueslist}%
    \!lop\!valueslist\to\!value
    \let\\=\cr\lineskiplimit=\maxdimen\lineskip=#1%
    \baselineskip=-1000pt\halign{\!lglue##\!rglue\cr \!value\!valueslist\cr}}%
  \ignorespaces}

% ** \lines [POSITIONING] {LINES}
% ** Builds a vertical array of the lines in LINES. Each line in LINES
% ** is terminated by a \cr.  In the resulting array, lines are
% ** centered by default, and positioned flush left (right) if 
% ** POSITIONING = l (r).  The lines in the array are subject to TeX's
% ** usual spacing rules: in particular the baselines are ordinarily an equal
% ** distance apart. The baseline of the array is the baseline of the
% ** the bottom line.
% ** See Subsection 2.3 of the manual.

\def\!lines[#1]#2{%
  \let\!lglue=\hfill \let\!rglue=\hfill
  \expandafter\let\csname !#1glue\endcsname=\relax
  \vbox{\halign{\!lglue##\!rglue\cr #2\crcr}}%
  \ignorespaces}

% ** \Lines [POSITIONING] {LINES}
% ** Like \lines, but the baseline of the array is the baseline of the
% ** top line.  See Subsection 2.3 of the manual.

\def\!Lines[#1]#2{%
  \let\!lglue=\hfill \let\!rglue=\hfill
  \expandafter\let\csname !#1glue\endcsname=\relax
  \vtop{\halign{\!lglue##\!rglue\cr #2\crcr}}%
  \ignorespaces}

% *********************************************
% *** PLOTTING (Things to do with plotting) ***
% *********************************************
 
% **  User commands
% **    \setplotsymbol ({PLOTSYMBOL} [ORIENTATION] <XSHIFT,YSHIFT>)
% **    \savelinesandcurves on "FILE_NAME"
% **    \dontsavelinesandcurves
% **    \writesavefile {MESSAGE}
% **    \replot {FILE_NAME}
 
% **  Internal command
% **    \!plot(XDIMEN,YDIMEN)
 
% **  \setplotsymbol ({PLOTSYMBOL} [ ] < , >)
% **  Save PLOTSYMBOL away in an hbox for use with curve plotting routines
% **  See Subsection 5.2 of the manual.
\def\setplotsymbol(#1#2){%
  \!setputobject{#1}{#2}
  \setbox\!plotsymbol=\box\!putobject%
  \!plotsymbolxshift=\!xshift 
  \!plotsymbolyshift=\!yshift 
  \ignorespaces}
 
\setplotsymbol({\fiverm .})%       ** initialize plotsymbol

% **  \!plot is either \!!plot (when no lines and curves are being saved) or
% **                   \!!!plot (when   lines and curves are being saved)
 
% **  \!!plot(XDIMEN,YDIMEN)
% **  Places the current plotsymbol a horizontal distance=XDIMEN-xorigin 
% **    and a vertical distance=YDIMEN-yorigin from the current
% **    reference point.  
\def\!!plot(#1,#2){%
  \!dimenA=-\!plotxorigin \advance \!dimenA by #1%    ** over
  \!dimenB=-\!plotyorigin \advance \!dimenB by #2%    ** up
  \kern\!dimenA\raise\!dimenB\copy\!plotsymbol\kern-\!dimenA%
  \ignorespaces}
 
% **  \!!!plot(XDIMEN,YDIMEN)
% **  Like \!!plot, but also saves the plot location in units of 
% **    scaled point, on file `replotfile'
\def\!!!plot(#1,#2){%
  \!dimenA=-\!plotxorigin \advance \!dimenA by #1%    ** over
  \!dimenB=-\!plotyorigin \advance \!dimenB by #2%    ** up
  \kern\!dimenA\raise\!dimenB\copy\!plotsymbol\kern-\!dimenA%
  \!countE=\!dimenA
  \!countF=\!dimenB
  \immediate\write\!replotfile{\the\!countE,\the\!countF.}%
  \ignorespaces}

% ** \savelinesandcurves on "FILE_NAME"
% **   Switch to save locations used for plotting lines and curves
% **   (No advantage in saving locations for solid lines; however
% **   replotting curve locations speeds things up by a factor of about 4. 
% ** \dontsavelinesandcurves
% **   Terminates \savelinesandcurves. The default.
% ** See Subsection 5.6 of the manual.
\def\savelinesandcurves on "#1" {%
  \immediate\closeout\!replotfile
  \immediate\openout\!replotfile=#1%
  \let\!plot=\!!!plot}

\def\dontsavelinesandcurves {%
  \let\!plot=\!!plot}
\dontsavelinesandcurves

% ** \writesavefile {MESSAGE}
% ** The message is preceded by a "%", so that it won't interfere
% ** with replotting.
% ** See Subsection 5.6 of the manual.
{\catcode`\%=11\xdef\!Commentsignal{%}}
\def\writesavefile#1 {%
  \immediate\write\!replotfile{\!Commentsignal #1}%
  \ignorespaces}

% ** \replot "FILE_NAME"
% **   Replots the locations saved earlier under \savelinesandcurves
% **   on "FILE_NAME"
% ** See Subsection 5.6 of the manual.
\def\replot"#1" {%
  \expandafter\!replot\input #1 /}
\def\!replot#1,#2. {%
  \!dimenA=#1sp
  \kern\!dimenA\raise#2sp\copy\!plotsymbol\kern-\!dimenA
  \futurelet\!nextchar\!!replot}
\def\!!replot{%
  \if /\!nextchar 
    \def\!next{\!finish}%
  \else
    \def\!next{\!replot}%
  \fi
  \!next}
% **************************************************
% ***  PYTHAGORAS  (Euclidean distance function) ***
% **************************************************

% ** User command:
% **   \placehypotenuse for <dimension1> and <dimension2> in <register> 

% ** Internal command:
% **   \!Pythag{X}{Y}{Z}
% **     Input X,Y are dimensions, or dimension registers.
% **     Output Z == sqrt(X**2+Y**2) must be a dimension register.
% **     Assumes that |X|+|Y| < 2048pt (about 28in).
 
% ** Without loss of generality, suppose  x>0, y>0.  Put s = x+y,
% **   z = sqrt(x**2+y**2). Then  z = s*f,  where  f = sqrt(t**2 + (1-t)**2)
% **   = sqrt((1+tau**2)/2), where  t = x/s  and  tau = 2(t-1/2) .
 
% ** Uses the \!divide macro (which uses registers \!dimenA--\!dimenD.
% ** Uses the \!removept macro   (e.g., 123.45pt --> 123.45)
% ** Uses registers \!dimenE--\!dimenI.
\def\!Pythag#1#2#3{%
  \!dimenE=#1\relax                                     
  \ifdim\!dimenE<\!zpt 
    \!dimenE=-\!dimenE 
  \fi%                                            ** dimE = |x|
  \!dimenF=#2\relax
  \ifdim\!dimenF<\!zpt 
    \!dimenF=-\!dimenF 
  \fi%                                            ** dimF = |y|
  \advance \!dimenF by \!dimenE%                  ** dimF = s = |x|+|y|
  \ifdim\!dimenF=\!zpt 
    \!dimenG=\!zpt%                               ** dimG = z = sqrt(x**2+y**2)
  \else 
    \!divide{8\!dimenE}\!dimenF\!dimenE%          ** now dimE = 8t = (8|x|)/s
    \advance\!dimenE by -4pt%                     ** 8tau = (8t-4)*2
      \!dimenE=2\!dimenE%                         **   (tau = 2*t - 1)
    \!removept\!dimenE\!!t%                       ** 8tau, without "pt"
    \!dimenE=\!!t\!dimenE%                        ** (8tau)**2, in pts
    \advance\!dimenE by 64pt%                     ** u = [64 + (8tau)**2]/2
    \divide \!dimenE by 2%                        **   [u = (8f)**2]
    \!dimenH=7pt%                                 ** initial guess g at sqrt(u)
    \!!Pythag\!!Pythag\!!Pythag%                  ** 3 iterations give sqrt(u)
    \!removept\!dimenH\!!t%                       ** 8f=sqrt(u), without "pt"
    \!dimenG=\!!t\!dimenF%                        ** z = (8f)*s/8
    \divide\!dimenG by 8
  \fi
  #3=\!dimenG
  \ignorespaces}

\def\!!Pythag{%                                   ** Newton-Raphson for sqrt
  \!divide\!dimenE\!dimenH\!dimenI%               ** v = u/g
  \advance\!dimenH by \!dimenI%                   ** g <-- (g + u/g)/2
    \divide\!dimenH by 2}

% **  \placehypotenuse for <XI> and <ETA> in <ZETA>
% **  See Subsection 9.3 of the manual.
\def\placehypotenuse for <#1> and <#2> in <#3> {%
  \!Pythag{#1}{#2}{#3}}

% **********************************************
% *** QUADRATIC ARC  (Draws a quadratic arc) ***
% **********************************************
 
% **  Internal command
% **    \!qjoin (XCOORD1,YCOORD1) (XCOORD2,YCOORD2)
 
% **  \!qjoin (XCOORD1,YCOORD1) (XCOORD2,YCOORD2)
% **  Draws an arc starting at the (last) point specified by the most recent
% **  \!qjoin, or \!ljoin, or \!start  and passing through (X_1,Y_1), (X_2,Y_2).
% **  Uses quadratic interpolation in both  x  and  y:  
% **    x(t), 0 <= t <= 1, interpolates  x_0, x_1, x_2  at  t=0, .5, 1
% **    y(t), 0 <= t <= 1, interpolates  y_0, y_1, y_2  at  t=0, .5, 1
 
\def\!qjoin (#1,#2) (#3,#4){%
  \advance\!intervalno by 1
  \!ifcoordmode
    \edef\!xmidpt{#1}\edef\!ymidpt{#2}%
  \else
    \!dimenA=#1\relax \edef\!xmidpt{\the\!dimenA}%
    \!dimenA=#2\relax \edef\!xmidpt{\the\!dimenA}%
  \fi
  \!xM=\!M{#1}\!xunit  \!yM=\!M{#2}\!yunit   \!rotateaboutpivot\!xM\!yM
  \!xE=\!M{#3}\!xunit  \!yE=\!M{#4}\!yunit   \!rotateaboutpivot\!xE\!yE
%
% ** Find coefficients for x(t)=a_x + b_x*t + c_x*t**2
  \!dimenA=\!xM  \advance \!dimenA by -\!xS%   ** dimA = I = xM - xS
  \!dimenB=\!xE  \advance \!dimenB by -\!xM%   ** dimB = II = xE-xM
  \!xB=3\!dimenA \advance \!xB by -\!dimenB%   ** b=3I-II
  \!xC=2\!dimenB \advance \!xC by -2\!dimenA%  ** c=2(II-I)
%
% ** Find coefficients for y(t)=y_x + b_y*t + c_y*t**2
  \!dimenA=\!yM  \advance \!dimenA by -\!yS%   
  \!dimenB=\!yE  \advance \!dimenB by -\!yM%  
  \!yB=3\!dimenA \advance \!yB by -\!dimenB%  
  \!yC=2\!dimenB \advance \!yC by -2\!dimenA% 
%
% ** Use Simpson's rule to calculate arc length over [0,1/2]:
% **   arc length = 1/2[1/6 f(0) + 4/6 f(1/4) + 1/6 f(1/2)]
% ** with f(t) = sqrt(x'(t)**2 + y'(t)**2).
  \!xprime=\!xB  \!yprime=\!yB%          ** x'(t) = b + 2ct
  \!dxprime=.5\!xC  \!dyprime=.5\!yC%    ** dt=1/4 ==> dx'(t) = c/2
  \!getf \!midarclength=\!dimenA
  \!getf \advance \!midarclength by 4\!dimenA
  \!getf \advance \!midarclength by \!dimenA
  \divide \!midarclength by 12
%
% ** Get arc length over [0,1].
  \!arclength=\!dimenA
  \!getf \advance \!arclength by 4\!dimenA
  \!getf \advance \!arclength by \!dimenA
  \divide \!arclength by 12%             ** Now have arc length over [1/2,1]
  \advance \!arclength by \!midarclength
  \global\advance \totalarclength by \!arclength
%
%
% ** Check to see if there's anything to plot in this interval
  \ifdim\!distacross>\!arclength 
    \advance \!distacross by -\!arclength%   ** nothing 
  \else
    \!initinverseinterp%  ** initialize for inverse interpolation on arc length
    \loop\ifdim\!distacross<\!arclength%     ** loop over points on arc 
      \!inverseinterp%    ** find  t  such that arc length[0,t] = distacross,
%                         **   using inverse quadratic interpolation
%                         ** now evaluate x(t)=(c*t + b)*t + a
      \!xpos=\!t\!xC \advance\!xpos by \!xB
        \!xpos=\!t\!xpos \advance \!xpos by \!xS
%                                             ** evaluate y(t)
      \!ypos=\!t\!yC \advance\!ypos by \!yB
        \!ypos=\!t\!ypos \advance \!ypos by \!yS
      \!plotifinbounds%                       ** plot point if in bounds
      \advance\!distacross \plotsymbolspacing%** advance arc length for next pt
      \!advancedashing%                       ** see "linear"
    \repeat  
    \advance \!distacross by -\!arclength%    ** prepare for next interval 
  \fi
  \!xS=\!xE%              ** shift ending points to starting points
  \!yS=\!yE
  \ignorespaces}

% ** \!getf -- Calculates sqrt(x'(t)**2 + y'(t)**2) and advances
% **   x'(t) and y'(t)
\def\!getf{\!Pythag\!xprime\!yprime\!dimenA%
  \advance\!xprime by \!dxprime
  \advance\!yprime by \!dyprime}

% ** \!initinverseinterp -- initializes for inverse quadratic interpolation
% ** of arc length provided  1/3 < midarclength/arclength < 2/3; otherwise
% ** initializes for inverse linear interpolation.
\def\!initinverseinterp{%
  \ifdim\!arclength>\!zpt
    \!divide{8\!midarclength}\!arclength\!dimenE% ** dimE=8w=8r/s, where  r 
%                                               **  = midarclength, s=arclength
% **  Test for  w  out of range:  w<1/3  or w>2/3
    \ifdim\!dimenE<\!wmin \!setinverselinear
    \else 
      \ifdim\!dimenE>\!wmax \!setinverselinear
      \else%                                    ** w  in range: initialize
        \def\!inverseinterp{\!inversequad}\ignorespaces
%
% **     Calculate the coefficients  \!beta  and  \!gamma  of the quadratic
% **                    t = \!beta*v + \!gamma*v**2
% **     taking the values  t=0, 1/2, 1  at  v=0, w==r/s, 1  respectively:
% **        \!beta = (1/2 - w**2)/[w(1-w)] 
% **        \!gamma = 1 - beta.
%
         \!removept\!dimenE\!Ew%           **  8w, without "pt"
         \!dimenF=-\!Ew\!dimenE%           **  -(8w)**2
         \advance\!dimenF by 32pt%         **  32 - (8w)**2
         \!dimenG=8pt 
         \advance\!dimenG by -\!dimenE%    **  8 - 8w
         \!dimenG=\!Ew\!dimenG%            **  (8w)*(8-8w)
         \!divide\!dimenF\!dimenG\!beta%   **  beta = (32-(8w)**2)/(8w(8-8w))
%                                          **       = (1/2 - w**2)/(w(1-w))
         \!gamma=1pt
         \advance \!gamma by -\!beta%      **  gamma = 1-beta
      \fi%       ** end of the \ifdim\!dimenE>\!wmax
    \fi%         ** end of the \ifdim\!dimenE<\!wmin
  \fi%           ** end of the \ifdim\!arclength>\!zpt
  \ignorespaces}

% ** For 0 <= t <= 1, let AL(t) = arclength[0,t]/arclength[0,1]; note
% ** AL(0)=0, AL(1/2)=midarclength/arclength, AL(1)=1.  This routine
% ** calculates an approximation to AL^{-1}(distance across/arclength),
% ** using the assumption that AL^{-1} is quadratic.  Specifically, 
% ** it finds  t  such that
% **    AL^{-1}(v) =. t = v*(\!beta + \!gamma*v)
% ** where  \!beta  and  \!gamma  are set by \!initinv, and where
% ** v=distance across/arclength
\def\!inversequad{%
  \!divide\!distacross\!arclength\!dimenG%   ** dimG = v = distacross/arclength
  \!removept\!dimenG\!v%                     ** v, without "pt"
  \!dimenG=\!v\!gamma%                       ** gamma*v
  \advance\!dimenG by \!beta%                ** beta + gamma*v
  \!dimenG=\!v\!dimenG%                      ** t = v*(beta + gamma*v)
  \!removept\!dimenG\!t}%                    ** t, without "pt"

% ** When  w <= 1/3  or  w >= 2/3, the following routine writes (using
% ** plain TEK's \wlog command) a warning message on the user's log file,
% ** and initializes for inverse linear interpolation on arc length.
\def\!setinverselinear{%
  \def\!inverseinterp{\!inverselinear}%
  \divide\!dimenE by 8 \!removept\!dimenE\!t
  \!countC=\!intervalno \multiply \!countC 2
  \!countB=\!countC     \advance \!countB -1
  \!countA=\!countB     \advance \!countA -1
  \wlog{\the\!countB th point (\!xmidpt,\!ymidpt) being plotted 
    doesn't lie in the}%
  \wlog{ middle third of the arc between the \the\!countA th 
    and \the\!countC th points:}%
  \wlog{ [arc length \the\!countA\space to \the\!countB]/[arc length 
    \the \!countA\space to \the\!countC]=\!t.}%
  \ignorespaces}
 
% **  Inverse linear interpolation
\def\!inverselinear{% 
  \!divide\!distacross\!arclength\!dimenG
  \!removept\!dimenG\!t}

% **************************************
% **  ROTATIONS  (Handles rotations) ***
% **************************************
 
% ** User commands
% **   \startrotation [by COS_OF_ANGLE SIN_OF_ANGLE] [about XPIVOT YPIVOT]
% **   \stoprotation

% **   \startrotation [by COS_OF_ANGLE SIN_OF_ANGLE] [about XPIVOT YPIVOT]
% ** Future (XCOORD,YCOORD)'s will be rotated about (XPIVOT,YPIVOT) 
% ** by the angle with the give COS and SIN. Both fields are optional.
% ** [COS,SIN] defaults to previous value, or (1,0).
% ** (XPIVOT,YPIVOT) defaults to previous value, or (0,0)
% ** You can't change the coordinate system in the scope of a rotation.
% ** See Subsection 9.1 of the manual.
\def\startrotation{%
  \let\!rotateaboutpivot=\!!rotateaboutpivot
  \let\!rotateonly=\!!rotateonly
  \!ifnextchar{b}{\!getsincos }%
    {\!getsincos by {\!cosrotationangle} {\!sinrotationangle} }}
\def\!getsincos by #1 #2 {%
  \edef\!cosrotationangle{#1}%
  \edef\!sinrotationangle{#2}%
  \!ifcoordmode 
    \let\!ROnext=\!ccheckforpivot
  \else
    \let\!ROnext=\!dcheckforpivot
  \fi
  \!ROnext}
\def\!ccheckforpivot{%
  \!ifnextchar{a}{\!cgetpivot}%
    {\!cgetpivot about {\!xpivotcoord} {\!ypivotcoord} }}
\def\!cgetpivot about #1 #2 {%
  \edef\!xpivotcoord{#1}%
  \edef\!ypivotcoord{#2}%
  \!xpivot=#1\!xunit  \!ypivot=#2\!yunit
  \ignorespaces}
\def\!dcheckforpivot{%
  \!ifnextchar{a}{\!dgetpivot}{\ignorespaces}}
\def\!dgetpivot about #1 #2 {%
  \!xpivot=#1\relax  \!ypivot=#2\relax
  \ignorespaces}

% ** Following terminates rotation.
% ** See Subsection 9.1 of the manual.
\def\stoprotation{%
  \let\!rotateaboutpivot=\!!!rotateaboutpivot
  \let\!rotateonly=\!!!rotateonly
  \ignorespaces}
 
% ** !!rotateaboutpivot{XREG}{YREG}
% ** XREG <-- xpvt + cos(angle)*(XREG-xpvt) - sin(angle)*(YREG-ypvt)
% ** YREG <-- ypvt + cos(angle)*(YREG-ypvt) + sin(angle)*(XREG-xpvt)
% ** XREG,YREG are dimension registers. Can't be \!dimenA to \!dimenD
\def\!!rotateaboutpivot#1#2{%
  \!dimenA=#1\relax  \advance\!dimenA -\!xpivot
  \!dimenB=#2\relax  \advance\!dimenB -\!ypivot
  \!dimenC=\!cosrotationangle\!dimenA
    \advance \!dimenC -\!sinrotationangle\!dimenB
  \!dimenD=\!cosrotationangle\!dimenB
    \advance \!dimenD  \!sinrotationangle\!dimenA
  \advance\!dimenC \!xpivot  \advance\!dimenD \!ypivot
  #1=\!dimenC  #2=\!dimenD
  \ignorespaces}

% ** \!!rotateonly{XREG}{YREG}
% ** Like \!!rotateaboutpivot, but with a pivot of  (0,0)
\def\!!rotateonly#1#2{%
  \!dimenA=#1\relax  \!dimenB=#2\relax 
  \!dimenC=\!cosrotationangle\!dimenA
    \advance \!dimenC -\!rotsign\!sinrotationangle\!dimenB
  \!dimenD=\!cosrotationangle\!dimenB
    \advance \!dimenD  \!rotsign\!sinrotationangle\!dimenA
  #1=\!dimenC  #2=\!dimenD
  \ignorespaces}
\def\!rotsign{}
\def\!!!rotateaboutpivot#1#2{\relax}
\def\!!!rotateonly#1#2{\relax}
\stoprotation

\def\!reverserotateonly#1#2{%
  \def\!rotsign{-}%
  \!rotateonly{#1}{#2}%
  \def\!rotsign{}%
  \ignorespaces}

\def\!getspan span <#1>{%
  \!dshade=#1\relax
  \!ifcoordmode 
    \let\!GRnext=\!GRccheckforAP
  \else
    \let\!GRnext=\!GRdcheckforAP
  \fi
  \!GRnext}
\def\!GRccheckforAP{%
  \!ifnextchar{p}{\!cgetanchor }
    {\!cgetanchor point at {\!xshadesave} {\!yshadesave} }}
\def\!cgetanchor point at #1 #2 {%
  \edef\!xshadesave{#1}\edef\!yshadesave{#2}%
  \!xshade=\!xshadesave\!xunit  \!yshade=\!yshadesave\!yunit
  \ignorespaces}
\def\!GRdcheckforAP{%
  \!ifnextchar{p}{\!dgetanchor}%
    {\ignorespaces}}
\def\!dgetanchor point at #1 #2 {%
  \!xshade=#1\relax  \!yshade=#2\relax
  \ignorespaces}

% **  \setshadesymbol  [<LS, RS, BS, TS>] ({SHADESYMBOL}
% **    <XDIMEN,YDIMEN> [ORIENTATION])
% **  Saves SHADESYMBOL away in an hbox for use with shading routines.
% **  A shade symbol will not be plotted if its plot position comes within
% **    distance LS of the left boundary,  RS of the right boundary,  TS of the
% **    top boundary,  BS of the bottom boundary.  These parameters have 
% **    default values that should work in most cases (see below).
% **    To override a default value, specify the replacement value
% **    in the appropriate subfield of the shrinkages field.
% **    0pt may be coded as  "z" (without the quotes).  To accept a
% **    default value, leave the field empty.  Thus
% **      [,z,,5pt]  sets  LS=default, RS=0pt, BS=default, TS=5pt .
% **    Skipping the shrinkages field accepts all the defaults.
% **  See Subsection 7.1 of the manual.
\def\setshadesymbol{%
  \!ifnextchar<{\!setshadesymbol}{\!setshadesymbol<,,,> }}

\def\!setshadesymbol <#1,#2,#3,#4> (#5#6){%
% **  set the shadesymbol
  \!setputobject{#5}{#6}%                        
  \setbox\!shadesymbol=\box\!putobject%
  \!shadesymbolxshift=\!xshift \!shadesymbolyshift=\!yshift
%
% **  set the shrinkages
  \!dimenA=\!xshift \advance\!dimenA \!smidge% ** default LS = xshift - smidge
  \!override\!dimenA{#1}\!lshrinkage%         
  \!dimenA=\!wd \advance \!dimenA -\!xshift%   ** default RS = width - xshift
    \advance\!dimenA \!smidge%                                  - smidge
    \!override\!dimenA{#2}\!rshrinkage
  \!dimenA=\!dp \advance \!dimenA \!yshift%    ** default BS = depth + yshift
    \advance\!dimenA \!smidge%                                  - smidge
    \!override\!dimenA{#3}\!bshrinkage
  \!dimenA=\!ht \advance \!dimenA -\!yshift%   ** default TS = height - yshift
    \advance\!dimenA \!smidge%                                  - smidge
    \!override\!dimenA{#4}\!tshrinkage
  \ignorespaces}
\def\!smidge{-.2pt}%

% ** \!override{NOMINAL DIMEN}{REPLACEMENT DIMEN}{DIMEN}
% ** Overrides the NOMINAL DIMEN by the REPLACEMENT DIMEN to produce DIMEN,
% ** according to the following rules:
% **   REPLACEMENT DIMEN empty: DIMEN <-- NOMINAL DIMEN
% **   REPLACEMENT DIMEN z:     DIMEN <-- 0pt
% **   otherwise:               DIMEN <-- REPLACEMENT DIMEN
% ** DIMEN must be a dimension register
\def\!override#1#2#3{%
  \edef\!!override{#2}% 
  \ifx \!!override\empty
    #3=#1\relax
  \else
    \if z\!!override
      #3=\!zpt
    \else
      \ifx \!!override\!blankz
        #3=\!zpt
      \else
        #3=#2\relax
      \fi
    \fi
  \fi
  \ignorespaces}
\def\!blankz{ z}

\setshadesymbol ({\fiverm .})%       ** initialize plotsymbol
%                                    ** \fivesy ^^B  is a small cross

% ** \!startvshade [at] (xS,ybS,ytS)
% ** Initiates vertical shading mode
\def\!startvshade#1(#2,#3,#4){%
  \let\!!xunit=\!xunit%
  \let\!!yunit=\!yunit%
  \let\!!xshade=\!xshade%
  \let\!!yshade=\!yshade%
  \def\!getshrinkages{\!vgetshrinkages}%
  \let\!setshadelocation=\!vsetshadelocation%
  \!xS=\!M{#2}\!!xunit
  \!ybS=\!M{#3}\!!yunit
  \!ytS=\!M{#4}\!!yunit
  \!shadexorigin=\!xorigin  \advance \!shadexorigin \!shadesymbolxshift
  \!shadeyorigin=\!yorigin  \advance \!shadeyorigin \!shadesymbolyshift
  \ignorespaces}
 
% ** \!starthshade [at] (yS,xlS,xrS)
% ** Initiates horizontal shading mode
\def\!starthshade#1(#2,#3,#4){%
  \let\!!xunit=\!yunit%
  \let\!!yunit=\!xunit%
  \let\!!xshade=\!yshade%
  \let\!!yshade=\!xshade%
  \def\!getshrinkages{\!hgetshrinkages}%
  \let\!setshadelocation=\!hsetshadelocation%
  \!xS=\!M{#2}\!!xunit
  \!ybS=\!M{#3}\!!yunit
  \!ytS=\!M{#4}\!!yunit
  \!shadexorigin=\!xorigin  \advance \!shadexorigin \!shadesymbolxshift
  \!shadeyorigin=\!yorigin  \advance \!shadeyorigin \!shadesymbolyshift
  \ignorespaces}

% **  \!lattice{ANCHOR}{SPAN}{LOCATION}{INDEX}{LATTICE LOCATION}
% **  Consider the lattice with points  ANCHOR + j*SPAN. This routine determines
% **  the index  k  of the smallest lattice point >= LOCATION, and sets
% **  LATTICE LOCATION = ANCHOR + k*SPAN.
% **  INDEX is assumed to be a count register, LATTICE LOCATION a dimen reg.
\def\!lattice#1#2#3#4#5{%
  \!dimenA=#1%                        ** dimA = ANCHOR
  \!dimenB=#2%                        ** dimB = SPAN  (assumed > 0pt)
  \!countB=\!dimenB%                  ** ctB  = SPAN, as a count
%
% ** Determine index of smallest lattice point >= LOCATION
  \!dimenC=#3%                        ** dimC = LOCATION
  \advance\!dimenC -\!dimenA%         ** now dimC = LOCATION-ANCHOR
  \!countA=\!dimenC%                  ** ctA = above, as a count
  \divide\!countA \!countB%           ** now ctA = desired index, if dimC <= 0
  \ifdim\!dimenC>\!zpt
    \!dimenD=\!countA\!dimenB%        ** (tentative k)*span
    \ifdim\!dimenD<\!dimenC%          ** if this is false, ctA = desired index
      \advance\!countA 1 %            ** if true, have to add 1
    \fi
  \fi
  \!dimenC=\!countA\!dimenB%          ** lattice location = anchor + ctA*span
    \advance\!dimenC \!dimenA
  #4=\!countA%                        ** the desired index
  #5=\!dimenC%                        ** corresponding lattice location
  \ignorespaces}

% ** \!qshade [with shrinkages] [[LS,RS,BS,TS]]
% ***** during vertical shading:
% **    [the region from (xS,ybS,ytS) to] (xM,ybM,ytM) [and] (xE,ybE,ytE)
% ** Shades the region {(x,y): xS <= x <= xE, yb(x) <= y <= yt(x)}, where 
% **   yb is the quadratic thru (xS,ybS) & (xM,ybM) & (xE,ybE)
% **   yt is the quadratic thru (xS,ytS) & (xM,ybM) & (xE,ytE)
% ** xS,ybS,ytS are either given by \!startvshade or carried over
% **   as the ending values of the immediately preceding \!qshade.
% ** For the interpretation of LS, RS, BS, & TS, see \setshadesymbol. The
% **   values set there can be overridden, for the course of this \!qshade
% **   only, in the same manner as overrides are specified for
% **   \setshadesymbol.
% ***** during horizontal shading:
% **    [the region from (yS,xlS,xrS) to] (yM,xlM,xrM) [and] (yE,xlE,xrE)
\def\!qshade#1(#2,#3,#4)#5(#6,#7,#8){%
  \!xM=\!M{#2}\!!xunit
  \!ybM=\!M{#3}\!!yunit
  \!ytM=\!M{#4}\!!yunit
  \!xE=\!M{#6}\!!xunit
  \!ybE=\!M{#7}\!!yunit
  \!ytE=\!M{#8}\!!yunit
  \!getcoeffs\!xS\!ybS\!xM\!ybM\!xE\!ybE\!ybB\!ybC%**Get coefficients B & C for
  \!getcoeffs\!xS\!ytS\!xM\!ytM\!xE\!ytE\!ytB\!ytC%**y=y0 + B(x-X0) + C(x-X0)**2
  \def\!getylimits{\!qgetylimits}%
  \!shade{#1}\ignorespaces}
 
% ** \!lshade ... (xE,ybE,ytE)
% ** This is like \!qshade, but the top and bottom boundaries are linear,
% ** rather than quadratic.
\def\!lshade#1(#2,#3,#4){%
  \!xE=\!M{#2}\!!xunit
  \!ybE=\!M{#3}\!!yunit
  \!ytE=\!M{#4}\!!yunit
  \!dimenE=\!xE  \advance \!dimenE -\!xS%   ** xE-xS
  \!dimenC=\!ytE \advance \!dimenC -\!ytS%  ** ytE-ytS
  \!divide\!dimenC\!dimenE\!ytB%            ** ytB = (ytE-ytS)/(xE-xS)
  \!dimenC=\!ybE \advance \!dimenC -\!ybS%  ** ybE-ybS
  \!divide\!dimenC\!dimenE\!ybB%            ** ybB = (ybE-ybS)/(xE-xS)
  \def\!getylimits{\!lgetylimits}%
  \!shade{#1}\ignorespaces}
 
% **  \!getcoeffs{X0}{Y0}{X1}{Y1}{X2}{Y2}{B}{C}
% **  Finds  B  and  C  such that the quadratic  y = Y0 + B(x-X0) + C(x-X0)**2
% **  passes through (X1,Y1) and (X2,Y2):  when X0=0=Y0, the formulas are:
% **                   B = S1 - X1*C,   C = (S2-S1)/X2
% **  with
% **                 S1 = Y1/X1,   S2 = (Y2-Y1)/(X2-X1).
\def\!getcoeffs#1#2#3#4#5#6#7#8{% 
  \!dimenC=#4\advance \!dimenC -#2%            ** dimC=Y1-Y0
  \!dimenE=#3\advance \!dimenE -#1%            ** dimE=X1-X0
  \!divide\!dimenC\!dimenE\!dimenF%            ** dimF=S1
  \!dimenC=#6\advance \!dimenC -#4%            ** dimC=Y2-Y1
  \!dimenH=#5\advance \!dimenH -#3%            ** dimH=X2-X1
  \!divide\!dimenC\!dimenH\!dimenG%            ** dimG=S2
  \advance\!dimenG -\!dimenF%                  ** dimG=S2-S1
  \advance \!dimenH \!dimenE%                  ** dimH=X2-X0
  \!divide\!dimenG\!dimenH#8%                  ** C=(S2-S1)/(X2-X0)
  \!removept#8\!t%                             ** C, without "pt"
  #7=-\!t\!dimenE%                             ** -C*(X1-X0)
  \advance #7\!dimenF%                         ** B=S1-C*(X1-X0)
  \ignorespaces}

\def\!shade#1{%
% ** Get LS,RS,BS,TS for this panel
  \!getshrinkages#1<,,,>\!nil% %       ** now effective LS=dimE, RS=dimF,
%                                      **   BS=dimG, TS=dimH
  \advance \!dimenE \!xS%              ** now dimE=xS+LS
  \!lattice\!!xshade\!dshade\!dimenE%  ** set parity=index of left-mst x-lattice
    \!parity\!xpos%                    **   point >= xS+LS, xpos=its location
  \!dimenF=-\!dimenF%                  ** set dimF=xE-RS
    \advance\!dimenF \!xE
  \!loop\!not{\ifdim\!xpos>\!dimenF}%  ** loop over x-lattice points <= xE-RS
    \!shadecolumn%                 
    \advance\!xpos \!dshade%           ** move over to next column
    \advance\!parity 1%                ** increase index of x-point
  \repeat
  \!xS=\!xE%                           ** shift ending values to starting values
  \!ybS=\!ybE
  \!ytS=\!ytE
  \ignorespaces}

\def\!vgetshrinkages#1<#2,#3,#4,#5>#6\!nil{%
  \!override\!lshrinkage{#2}\!dimenE
  \!override\!rshrinkage{#3}\!dimenF
  \!override\!bshrinkage{#4}\!dimenG
  \!override\!tshrinkage{#5}\!dimenH
  \ignorespaces}
\def\!hgetshrinkages#1<#2,#3,#4,#5>#6\!nil{%
  \!override\!lshrinkage{#2}\!dimenG
  \!override\!rshrinkage{#3}\!dimenH
  \!override\!bshrinkage{#4}\!dimenE
  \!override\!tshrinkage{#5}\!dimenF
  \ignorespaces}

\def\!shadecolumn{%
  \!dxpos=\!xpos
  \advance\!dxpos -\!xS%            ** dx = x - xS
  \!removept\!dxpos\!dx%            ** ditto, without "pt"
  \!getylimits%                     ** get top and bottom y-values
  \advance\!ytpos -\!dimenH%        ** less TS
  \advance\!ybpos \!dimenG%         ** plus BS
  \!yloc=\!!yshade%                 ** get anchor point for this column
  \ifodd\!parity 
     \advance\!yloc \!dshade
  \fi
  \!lattice\!yloc{2\!dshade}\!ybpos%
    \!countA\!ypos%                 ** ypos=smallest y point for this column
  \!dimenA=-\!shadexorigin \advance \!dimenA \!xpos%      ** over
  \loop\!not{\ifdim\!ypos>\!ytpos}% ** loop over ypos <= yt(t)
    \!setshadelocation%             ** vmode: xloc=xpos, yloc=ypos 
%                                   ** hmode: xloc=ypos, yloc=xpos 
    \!rotateaboutpivot\!xloc\!yloc%
    \!dimenA=-\!shadexorigin \advance \!dimenA \!xloc%    ** over
    \!dimenB=-\!shadeyorigin \advance \!dimenB \!yloc%    ** up
    \kern\!dimenA \raise\!dimenB\copy\!shadesymbol \kern-\!dimenA
    \advance\!ypos 2\!dshade
  \repeat
  \ignorespaces}
 
\def\!qgetylimits{%
  \!dimenA=\!dx\!ytC              
  \advance\!dimenA \!ytB%         ** yt(t)=ytS + dx*(Bt + dx*Ct)
  \!ytpos=\!dx\!dimenA
  \advance\!ytpos \!ytS
  \!dimenA=\!dx\!ybC              
  \advance\!dimenA \!ybB%         ** yb(t)=ybS + dx*(Bb + dx*Cb)
  \!ybpos=\!dx\!dimenA
  \advance\!ybpos \!ybS}
 
\def\!lgetylimits{%
  \!ytpos=\!dx\!ytB%              ** yt(t)=ytS + dx*Bt
  \advance\!ytpos \!ytS
  \!ybpos=\!dx\!ybB%              ** yb(t)=ybS + dx*Bb
  \advance\!ybpos \!ybS}
 
\def\!vsetshadelocation{%         ** vmode: xloc=xpos, yloc=ypos 
  \!xloc=\!xpos
  \!yloc=\!ypos}
\def\!hsetshadelocation{%         ** hmode: xloc=ypos, yloc=xpos 
  \!xloc=\!ypos
  \!yloc=\!xpos}

% **************************************
% *** TICKS  (Draws ticks on graphs) ***
% **************************************

% ** User commands
% **   \ticksout
% **   \ticksin
% **   \gridlines
% **   \nogridlines
% **   \loggedticks
% **   \unloggesticks
% ** See Subsection 3.4 of the manual

% ** The following is an option of the \axis command
% **   ticks 
% **     [in] [out] 
% **     [long] [short] [length <LENGTH>] 
% **     [width <WIDTH>]
% **     [andacross] [butnotacross] 
% **     [logged] [unlogged] 
% **     [unlabeled] [numbered] [withvalues VALUE1 VALUE2 ... VALUEk / ]
% **     [quantity Q] [at LOC1 LOC2 ... LOCk / ] [from LOC1 to LOC2 by
% **       LOC_INCREMENT]
% ** See Subsection 3.2 of the manual for the rules.

% ** The various options of the  tick  field are processed by the
% ** \!nextkeyword  command defined below.
% ** For example, `\!nextkeyword short '  expands to  `\!ticksshort',
% ** while `\!nextkeyword withvalues' expands to `\!tickswithvalues'.

\def\!axisticks {%
  \def\!nextkeyword##1 {%
    \expandafter\ifx\csname !ticks##1\endcsname \relax
      \def\!next{\!fixkeyword{##1}}%
    \else
      \def\!next{\csname !ticks##1\endcsname}%
    \fi
    \!next}%
  \!axissetup
    \def\!axissetup{\relax}%
  \edef\!ticksinoutsign{\!ticksinoutSign}%
  \!ticklength=\longticklength
  \!tickwidth=\linethickness
  \!gridlinestatus
  \!setticktransform
  \!maketick
  \!tickcase=0
  \def\!LTlist{}%
  \!nextkeyword}

\def\ticksout{%
  \def\!ticksinoutSign{+}}

\ticksout

\def\nogridlines{%
  \def\!gridlinestatus{\!gridlinestoofalse}}
\nogridlines

\def\loggedticks{%
  \def\!setticktransform{\let\!ticktransform=\!logten}}
\def\unloggedticks{%
  \def\!setticktransform{\let\!ticktransform=\!donothing}}
\def\!donothing#1#2{\def#2{#1}}
\unloggedticks

% ** \!ticks/ : terminates read of tick options
\expandafter\def\csname !ticks/\endcsname{%
  \!not {\ifx \!LTlist\empty}
    \!placetickvalues
  \fi
  \def\!tickvalueslist{}%
  \def\!LTlist{}%
  \expandafter\csname !axis/\endcsname}

\def\!maketick{%
  \setbox\!boxA=\hbox{%
    \beginpicture
      \!setdimenmode
      \setcoordinatesystem point at {\!zpt} {\!zpt}   
      \linethickness=\!tickwidth
      \ifdim\!ticklength>\!zpt
        \putrule from {\!zpt} {\!zpt} to
          {\!ticksinoutsign\!tickxsign\!ticklength}
          {\!ticksinoutsign\!tickysign\!ticklength}
      \fi
      \if!gridlinestoo
        \putrule from {\!zpt} {\!zpt} to
          {-\!tickxsign\!xaxislength} {-\!tickysign\!yaxislength}
      \fi
    \endpicturesave <\!Xsave,\!Ysave>}%
    \wd\!boxA=\!zpt}
  
\def\!ticksin{%
  \def\!ticksinoutsign{-}%
  \!maketick
  \!nextkeyword}

\def\!ticksout{%
  \def\!ticksinoutsign{+}%
  \!maketick
  \!nextkeyword}

\def\!tickslength<#1> {%
  \!ticklength=#1\relax
  \!maketick
  \!nextkeyword}

\def\!tickslong{%
  \!tickslength<\longticklength> }

\def\!ticksshort{%
  \!tickslength<\shortticklength> }

\def\!tickswidth<#1> {%
  \!tickwidth=#1\relax
  \!maketick
  \!nextkeyword}

\def\!ticksandacross{%
  \!gridlinestootrue
  \!maketick
  \!nextkeyword}

\def\!ticksbutnotacross{%
  \!gridlinestoofalse
  \!maketick
  \!nextkeyword}

\def\!tickslogged{%
  \let\!ticktransform=\!logten
  \!nextkeyword}

\def\!ticksunlogged{%
  \let\!ticktransform=\!donothing
  \!nextkeyword}

\def\!ticksunlabeled{%
  \!tickcase=0
  \!nextkeyword}

\def\!ticksnumbered{%
  \!tickcase=1
  \!nextkeyword}

\def\!tickswithvalues#1/ {%
  \edef\!tickvalueslist{#1! /}%
  \!tickcase=2
  \!nextkeyword}

\def\!ticksquantity#1 {%
  \ifnum #1>1
    \!updatetickoffset
    \!countA=#1\relax
    \advance \!countA -1
    \!ticklocationincr=\!axisLength
      \divide \!ticklocationincr \!countA
    \!ticklocation=\!axisstart
    \loop \!not{\ifdim \!ticklocation>\!axisend}
      \!placetick\!ticklocation
      \ifcase\!tickcase
          \relax %  Case 0: no labels
        \or
          \relax %  Case 1: numbered -- not available here
        \or
          \expandafter\!gettickvaluefrom\!tickvalueslist
          \edef\!tickfield{{\the\!ticklocation}{\!value}}%
          \expandafter\!listaddon\expandafter{\!tickfield}\!LTlist%
      \fi
      \advance \!ticklocation \!ticklocationincr
    \repeat
  \fi
  \!nextkeyword}

\def\!ticksat#1 {%
  \!updatetickoffset
  \edef\!Loc{#1}%
  \if /\!Loc
    \def\next{\!nextkeyword}%
  \else
    \!ticksincommon
    \def\next{\!ticksat}%
  \fi
  \next}    
      
\def\!ticksfrom#1 to #2 by #3 {%
  \!updatetickoffset
  \edef\!arg{#3}%
  \expandafter\!separate\!arg\!nil
  \!scalefactor=1
  \expandafter\!countfigures\!arg/
  \edef\!arg{#1}%
  \!scaleup\!arg by\!scalefactor to\!countE
  \edef\!arg{#2}%
  \!scaleup\!arg by\!scalefactor to\!countF
  \edef\!arg{#3}%
  \!scaleup\!arg by\!scalefactor to\!countG
  \loop \!not{\ifnum\!countE>\!countF}
    \ifnum\!scalefactor=1
      \edef\!Loc{\the\!countE}%
    \else
      \!scaledown\!countE by\!scalefactor to\!Loc
    \fi
    \!ticksincommon
    \advance \!countE \!countG
  \repeat
  \!nextkeyword}

\def\!updatetickoffset{%
  \!dimenA=\!ticksinoutsign\!ticklength
  \ifdim \!dimenA>\!offset
    \!offset=\!dimenA
  \fi}

\def\!placetick#1{%
  \if!xswitch
    \!xpos=#1\relax
    \!ypos=\!axisylevel
  \else
    \!xpos=\!axisxlevel
    \!ypos=#1\relax
  \fi
  \advance\!xpos \!Xsave
  \advance\!ypos \!Ysave
  \kern\!xpos\raise\!ypos\copy\!boxA\kern-\!xpos
  \ignorespaces}

\def\!gettickvaluefrom#1 #2 /{%
  \edef\!value{#1}%
  \edef\!tickvalueslist{#2 /}%
  \ifx \!tickvalueslist\!endtickvaluelist
    \!tickcase=0
  \fi}
\def\!endtickvaluelist{! /}

\def\!ticksincommon{%
  \!ticktransform\!Loc\!t
  \!ticklocation=\!t\!!unit
  \advance\!ticklocation -\!!origin
  \!placetick\!ticklocation
  \ifcase\!tickcase
    \relax % Case 0: no labels
  \or %      Case 1: numbered
    \ifdim\!ticklocation<-\!!origin
      \edef\!Loc{$\!Loc$}%
    \fi
    \edef\!tickfield{{\the\!ticklocation}{\!Loc}}%
    \expandafter\!listaddon\expandafter{\!tickfield}\!LTlist%
  \or %      Case 2: labeled
    \expandafter\!gettickvaluefrom\!tickvalueslist
    \edef\!tickfield{{\the\!ticklocation}{\!value}}%
    \expandafter\!listaddon\expandafter{\!tickfield}\!LTlist%
  \fi}

\def\!separate#1\!nil{%
  \!ifnextchar{-}{\!!separate}{\!!!separate}#1\!nil}
\def\!!separate-#1\!nil{%
  \def\!sign{-}%
  \!!!!separate#1..\!nil}
\def\!!!separate#1\!nil{%
  \def\!sign{+}%
  \!!!!separate#1..\!nil}
\def\!!!!separate#1.#2.#3\!nil{%
  \def\!arg{#1}%
  \ifx\!arg\!empty
    \!countA=0
  \else
    \!countA=\!arg
  \fi
  \def\!arg{#2}%
  \ifx\!arg\!empty
    \!countB=0
  \else
    \!countB=\!arg
  \fi}
 
\def\!countfigures#1{%
  \if #1/%
    \def\!next{\ignorespaces}%
  \else
    \multiply\!scalefactor 10
    \def\!next{\!countfigures}%
  \fi
  \!next}

\def\!scaleup#1by#2to#3{%
  \expandafter\!separate#1\!nil
  \multiply\!countA #2\relax
  \advance\!countA \!countB
  \if -\!sign
    \!countA=-\!countA
  \fi
  #3=\!countA
  \ignorespaces}

\def\!scaledown#1by#2to#3{%
  \!countA=#1\relax%                          ** get original #
  \ifnum \!countA<0 %                         ** take abs value,
    \def\!sign{-}%                            **   remember sign
    \!countA=-\!countA
  \else
    \def\!sign{}%
  \fi
  \!countB=\!countA%                          ** copy |#|
  \divide\!countB #2\relax%                   ** integer part (|#|/sf)
  \!countC=\!countB%                          ** get sf * (|#|/sf)
    \multiply\!countC #2\relax
  \advance \!countA -\!countC%                ** ctA is now remainder
  \edef#3{\!sign\the\!countB.}%               ** +- integerpart.
  \!countC=\!countA %                         ** Tack on proper number
  \ifnum\!countC=0 %                          **   of zeros after .
    \!countC=1
  \fi
  \multiply\!countC 10
  \!loop \ifnum #2>\!countC
    \edef#3{#3\!zero}%
    \multiply\!countC 10
  \repeat
  \edef#3{#3\the\!countA}%                    ** Add on rest of remainder
  \ignorespaces}

\def\!placetickvalues{%
  \advance\!offset \tickstovaluesleading
  \if!xswitch
    \setbox\!boxA=\hbox{%
      \def\\##1##2{%
        \!dimenput {##2} [B] (##1,\!axisylevel)}%
      \beginpicture 
        \!LTlist
      \endpicturesave <\!Xsave,\!Ysave>}%
    \!dimenA=\!axisylevel
      \advance\!dimenA -\!Ysave
      \advance\!dimenA \!tickysign\!offset
      \if -\!tickysign
        \advance\!dimenA -\ht\!boxA
      \else
        \advance\!dimenA  \dp\!boxA
      \fi
    \advance\!offset \ht\!boxA 
      \advance\!offset \dp\!boxA
    \!dimenput {\box\!boxA} [Bl] <\!Xsave,\!Ysave> (\!zpt,\!dimenA)
  \else
    \setbox\!boxA=\hbox{%
      \def\\##1##2{%
        \!dimenput {##2} [r] (\!axisxlevel,##1)}%
      \beginpicture 
        \!LTlist
      \endpicturesave <\!Xsave,\!Ysave>}%
    \!dimenA=\!axisxlevel
      \advance\!dimenA -\!Xsave
      \advance\!dimenA \!tickxsign\!offset
      \if -\!tickxsign
        \advance\!dimenA -\wd\!boxA
      \fi
    \advance\!offset \wd\!boxA
    \!dimenput {\box\!boxA} [Bl] <\!Xsave,\!Ysave> (\!dimenA,\!zpt)
  \fi}

\normalgraphs
\catcode`!=12 %  *****  THIS MUST NEVER BE OMITTED

\def\FOOTNOTE{\footnote*}

\author{Daniel Allcock\FOOTNOTE{Supported by National Science
Foundation Graduate and Postdoctoral Fellowships.}}

\title{New complex- and quaternion-hyperbolic reflection groups}

\date{1 June 1999}
%\date{30 May 1999}
%\date{9 September 1997}
%
% Commenced 20 February 1996

\subject{22E40 (11F06, 11E39, 53C35)}

\note{Dedicated to my father John Allcock, 1940--1991}

% small enough font to dodge overfull hbox
\font\titlefont=cmb10 scaled \magstep2
\plaintitlepage
%

%dashes
%
\def\emdash{---}

%
%----------Options-(comment-or-uncomment-as-desired)-------------
%
\magnification=1000
%
%    for doublespacing
%\baselineskip=20pt
%
\symbolictags

%    eliminates the nasty overfullbox bars:
\def\hideboxes{\overfullrule=0pt}
\hideboxes 

%-----------------beginning-of-document---------------------------
%

% Definitions for cross-references
%
% Do not tamper with the following line.
% ---begin deftags---
\deftag{1}{sec-introduction}
\deftag{2}{sec-lattice-background}
\deftag{2.1}{eq-def-of-reflection}
\deftag{3}{sec-lattice-reference}
\deftag{3.1}{thm-barnes-wall-properties}
\deftag{4}{sec-hyperbolic-space}
\deftag{4.1}{eq-half-space-coordinates}
\deftag{4.2}{eq-def-of-translation}
\deftag{4.3}{eq-product-of-translations}
\deftag{4.4}{eq-inverse-of-translation}
\deftag{4.5}{eq-commutator-of-translations}
\deftag{4.6}{eq-conjugate-of-translation}
\deftag{5}{sec-reflections}
\deftag{5.1}{thm-stabilizer-of-rho}
\deftag{5.1}{eq-got-translation}
\deftag{5.2}{thm-gaussian-long-roots-height-reduction}
\deftag{5.2}{eq-coordinates-of-r}
\deftag{5.3}{eq-bar-1}
\deftag{5.4}{eq-bar-2}
\deftag{5.5}{eq-bar-3}
\deftag{5.6}{eq-bar-4}
\deftag{5.7}{eq-bar-5}
\deftag{5.3}{thm-short-root-height-reduction}
\deftag{5.8}{eq-coordinates-of-r-2}
\deftag{5.9}{eq-ip-with-r'-2}
\deftag{5.10}{eq-height-of-v'-2}
\deftag{5.11}{eq-foo3}
\deftag{5.12}{eq-relation-between-xi's}
\deftag{5.13}{eq-foo4}
\deftag{5.14}{eq-relation-between-xi's-2}
\deftag{5.1}{tab-table-of-heights}
\deftag{6}{sec-reflection-groups}
\deftag{6.1}{thm-cov-radius-implies-reflective}
\deftag{6.2}{thm-gaussian-lattices-reflective}
\deftag{6.3}{thm-examples-of-reflective-lattices}
\deftag{6.4}{thm-basics-for-eisenstein-lattices}
\deftag{6.5}{thm-braiding-reflections}
\deftag{6.6}{thm-eI2,1-eI3,1-eI4,1}
\deftag{6.7}{thm-basics-for-Hurwitz-lattices}
\deftag{6.8}{thm-qI2,1-qI3,1}
\deftag{6.9}{thm-linear-independence-descends}
\deftag{6.10}{thm-coxeter-todd-well-covered}
\deftag{6.11}{thm-eI7,1}
\deftag{6.12}{thm-barnes-wall-well-covered}
\deftag{6.13}{thm-leech-holes-in-BW}
\deftag{6.1}{eq-dimension-counting}
\deftag{6.14}{thm-qI5,1}
\deftag{7}{sec-selfdual-classification}
\deftag{7.1}{thm-indefinite-selfdual-classification}
\deftag{7.2}{thm-classification-of-selfdual-lattices}
\deftag{8}{sec-comparison}
\deftag{8.1}{thm-classify-honest-reflections}
\deftag{8.2}{thm-projective-reflections-are-honest}
\deftag{8.3}{thm-cone-angles}
\deftag{8.4}{thm-my-groups-are-new}
% ---end deftags---
% Do not tamper with the previous line or the following line.
% ---begin defcites---
\defcite{1}{allcock:ch13}
\defcite{2}{allcock:ch4-cubic-moduli}
\defcite{3}{allcock:ch4-cubic-moduli-announcement}
\defcite{4}{allcock:oh2}
\defcite{5}{bachoc:selfdual-Hurwitz-lattices}
\defcite{6}{bachoc:8-dimensional-Hurwitz-lattices}
\defcite{7}{reb:lzl}
\defcite{8}{reb:like-leech}
\defcite{9}{borel:arithmetic-groups}
\defcite{10}{jhc:26dim}
\defcite{11}{ATLAS}
\defcite{12}{jhc:leech-radius}
\defcite{13}{jhc:laminated-lattices}
\defcite{14}{jhc:Coxeter-Todd}
\defcite{15}{splag}
\defcite{16}{coxeter:quotients-of-braid-groups}
\defcite{17}{cox:gens-and-rlns}
\defcite{18}{deligne:monodromy-of-hypergeometric-functions}
\defcite{19}{feit:unimodular-Eisenstein}
\defcite{20}{milnor:bilinear-forms}
\defcite{21}{mostow:remarkable-polyhedra}
\defcite{22}{mostow:picard-lattices-from-half-inegral-conditions}
\defcite{23}{mostow:monodromy-groups-on-the-complex-n-ball}
\defcite{24}{thurston:shapes-of-polyhedra}
\defcite{25}{vin:groups-of-quad-forms}
\defcite{26}{vin:unimodular-quad-forms}
\defcite{27}{vin:I18.1-and-I19.1}
% ---end defcites---
% Do not tamper with the preceding line.

\def\w{\omega}
\def\wbar{\bar{\w}}

\def\quaternionic{{\bbb H}}
\def\K{{\bbb K\kern.05em}}	
\def\O{{\bbb O}}
\def\eisen {\cale} % Eisenstein integers
\def\gauss {\calg} % Gaussian integers
\def\Hurwitz{\calh}
\def\ring{\calr}

\def\H{\quaternionic}

\def\reflec{\mathop{\rm Reflec}\nolimits} % group generated by reflections 
\def\shr{\mathop{\reflec}_0\nolimits} % group generated by
					 % reflections in short
					 % roots 
\def\paut{P\!\mathop{\rm Aut}\nolimits}

\def\GL{\mathop{\rm GL}\nolimits}
\def\SL{\mathop{\rm SL}\nolimits}
\def\height{{\rm ht\,}}
\def\diag{\hbox{diag}\,}
\def\spanof#1{\langle#1\rangle}
\def\cp{\C P^}
\def\kp{\K P^}
\def\cpnm{\C P^{n-1}}
\def\kpnp{\K P^{n+1}}
\def\khnp{\K H^{n+1}}
\def\dkhnp{\partial\khnp}

\def\hyperbolic{H}
\def\rh {\R\hyperbolic^}
\def\ch {\C\hyperbolic^}	
\def\qh{\quaternionic\hyperbolic^} 

\def\qhn{\qh n}
\def\chn{\ch n}

\def\rhnp {\rh{n+1}}

\def\I#1,#2{I_{#1,#2}}
\def\II#1,#2{I\negspace I_{#1,#2}}
\def\oh{\O\hyperbolic^}

\def\trianglegroup#1#2#3{(#1,#2,#3)}

\def \pu #1,#2{PU(#1,#2)}
\def\ip#1#2{\langle#1{|}#2\rangle}
\def\Eisenstein {\eisen}
\def\Gaussian{\gauss}
\def\cct{\Lambda_6^\Eisenstein} % complex Coxeter-Todd lattice
 % real form of Coxeter-Todd lattice 
\def\rbwl{BW_{16}} % real form of barnes-wall lattice
\def\qbwl{\Lambda_4^\Hurwitz} % quaternionic Barnes-Wall lattice 
\def\leech{\Lambda_{24}}
 % quaternionic leech lattice
\def\geeight{E_8^\Gaussian}	% Gaussian E8 lattice 
\def\qeeight{E_8^\Hurwitz}	%quaternionic E8 lattice

\def\eI#1,#2{\I#1,#2^{\eisen}}   % Lorentzian lattices of
				 % various types
\def\eII#1,#2{\II#1,#2^{\,\eisen}}
\def\qI#1,#2{\I#1,#2^{\,\Hurwitz}}
\def\qII#1,#2{\II#1,#2^{\,\Hurwitz}}
\def\gI#1,#2{\I#1,#2^{\,\gauss}}
\def\gII#1,#2{\II#1,#2^{\,\gauss}}
\def\zI#1,#2{\I#1,#2^{\,\Z}}
\def\rI#1,#2{\I#1,#2^{\,\ring}}
\def\rII#1,#2{\II#1,#2^{\,\ring}}
\def\zI#1,#2{\I#1,#2^{\,\Z}}
\def\zII#1,#2{\II#1,#2^{\,\Z}}

\def\vec#1#2{(#1_1,\ldots,#1_{#2})}
\def\vector#1#2{\vec{#1}{#2}}

% imaginary constants
\def\i{{\rm i}}
\def\j{{\rm j}}
\def\k{{\rm k}}

\def\Cbar{\bar C}
\def\G{\Gamma}
\def\negspace{\!}

\abstract
We consider the automorphism groups of various Lorentzian
lattices over the Eisenstein, Gaussian, and Hurwitz integers,
and in some of them we find reflection groups of finite index.
These provide  new finite-covolume
reflection groups acting on complex and quaternionic hyperbolic
spaces.  Specifically, we provide groups
acting on $\chn$ for all $n<6$ and $n=7$, and on $\qhn$ for
$n=1,2,3,$ and $5$. We compare our groups to those discovered
by Deligne and Mostow and show that our largest 
examples are new.  For many of these Lorentzian
lattices we show that the entire symmetry group is generated by
reflections, and obtain a description of the group in terms of
the combinatorics of a lower-dimensional positive-definite
lattice.  The techniques needed for our lower-dimensional
examples are elementary, but to construct our best examples we
also need certain facts about the Leech lattice.  We give a
new and geometric proof of the classifications of selfdual
Eisenstein lattices of dimension $\leq6$ and of selfdual Hurwitz
lattices of dimension $\leq4$.

\section{\Tag{sec-introduction}. Introduction}

In this paper we carry out complex and quaternionic analogues of
some of Vinberg's extensive study of reflection groups on real
hyperbolic space. In \cite{vin:groups-of-quad-forms} and
\cite{vin:unimodular-quad-forms} he investigated the 
symmetry groups of the integral quadratic forms
$\diag[-1,+1,\ldots,+1]$, or equivalently the Lorentzian
lattices $\I n,1$. He was able to describe these groups very
concretely for $n\leq17$, and extensions of his work by Vinberg
and Kaplinskaja \cite{vin:I18.1-and-I19.1} and Borcherds
\cite{reb:lzl} provide similar
descriptions for all $n\leq23$. In particular, the subgroup of
$\aut\I n,1$ generated by reflections has finite index just when
$n\leq19$.

In this paper, we study the symmetry groups of Lorentzian lattices
over the rings $\Gaussian$ and $\Eisenstein$ of Gaussian and
Eisenstein integers and the ring $\Hurwitz$ of Hurwitz integers
(a discrete subring of the skew field $\H$ of quaternions). Most
of the paper is devoted to the most natural of such lattices,
the selfdual ones. The symmetry groups of these lattices provide
a large number of discrete groups generated by reflections and
acting with finite-volume quotient on  the hyperbolic spaces
$\ch n$ and $\qh n$. We construct a total of 19 such groups,
including groups acting on $\ch7$ and $\qh5$. 
At least one of our groups has been discovered before, in the
work of Deligne and Mostow
\cite{deligne:monodromy-of-hypergeometric-functions}, Mostow
\cite{mostow:picard-lattices-from-half-inegral-conditions} and
Thurston \cite{thurston:shapes-of-polyhedra}, but our largest
examples are new.  To the author's knowledge,
quaternion-hyperbolic reflection groups not been studied before.

Our results and techniques have found important
application in work of the author, J. Carlson and D. Toledo
\cite{allcock:ch4-cubic-moduli}, 
\cite{allcock:ch4-cubic-moduli-announcement} 
on the moduli space of complex cubic surfaces. Namely, this
space is isomorphic to the Satake compactification of the
quotient of $\ch4$ by one of the reflection groups studied
here. Furthermore, the moduli space of ``marked'' cubic surfaces
may be realized as the Satake compactification of the quotient
of $\ch4$ by a congruence subgroup, which is also a reflection
group in its own right.
\goodbreak

The techniques used by Vinberg and others for the real
hyperbolic case rely heavily on the fact that if a discrete
group $G$ is generated by reflections of $\rh n$, then the
mirrors of the reflections of $G$ chop $\rh n$ into pieces and each
piece may be taken as a fundamental domain for $G$. Work with
complex or quaternionic reflection groups is much more
complicated, since hyperplanes have real codimension 2 or 4, and
so the mirrors  
fail to chop hyperbolic space into pieces. Our solution to
this problem is to avoid fundamental domains altogether.
Each of our groups is defined as the
subgroup $\reflec L$ of $\aut L$ generated by reflections, where
$L$ is a Lorentzian lattice over $\Gaussian$, $\Eisenstein$ or
$\Hurwitz$. (A Lorentzian lattice is a free module equipped with
a Hermitian form of signature $-+\cdots+$.) 
Since $\aut L$ is an arithmetic group,
to show that the quotient of $\ch n$ or $\qh n$ by $\reflec L$ has finite
volume it suffices to show that $\reflec L$ has finite index in
$\aut L$. In this case we say that $L$ is reflective. Our basic
strategy for proving a lattice $L$ to be reflective is to prove first that
$\reflec L$ acts with only finitely many orbits on the vectors
of $L$ of norm $0$, and second that the stabilizer in $\reflec
L$ of one such vector has finite index in the stabilizer in
$\aut L$. We work mostly arithmetically, avoiding use of
such tools as the bisectors introduced by Mostow for his study
\cite{mostow:remarkable-polyhedra} of 
reflection groups on $\ch2$.

However, there are certain steps in our constructions where
geometric ideas play a key role. We express each of our
Lorentzian lattices $L$ in the form $\Lambda\oplus\II1,1$, where
$\Lambda$ is positive-definite and $\II1,1$ is a certain
2-dimensional lattice, the ``hyperbolic plane'', with inner
product matrix $\smallmatrix0110$. It turns out
that this description of $L$ allows one to easily write down a
large collection of reflections of $L$, parameterized by (a
central extension of) the lattice $\Lambda$. It turns out that
if $\Lambda$  has enough vectors of norms~1 and~2, and provides a
good covering of  Euclidean space by balls, then one can
automatically deduce that $L$ 
is reflective. This implication is the content of
Thm.~\tag{thm-cov-radius-implies-reflective}. The rest of
Section~\tag{sec-reflection-groups} is devoted to the
application of this theorem and related ideas in the
study of various examples. In particular, we prove that each of
the selfdual Lorentzian lattices
$$
\vbox{\halign{$#$\hfil&\qquad$#$\hfil\cr
\eI n,1&n=1,2,3,4,7,\cr
\noalign{\smallskip}\cr
\gII n,1&n=1,5,\cr
\noalign{\smallskip}\cr
\qI n,1&n=1,2,3,5\cr}}
$$
is reflective. (These lattices are defined in
Section~\tag{sec-lattice-reference} and characterized in
Thm.~\tag{thm-indefinite-selfdual-classification}.) For some
of these, we obtain more detailed information. In particular, we
prove that $\reflec\eI n,1=\aut\eI n,1$ for $n=2$, $3$, $4$ or
$7$ and that $\reflec\qI n,1$ has index at most $4$ in $\aut\qI
n,1$ for $n=2$, $3$ or 
$5$. We also give explicit
descriptions of the reflection groups of $\rI1,1$ and
$\gII1,1$ as subgroups of certain Coxeter groups, acting on
$\ch1\isomorphism \rh2$ and  $\qh1\isomorphism\rh4$. 

We note that the geometric ideas used here, namely that
good coverings of Euclidean space lead to hyperbolic reflection
groups, apply even when $\C$ or $\H$ is replaced by the
nonassociative field $\O$ of octaves (or octonions or Cayley
numbers). In \cite{allcock:oh2} we constructed two octave
reflection groups acting on $\oh2$ and one acting on
$\oh1\isomorphism \rh8$, and interpreted these groups as the
stabilizers of `lattices' over a certain discrete subring of
$\O$.

We provide background information on lattices in
Section~\tag{sec-lattice-background} and examples of them in
Section~\tag{sec-lattice-reference}; the latter should be
referred to only as needed. Section~\tag{sec-hyperbolic-space}
establishes our conventions regarding hyperbolic geometry. In
Section~\tag{sec-reflections} we relate certain geometric
properties of a positive-definite lattice $\Lambda$ to the
reflection group of $\Lambda\oplus\II1,1$ and lay the
foundations for Section~\tag{sec-reflection-groups}, where we
construct all of our examples. 
In Section~\tag{sec-selfdual-classification} we
explain the correspondence between primitive isotropic
sublattices of $\I n+1,1$ and positive-definite selfdual
lattices of dimension $n$. We use this correspondence to provide
a quick 
geometric proof of the classification of selfdual
lattices over $\Eisenstein$ and $\Hurwitz$ in dimensions
$\leq 6$ and $\leq 4$, respectively. The only examples besides the 
lattices $\Eisenstein^n$ and $\Hurwitz^n$ 
are the Coxeter-Todd lattice $\cct$ and a quaternionic
form $\qbwl$ of the Barnes-Wall lattice.  
In Section~\tag{sec-comparison} we show that our largest three
groups,
namely $\reflec L$ for $L=\eI7,1$, $\eI4,1$ and $\qII5,1$,  
are not among the 94 groups constructed in
\cite{deligne:monodromy-of-hypergeometric-functions},
\cite{mostow:picard-lattices-from-half-inegral-conditions} and
\cite{thurston:shapes-of-polyhedra}. We also sketch a proof that
$\reflec\eI3,1$ does appear among these groups.

The easiest route  to a new reflection group
is our study of $\gII5,1=\geeight\oplus\gII1,1$, which acts
on $\ch5$. This requires only Lemmas~\tag{thm-stabilizer-of-rho} and~\tag{thm-gaussian-long-roots-height-reduction} and
Thm.~\tag{thm-gaussian-lattices-reflective}. Most of our other examples require the more
complicated Lemma~\tag{thm-short-root-height-reduction} in place of~\tag{thm-gaussian-long-roots-height-reduction}. The arguments for
$\eI7,1$ and $\qI5,1$ also require fairly involved
space-covering arguments, involving embeddings of the
Coxeter-Todd and Barnes-Wall lattices into the famous Leech
lattice $\leech$. It is pleasing that $\leech$ plays a role
here, because our basic approach was inspired by Conway's
elegant description \cite{jhc:26dim} of the isometry group of the
$\Z$-lattice $\II25,1=\leech\oplus\II1,1$ in terms of the
combinatorics of $\leech$. The embeddings of the Coxeter-Todd
and Barnes-Wall lattices into $\leech$ have also been used by
Borcherds \cite{reb:like-leech} to produce interesting {\it real} hyperbolic
reflection groups, acting on $\rh{13}$ and $\rh{17}$. Finally,
the Leech lattice plays a much more direct role in 
\cite{allcock:ch13}, which constructs several other complex and quaternionic
hyperbolic reflection groups, including one on $\ch{13}$ and one
on $\qh7$.

Most of this paper is derived from the author's Ph.D. thesis at
Berkeley; he would like to thank his dissertation advisor,
R. Borcherds, for his interest and suggestions\emdash in particular
for suggesting that the quaternionic Barnes-Wall lattice would
provide a reflection group on $\qh5$.

\section{\Tag{sec-lattice-background}. Lattices}

We denote by $\ring$ any one of the rings $\gauss$, $\eisen$,
and $\Hurwitz$\emdash the Eisenstein, Gaussian, and Hurwitz
integers. That is, $\gauss=\Z[\i ]$ and $\eisen=\Z[\w]$, where
$\w=(-1+\sqrt{-3})/2$ is a primitive cube root of unity.  The
ring $\Hurwitz$ is the integral span of its $24$ units $\pm1$,
$\pm \i $, $\pm \j $, $\pm \k $ and $(\pm1\pm \i \pm \j \pm \k )/2$ in the
skew field $\H$ of quaternions. We write $\K$ for the field
($\C$ or $\H\,$) naturally containing $\ring$.  Conjugation
$x\mapsto\bar x$ denotes complex or quaternionic conjugation, as
appropriate.  For any element $x$ of $\K$, we write $\re
x=(x+\bar x)/2$ and $\im x=(x-\bar x)/2$ for the real and
imaginary parts of $x$, and say that $x$ is imaginary if $\re
x=0$. If $X\sset\K$ then we write $\im X$ for the set of  imaginary
elements of $X$. For any $x\in\K$, $x\bar x$ is a positive real
number, and the absolute value $|x|$ of
$x$ is defined to be $(x\bar x)^{1/2}$.  It is convenient
to define the element $\theta=\w-\bar\w=\sqrt{-3}$ of
$\eisen$. We will sometimes also consider $\w$ and $\theta$ as
elements of $\Hurwitz$, via the embedding
$\eisen\rightarrow\Hurwitz$ defined  by
$\w\mapsto(-1+\i +\j +\k )/2$ or equally well by $\theta\mapsto \i +\j +\k $.
 
A {lattice} $\Lambda$ over $\ring$ is a free (right) module over
$\ring$ equipped with a Hermitian form, which is to say a
$\Z$-bilinear pairing (the inner product) $\ip{\cdot}{\cdot}:
\Lambda\times \Lambda\rightarrow \K$ such that
$$
\ip{x}{y}=\overline{\ip{y}{x}}\qquad\hbox{and}\qquad
\ip{x}{y\a}=\ip{x}{y}\a
$$
for all $x,y\in \Lambda$ and 
$\a\in\ring$. 
We use right modules and right-linear Hermitian forms so that
lattice automorphisms can be described by matrices acting on the
left. 
A Hermitian form on a (right) vector space
over $\K$ is defined similarly.
Section~\tag{sec-lattice-reference} defines a number of
interesting lattices and lists some of their
properties. Sometimes we indicate that a lattice $\Lambda$ is an
$\ring$-lattice by writing $\Lambda^\ring$ or somesuch.

If $S\sset \Lambda$ then we denote by $S^\perp$ its orthogonal
complement: those elements of $\Lambda$ whose inner products
with all elements of $S$ vanish.  We say that $\Lambda$ is nonsingular if
$\Lambda^\perp=\{0\}$ and that $\Lambda$ is
integral if for all $x,y\in \Lambda$, the inner product
$\ip{x}{y}$ lies in $\ring$.  All lattices we consider will
be integral and nonsingular unless otherwise specified.  The
dual $\Lambda^*$ of $\Lambda$ is the set of all $\ring$-linear maps
from $\Lambda$ to $\ring$.  An integral lattice $\Lambda$ is
called {selfdual} if the natural map from $\Lambda$ to
$\Lambda^*$ is onto. A selfdual lattice is sometimes called
`unimodular', because the matrix of inner products of any
basis for $\Lambda$ has determinant $\pm1$; we use `selfdual' to avoid
discussing determinants of quaternionic matrices.

The {norm} of a vector $v\in V$ is defined to be
$v^2=\ip{v}{v}$; some authors call this the squared norm of $v$,
but our convention is better for indefinite forms.  We say
that $v$ is isotropic, or null, if $v^2=0$. A lattice is
isotropic, or null, if each of its elements is. A lattice is
called even if each of its elements has even norm and odd
otherwise. A sublattice $\Lambda'$ of $\Lambda$ is called
primitive if $\Lambda'=\Lambda\cap(\Lambda'\tensor\Q)$. A vector
$v$ of $\Lambda$ is called primitive if $v=w\a$ for $w\in
\Lambda$ and $\a\in\ring$ implies that $\a$ is a unit. Because
the rings $\Gaussian$, $\Eisenstein$ and $\Hurwitz$ are
principal ideal domains, a nonzero vector is primitive if and only if its
$\ring$-span is primitive as a sublattice. We will sometimes
write $\spanof{v}$ for the $\ring$-span of $v\in \Lambda$. 

We sometimes define an $\ring$-lattice by describing a Hermitian
form on $\ring^n$. We do this by giving an $n\times n$ matrix
$(\phi_{ij})$ with entries in $\ring$ such that
$\overline{\phi_{ij}}=\phi_{ji}$. Then the Hermitian form is
given by
$$
\ip{\vector{x}{n}}{\vector{y}{n}}=\sum_{i,j=1}^{n}\bar{x}_i\phi_{ij}y_j\;.
$$
We may also view a lattice as a subset of the vector space
$V=\Lambda\tensor\R$ over the field $\K=\R\tensor\ring$. The
Hermitian form on $\Lambda$ gives rise to one on $V$.  If
$\Lambda$ is nonsingular then $\Lambda^*$ may be identified with
the set of vectors in $V$ having $\ring$-integral inner product
with each element of $\Lambda$.
Every nonsingular Hermitian form on a vector space $V$ over $\K$
is equivalent under $\GL(V)$ to one given by a diagonal matrix,
with each diagonal entry being $\pm1$.  
The signature of $\Phi$ is the ordered pair $(n,m)$ where $n$
(resp. $m$) is the number of $+1$'s (resp. $-1$'s). This
characterizes $\Phi$ up to equivalence under $\GL(V)$. We write
$\K^{n,m}$ for the vector space $\K^{n+m}$ equipped with the
standard Hermitian form of signature $(n,m)$; the isometry group
of $\K^{n,m}$ is the unitary group $U(n,m;\K)$.
The term ``Lorentzian'' is applied to various
concepts in the study of real Minkowski space $\R^{n,1}$. By
analogy with this we call a lattice Lorentzian
if its signature is $(n,1)$.
Any isotropic sublattice of a Lorentzian lattice has
dimension $\le1$. 

If $\Lambda$ is positive-definite then $\Lambda\tensor\R$ is a
copy of Euclidean space under the metric
$d(x,y)=\sqrt{(x-y)^2}$. Points of $\Lambda\tensor\R$ at maximal
distance from $\Lambda$ are called deep holes of $\Lambda$. The
maximal distance is called the covering radius of $\Lambda$,
because closed balls of that radius placed at lattice points
exactly cover $\Lambda\tensor\R$. The lattice points nearest a
deep hole are called the vertices of the hole.  The covering
radii of the $\Z$-lattices $\im\Gaussian$, $\im\Eisenstein$ and
$\im\Hurwitz$ are $1/2$, $\sqrt3/2$ and $\sqrt3/2$,
respectively. The first two are obvious and the last follows
because $\im\Hurwitz$ is the 3-dimensional lattice spanned
by $\i $, $\j $ and $\k $.  Any two deep holes of $\im\ring$ are
equivalent under translation by some element of $\im\ring$.

Suppose that $V$ is a Hermitian vector space over $\K$, $\xi\in\K$ is a root of
unity and $v\in V$ has nonzero norm.  We define the
$\xi$-reflection in $v$ to be the map
$$
v\mapsto v- r(1-\xi){\ip{r}{v}\over r^2}\;.
\eqno\eqTag{eq-def-of-reflection}
$$
This is an isometry of the right vector space $V$ which fixes $r^\perp$
pointwise and carries $r$ to $r\xi$.  ({\it Warning:\/} if
$\K=\H$ then although the reflection acts by right scalar
multiplication on $r$, it does not act this way on  the
entire $\H$-span of $r$. This is due to the noncommutativity of
multiplication in $\H$.)  Unless otherwise specified, we will use the term
``reflection'' to mean ``reflection in a vector of positive
norm''. Under the conventions of
Section~\tag{sec-hyperbolic-space}, 
$(-1)$-reflections in negative norm vectors act on hyperbolic space as
inversions in points, rather than by reflections in
hyperplanes. This is why we focus on positive-norm vectors.
We call $r^\perp$ the 
mirror of the reflection.  Reflections of order 2, 3,$\,\ldots$
are sometimes called biflections, triflections, etc.  A
$\xi$-reflection is a biflection just if $\xi=-1$; in this case
we recover the classical notion of a reflection.

Suppose $L$ is an integral lattice.  If $v\in L$ has norm $1$
(resp. $2$) then we say that $v$ is a short (resp. long) root of
$L$.  Inspection of Eq.~\eqtag{eq-def-of-reflection} reveals that if
$\xi$ is a unit of $\ring$ then $\xi$-reflection in any short
root of $L$ preserves $L$.  Furthermore, biflections in long
roots of $L$ also preserve $L$.  We define the reflection group
$\reflec L$ to be the subgroup of $\aut L$ generated by
reflections (in positive-norm vectors), and we say that $L$ is
reflective if $\reflec L$ 
has finite index in $\aut L$.  In general, a group generated by
reflections is called a reflection group. Since $\aut L$ is an
arithmetic subgroup of the semisimple real Lie group
$U(L\tensor\R\,;\K)$, a theorem of Borel and Harish-Chandra
\cite{borel:arithmetic-groups} implies that it has finite
covolume. It follows that $L$ is reflective if and only if
$\reflec L$ also has finite covolume. It may happen that
$\reflec L$ contains reflections other than those in its roots,
but we will not use them.  

\section{\Tag{sec-lattice-reference}. Reference: examples of lattices}
 
This section contains background information on the various
complex and
quaternionic lattices we will use; it should be
referred to only as necessary. We briefly define each
lattice, list a few important
properties and give references to the literature. The main
source is \ecite{splag}{Chap.~4}.  All lattices described
here are integral.  When lattices are described as subsets of
$\K^n$ it should be understood that  the Hermitian form is
$\ip{\vector{x}{n}}{\vector{y}{n}}=\sum\bar{x}_iy_i$.

The simplest lattice is $\ring^n$, which is obviously selfdual.
Its symmetry group contains the left-multiplica\-tion by each
diagonal matrix all of whose diagonal entries are units of
$\ring$. It is easy to see that the group of these coincides with the
group generated by the reflections in the short roots. Adjoining
to this group the permutations of coordinates, which are
generated by biflections in long roots such as
$(1,-1,0,\ldots,0)$, we see that $\aut\ring^n$ is a reflection
group.

If $\Lambda$ is a lattice then its {real form} is the
$\Z$-module $\Lambda$ equipped with the inner product
$(x,y)=\re\ip{x}{y}$.  Here are three forms of the $E_8$ root
lattice:

$$\displaylines{ E_8={1\over2}\set{\vec{x}{8}\in\Z^8 }{
x_i\equiv x_j \pmod{2},\;\;
\sum x_i\congruent 0 \pmod{4}}\;, \cr
E_8^\gauss = {1\over1+\i }\set{\vec{x}{4}\in\gauss^4 }{ x_i\equiv
x_j
\pmod{1+\i },\;\; \sum x_i\congruent 0\pmod{2}}, \cr
\qeeight = \set{(x_1,x_2)\in\Hurwitz^2 }{
x_1+x_2\congruent 0\pmod{1+\i}}.
\cr}$$
It is straightforward to identify the real forms of these
lattices with each other; each has covering radius 1 and
minimal norm 2.  Often the dimension of a lattice is indicated
by a subscript. Unfortunately, this sometimes
refers to its dimension as a $\Z$-lattice and sometimes to its
dimension as an
$\ring$-lattice. There seems to be no universal solution to this
notational problem.

Another set of useful even Gaussian lattices are
$$
D_{2n}^\gauss=\set{\vector{x}{n}\in\gauss^n}{
\sum x_i\congruent 0\pmod{1+\i }}\;, 
$$
whose real forms are the $D_{2n}$
root lattices.
The $D_4$ lattice is also the real form of $\Hurwitz$,
scaled up by a factor of $2^{1/2}$. The covering radius of $D_{2n}$
is 
$(n/2)^{1/2}$. 

The Eisenstein lattice 
$$ 
D_3(\theta)=\set{(x,y,z)\in\Eisenstein^3}{x+y+z\congruent0\pmod\theta}
$$
is one of the lattices $D_n(\sqrt{-3})$  introduced by 
Feit in \cite{feit:unimodular-Eisenstein}.
It has $54$ long roots and $72$ vectors  of norm
$3$; biflections in the former and triflections in the latter
preserve the lattice. Its covering radius is
$1$; this can be seen as follows. 
According to \ecite{splag}{p. 126}, the real form of the lattice
$\set{(x,y,z)\in\Eisenstein^3}{x\congruent y\congruent z\ (\mod \theta)}$
is the $E_6$ root lattice scaled up by $(3/2)^{1/2}$. This
identification can be used to show that the real form of
$D_3(\theta)$ is the real form of $E_6^*$ scaled up by
$(3/2)^{1/2}$, where $E_6^*$ is the dual (over $\Z$) of $E_6$.
By \ecite{splag}{p. 127}, the
covering radius of $E_6^*$ is $(2/3)^{1/2}$, so the
covering radius of $D_3(\theta)$ is $1$.

The Coxeter-Todd lattice $\cct$ is a selfdual
$\Eisenstein$-lattice that is spanned by its long roots, which
are also its  minimal vectors. It is discussed at length in
\cite{jhc:Coxeter-Todd};
we quote just one of the definitions given there.
$$
\cct={1\over\theta}\set{\vec{x}{6}\in\eisen }{
x_i\equiv x_j
\pmod{\theta},\;\; \sum x_i\congruent 0\pmod{3}}\;.
$$
%Its covering radius is $\sqrt{4/3}$, but the proofs of 
%Thms.~\tag{thm-eI7,1} and \tag{thm-qI5,1} show that the Leech
%holes of $\cct$ are more important for our purposes than the
%deep holes.
It automorphism group is the  finite
complex reflection group $6{\cdot}U_4(3){:}2$, and $\cct$
shares many
interesting properties with $E_8$ and the Leech lattice $\leech$. We refer to 
\cite{jhc:Coxeter-Todd} for details.

The quaternionic Barnes-Wall lattice is
$$
\qbwl={1\over 1+\i }\set{\vector{x}{4}\in\Hurwitz }{
x_i\congruent x_j \pmod{(1+\i )\Hurwitz},\;\;
\sum x_i\in2\Hurwitz}\;.
$$
We may recognize the real form of $2^{1/2}\qbwl$ by identifying
the vector 
$$
(a_1+b_1\i +c_1\j +d_1\k ,\ldots,a_4+b_4\i +c_4\j +d_4\k )
$$
with the vector in $\R^{16}$ whose coordinates we arrange in the
square array
% all this faffing about is to make sure that the array is
% square no matter which fonts are used
%
% user parameters
\newdimen\vpadding
\newdimen\hpadding
\vpadding=0pt
\hpadding=5pt
\newbox\boxa
\newbox\boxb
\newbox\boxc
\newbox\boxd
\setbox\boxa=\hbox{$\displaystyle\matrix{a_1&a_2\cr d_1&c_2\cr}$}
\setbox\boxb=\hbox{$\displaystyle\matrix{a_3&a_4\cr d_3&c_4\cr}$}
\setbox\boxc=\hbox{$\displaystyle\matrix{b_1&d_2\cr c_1&b_2\cr}$}
\setbox\boxd=\hbox{$\displaystyle\matrix{b_3&d_4\cr c_3&b_4\cr}$}
\newdimen\maxx
\newdimen\maxy
\newdimen\maxdim
\maxx=\wd\boxa
\ifnum\maxx<\wd\boxb \maxx=\wd\boxb\fi
\ifnum\maxx<\wd\boxc \maxx=\wd\boxc\fi
\ifnum\maxx<\wd\boxd \maxx=\wd\boxd\fi
\advance\maxx by \hpadding
\maxy=\ht\boxa
\ifnum\maxy<\ht\boxb \maxy=\ht\boxb\fi
\ifnum\maxy<\ht\boxc \maxy=\ht\boxc\fi
\ifnum\maxy<\ht\boxd \maxy=\ht\boxd\fi
\advance\maxy by \vpadding
\ifnum\maxx<\maxy \maxdim=\maxy \else \maxdim=\maxx\fi
\def\pad#1{\hbox to\maxdim{\hfil\vbox to\maxdim{\vfil #1\vfil}\hfil}}
$$
{4\over\sqrt8}\;
\vcenter{\halign{\vrule\pad{\box#}\vrule&\pad{\box#}\vrule\cr
\noalign{\hrule}%
\boxa&\boxb\cr
\noalign{\hrule}%
\boxc&\boxd\cr
\noalign{\hrule}}}
$$
where the inner product is the usual one on $\R^{16}$. This
array may be taken to be (say) the left 4 columns of the
$4\times6$ array in the MOG
description \ecite{ATLAS}{p.~97} of the Leech lattice
$\leech$, and then the real form of $2^{1/2}\qbwl$ is visibly
the real Barnes-Wall lattice $\rbwl$ 
\ecite{splag}{Chap.~4}. 

\beginproclaim Theorem
\Tag{thm-barnes-wall-properties}.
The  lattice $\qbwl$ is selfdual and spanned by its
minimal vectors, which have norm
$2$. Its automorphism group is generated by the biflections in
its minimal vectors. Each class of $\qbwl$ modulo $\qbwl(1+\i )$ is
represented by a vector of norm at most $3$. 
The deep holes of $\qbwl$ coincide with the set 
$\set{\lambda(1+\i
)^{-1}}{\lambda\in\qbwl,\lambda^2\congruent1(\mod~2)}$.\endproclaim  

\beginproof{Proof:}
Proofs of all claims except the  last appear in
\ecite{bachoc:selfdual-Hurwitz-lattices}{Sect.~4.6}.  
Most of the rest of the work has been done for us by Conway and Sloane
\ecite{jhc:laminated-lattices}{Sect.~5}. They
showed that the deep
holes of $\rbwl$ nearest $0$ are the halves of certain
vectors $v\in\rbwl$ of norm 12, and further that 
such $v$ are not congruent modulo 2 to minimal vectors of $\rbwl$.
(They write $\Lambda_{16}$ for  $\rbwl$.)
After rescaling, we find that the deep holes of $\qbwl$ nearest
$0$ are the
halves of certain elements $v$ of norm $6$ in $\qbwl$. Since
each such $v$ has even norm 
and is not congruent modulo $2=-(1+\i)^2$ to any root,
it must map to  $0$ in $\qbwl/\qbwl(1+\i )$. Therefore
$v=\lambda(1+\i )$ for some $\lambda$ of norm $3$ in $\qbwl$ and so the
deep holes nearest $0$ have the form
$v/2=\lambda(1+\i )/2=(\lambda\i)(1+\i)^{-1}$. 

The deep holes of $\qbwl$ are the translates
by lattice vectors of the deep holes nearest zero. That is, the
set of deep holes coincides with the set
$$ 
\set{\lambda(1+\i )^{-1}}{\hbox{$\lambda\in\qbwl$ is congruent modulo
$1+\i $ to a norm $3$ lattice vector}}\;.
$$
The norms of any two lattice vectors that are
congruent modulo $1+\i $ have the same parity. Since each lattice
vector is congruent to some vector of norm $0$, $2$ or $3$, the
set above coincides with the one  in the statement of the theorem.
\endproof

%
% INTERESTING BUT TANGENTIAL FACTS ABOUT BARNES-WALL LATTICE 
%
%Note that $\aut\rbwl$ is much larger
%than $\aut\qbwl$, having  order
%$2^{21}3^55^27=89\,181\,388\,800$ and structure $2^{1+8}_+\cdot
%O_8^+(2)$.  
%Embedding our array of
%coordinates in the MOG array also identifies $2^{1/2}\cdot\qbwl$
%as an $\Hurwitz$-sublattice of
%the quaternionic Leech lattice $\qleech$, which is described in 
%\oldcite{raw:quaternionic-leech}. The automorphism group of
%$\qleech$ is $2\cdot G_2(4)$, and inspection of the list of maximal
%subgroups of $G_2(4)$ in \oldcite{ATLAS} or
%\oldcite{raw:quaternionic-leech} shows that $O_6^-(2)$ is not
%involved in $G_2(4)$, and hence fairly few automorphisms of
%$\qbwl$ extend to $\qleech$. This contrasts with the real case;
%every automorphism of $\rbwl$ extends to $\leech$ (p. 131, \oldcite{splag}). 
%
%Finally, it is amusing to note that $\qbwl$ may be defined as
%the sublattice of $\Hurwitz\tensor\Hurwitz$ spanned by the
%vectors of the forms $(\hbox{norm }1)\tensor(\hbox{norm }2)$ and
%$(\hbox{norm }2)\tensor(\hbox{norm }1)$, where
%$\Hurwitz\tensor\Hurwitz$ is equipped with the inner product
%$$
%\ip{a\tensor b}{c\tensor d}=\re(\bar ac)\bar bd\;.
%$$

Now we describe some indefinite selfdual lattices. 
The lattice  $\rI n,m$
is the $\ring$-module $\ring^{n+m}$ equipped with the
inner product given by the diagonal matrix 
$$ 
\diag[+1,\ldots,+1,-1,\ldots,-1] 
$$
with $n$ (resp. $m$) $+1$'s (resp. $-1$'s). 
The lattice $\rII1,1$ is the module $\ring^2$ with 
inner product matrix $\smallmatrix0110$. If $\ring=\Eisenstein$
or $\Hurwitz$ then $\rII1,1\isomorphism \rI1,1$ because one can
find a norm $1$ vector in the former lattice. If
$\ring=\Gaussian$ then $\II1,1$ is even, whereas $\I1,1$ is
odd. We define the Gaussian lattices $\gII 4m+n,n$ to be the
lattices 
$$ 
\gII 4m+n,n=\geeight\oplus\cdots\oplus\geeight\oplus
	\gII1,1\oplus\cdots\oplus\gII1,1\;, 
$$
where there are $m$ summands $\geeight$ and $n$ summands
$\gII1,1$. These lattices are even and selfdual.
By
Thm.~\tag{thm-indefinite-selfdual-classification},
every indefinite selfdual lattice over $\ring$ appears among the
examples just given. In particular,
$\cct\oplus\eII1,1\isomorphism\eI7,1$ and
$\qbwl\oplus\qII1,1\isomorphism\qI5,1$.

\section{\Tag{sec-hyperbolic-space}. Hyperbolic space}

The
hyperbolic space
$\khnp$ ($n\geq0$) is defined as the image in projective
space $\kpnp$ of the
set of vectors of negative norm in $\K^{n+1,1}$; its boundary
$\dkhnp$ is the image of the (nonzero) null vectors. We
write elements of $\K^{n+1,1}$ in the
form $(\lambda;\mu,\nu)$ with $\lambda\in\K^{n,0}$ and
$\mu,\nu\in\K$, with inner product  
$$
\ip{(\lambda_1;\mu_1,\nu_1)}{(\lambda_2;\mu_2,\nu_2)}=
\ip{\lambda_1}{\lambda_2} + \bar{\mu}_1\nu_2 + \bar{\nu}_1\mu_2\;.
$$
This corresponds to a decomposition
$\K^{n+1,1}\isomorphism\K^{n,0}\oplus\smallmatrix{0}{1}{1}{0}$. 
We will often refer to points in projective space by naming
vectors in the underlying vector space.

It is convenient to distinguish  the
isotropic vector $(0;0,1)$ and give it the name $\rho$.  
Every point of $\khnp\cup\dkhnp$ except
$\rho$ has a unique preimage in $\K^{n+1,1}$ with inner product
$1$ with $\rho$, and so we may make the identifications
$$
\eqalign{
\khnp=\{(\lambda;1,z):\lambda\in\K^n,\;
\lambda^2+2\re(z)>0\}\;.\cr
\dkhnp\setminus\{\rho\}=\{(\lambda;1,z):\lambda\in\K^n,\;
\lambda^2+2\re(z)=0\}\;.\cr
}\eqno\eqTag{eq-half-space-coordinates}
$$ 
We define the height of a vector $v\in\K^{n+1,1}$ to be
$\height v=\ip{\rho}{v}$. For $v=(\lambda;\mu,\nu)$ this is simply
$\mu$. For vectors of any fixed norm, the height function
measures how far away from $\rho$ the corresponding points in
projective space are; the smaller the height, the
closer to $\rho$. We will sometimes say that a vector $v'$  has
height  less than that of another vector $v$. By this we
will mean $|\height v'|<|\height v|$.

We say that the vector  $(\lambda;\mu,\nu)$ of height $\mu\neq0$ lies
over $\lambda\mu^{-1}\in\K^{n}$. It is obvious that all the
scalar multiples of any given vector of nonzero height lie over
the same point of $\K^n$, so we may think of points in
projective space (except for those in $\rho^\perp$) as lying
over elements of $\K^n$. The geometric content of this
definition is that the lines in $\kpnp$ passing through $\rho$
and meeting $\khnp$ are in one-to-one correspondence with 
$\K^n$. The points in the line associated to
$\lambda\in\K^n$ are  the scalar multiples of those of the
form $(\lambda;1,z)$ with $z\in\K$, which are precisely the points of
$\kp{n+1}$ lying over $\lambda$. We gave two special cases in
Eq.~\eqtag{eq-half-space-coordinates}. In particular, the family of
height one isotropic vectors  lying over $\lambda$ is
parameterized by the elements of $\im\K$. This description of
$\dkhnp\setminus\{\rho\}$ as a bundle over $\K^n$ with fiber
$\im\K$ will help us relate the properties of lattices in $\K^n$
to properties of groups acting on $\khnp$.

The subgroup  of $U(n+1,1;\K)$ fixing $\rho$ contains
transformations $T_{x,z}$  (with $x\in\K^n$,
$z\in\im\K$) defined by
$$\eqalignno{
\rho & \mapsto  \rho  \cr 
\llap{$T_{x,z}$:\qquad}
(0;1,0) & \mapsto  (x;1,z-x^2/2)&\eqTag{eq-def-of-translation}  \cr
(\lambda;0,0) & \mapsto  (\lambda;0,-\ip{x}{\lambda}) 
\rlap{\qquad for each  $\lambda\in\K^n$.} 
\cr}$$
(The map is defined in terms of some unspecified but
fixed inner product on $\K^n$.)
We call these maps  translations. If we regard
elements 
of $\K^{n+1,1}$
as column vectors then  $T_{x,z}$ acts by multiplication on the
left by the matrix
$$
\pmatrix{
I_n 	& x 		& 0 	\cr 
0 	& 1		& 0	\cr
-x^*	& z-x^2/2	& 1	\cr}\;.
$$
We have  written  $x^*$
for  the linear function $y\mapsto\ip{x}{y}$ 
on $\K^{n,0}$ defined by $x$. We have the relations
$$\displaylines{
T_{x,z}\circ T_{x',z'} = 
T_{x+x',z+z'+\im\ip{x'}{x}} 
\hfil\llap{\eqTag{eq-product-of-translations}}\hfilneg \cr
T_{x,z}^{-1} = T_{-x,-z}
\hfil\llap{\eqTag{eq-inverse-of-translation}}\hfilneg \cr 
T_{x,z}^{-1}\circ T_{x',z'}^{-1}\circ T_{x,z}\circ T_{x',z'} =
T_{0,2\im\ip{x'}{x}}\;,
\hfil\llap{\eqTag{eq-commutator-of-translations}}\hfilneg
\cr}$$
which are most easily verified in the order listed.
These relations make it
clear that  the
translations form a group and that its
center and commutator subgroup coincide and consist of the
$T_{0,z}$. We call elements of this subgroup
central translations. The translations form 
a (complex or quaternionic) Heisenberg group which  acts
freely and  transitively
on $\dkhnp\setminus\{\rho\}$. 
%In Section~\tag{sec-26dim} our arrangement
%of mirrors was invariant under a group $\leech\isomorphism\Z^{24}$;
%these symmetries are translations in the current sense. 
If $v\in\K^{n+1,1}$ lies over $\lambda\in\K^n$ then $T_{x,z}(v)$
lies over $\lambda+x$. That is, the translations act in the
natural way (by translations) on the points of $\K^n$
over which vectors in 
$\K^{n+1,1}$ lie. 

We note that these constructions all make sense when $\K=\R$,
and even simplify. Since $\im\R=0$, the translations form
an abelian group, which is just the obvious set of translations
in the usual upper half-space model for $\rhnp$. The 
obvious projection map from the upper half-space to $\R^n$ 
carries points of $\rhnp$ to the points of $\R^n$ over which
they lie, in the sense defined above. This is the source of the
terminology. 

The simultaneous stabilizer of $(0;1,0)$ and $(0;0,1)$
is the  unitary group $U(n,0;\K)$, which fixes pointwise the
second summand of the decomposition
$\K^{n+1,1}=\K^{n,0}\oplus\K^{1,1}$. If $S$ is an
element of this unitary group then matrix computations reveal 
$$
S\circ T_{x,z}\circ S^{-1}=T_{Sx,z}\;.
\eqno\eqTag{eq-conjugate-of-translation}
$$
%This is useful in its own right and also
%shows that the group of translations is  normal in the full
%stabilizer of   
%$\rho$.

\section{\Tag{sec-reflections}. Reflections in Lorentzian
lattices}

The Lorentzian lattices we will consider all have the form
$\Lambda\oplus\II1,1$, where $\Lambda$ is a
positive-definite $\ring$-lattice and $\II1,1$ is 
the 2-dimensional selfdual lattice defined
by the matrix
$\II1,1=\smallmatrix0110$.
In general we will write $L$ for a Lorentzian lattice
$\Lambda\oplus\II1,1$, where $\Lambda$ and even $\ring$ may be
left unspecified, except that $\Lambda$ will always be
positive-definite.  
We write elements of $L=\Lambda\oplus\II1,1$ in the form
$(\lambda;\mu,\nu)$ with $\lambda\in \Lambda$ and
$\mu,\nu\in\ring$. This embeds $L$ in the description of
$\K^{n+1,1}$ given in Section~\tag{sec-hyperbolic-space} and
allows us to transfer 
to $L$ several important concepts defined there. In particular,
$\rho=(0;0,1)$ is an element of $L$ and we define the height of
elements of $L$ as before. For $v\in L$ of nonzero height, 
the point of $\Lambda\tensor\R$ over which $v$ lies is in
$\Lambda\tensor\Q$ but not necessary in $\Lambda$.

There are two basic ideas in this section. First,  this
description of $L$ provides a  way to write down a large
collection of reflections of $L$, essentially parameterized by
the elements of a discrete Heisenberg group of translations. The
second idea is that if $r$ is a root of $L$ and $v$ is a null
vector in $\K^{n+1,1}$, and if $r$ and $v$ lie over points of
$\K^n$ that are sufficiently close, then by applying a
reflection of $L$ one can reduce the height of $v$. (This
reflection might be in some root other than $r$.)

Both of these ideas can be found in the simpler setting of real
hyperbolic space,  in Conway's study \cite{jhc:26dim} of the
automorphism group of the Lorentzian $\Z$-lattice
$\II25,1=\leech\oplus\II1,1$. Here $\leech$ is the Leech
lattice, and Conway found a set of reflections permuted freely by a
group of translations naturally isomorphic to the additive group
of $\leech$. By using facts about the covering radius of
$\leech$ together with the second idea described above, he was
able to prove that these reflections generate the entire
reflection group of $\II25,1$.

The major complication in transferring this approach to
our setting is that the discrete group of translations is no
longer a copy of $\Lambda$ but a central extension of $\Lambda$
by $\im\ring$. This  issue  dramatically complicates the
precise formulation (Table~\tag{tab-table-of-heights}) of the second main idea.
For example, it is complicated to state exactly what happens
when one can't {\it quite\/} reduce the height of $v\in\K^{n+1,1}$ by
using a reflection.

We begin by finding the translations in $\aut L$ and showing that
under simple conditions, $\reflec L$ contains a large number of
them. The translation $T_{x,z}$
preserves $L$ just if $x\in \Lambda$ and $z-x^2/2\in\ring$. If
$\ring=\Eisenstein$ or $\Hurwitz$ then for any given $x\in \Lambda$
we may  choose $z\in\im\K$ such that $T_{x,z}\in\aut L$, by
taking $z=0$ or $\theta/2$ according as $x^2$ is even or
odd. If $\ring=\Gaussian$ then such a $z$ exists if and only if
$x^2$ is even;  $z$ may then be taken to be zero.  
The different rings behave differently because
$\Eisenstein$ and $\Hurwitz$ contain elements with half-integral
real parts, while $\Gaussian$ does not.  
All the central translations $T_{0,z}$ with
$z\in\ring$  lie in $\aut L$\emdash they fix $\Lambda$ pointwise and act by
isometries of  $\II1,1$.
The assertions of the next lemma are precise
formulations of the idea that if $\aut\Lambda$ contains many 
reflections then $\reflec L$ contains many translations.

\beginproclaim Lemma
{\Tag{thm-stabilizer-of-rho}}. 
Let $L=\Lambda\oplus\II1,1$ for some  positive-definite
$\ring$-lattice $\Lambda$. Define
$$\displaylines{
\Lambda_0=\set{x\in \Lambda}{T_{x,z}\in\reflec L
\hbox{\rm\  for some }z\in\im\K}
\hbox{ and}\cr
\cals=\{z\in\im\ring \mid T_{0,z}\in\reflec L\}\;. 
\cr}$$
\item{\rom1} If $\ring=\Eisenstein$ or $\Hurwitz$ then
$\Lambda_0$ contains the short roots of $\Lambda$.
\item{\rom2} If $r$ is a long root of $\Lambda$ then $2r\in
\Lambda_0$. Furthermore, 
if $r$ has inner product $1$ with some element of $\Lambda$ then
$r$ itself lies in $\Lambda_0$.
\item{\rom3} $\cals$ contains the integral span of the elements
$2\im\ip{x}{y}$ with $x,y\in\Lambda_0$.
\item{\rom4} If the roots of $\Lambda\neq\{0\}$ span $\Lambda$
up to finite index then the stabilizer of $\rho$ in $\reflec L$
has finite index in the stabilizer in $\aut L$.
\par
\endproclaim 

%\remark{Remarks:}
%By Eq.~\eqtag{eq-product-of-translations}
%and \eqtag{eq-inverse-of-translation}, $\Lambda_0$ is closed
%under addition and negation, 
%so  assertion
%\rom1 makes sense.
%We will not need  an analogue of \rom2 for short roots in
%Gaussian lattices.
%In Thm.~\tag{thm-rho-stabilizer-finite-index-2}, similar but
%stronger hypotheses are used to 
%obtain similar but stronger conclusions.

\beginproof{Proof:}
Let $R$ be a $\xi$-reflection of $\Lambda$ with mirror $M$. We
regard $R$ as acting on $L$, fixing the summand $\II1,1$ 
pointwise. If $T_{x,z}\in\aut L$ then  
$T_{x,z}^{-1}\circ R\circ T_{x,z}\in\reflec L$. By
Eqs.~\eqtag{eq-inverse-of-translation}, 
\eqtag{eq-conjugate-of-translation} and
\eqtag{eq-product-of-translations}, 
$$ 
T_{x,z}^{-1}\circ R\circ T_{x,z}\circ R^{-1}=
T_{-x,-z}\circ T_{Rx,z} = T_{Rx-x,-\im\ip{Rx}{x}}\;,
$$
proving that
$$
Rx-x\in \Lambda_0 
\eqno\eqTag{eq-got-translation}
$$
for all $x\in \Lambda$ and all reflections $R$ of $\Lambda$.
The geometric picture behind this computation is
that both $M$ and its translate by $T_{x,z}^{-1}$ pass 
through $\rho$ and are parallel there;
we have  constructed a translation out of
reflections in two parallel mirrors.

\rom1
If $r$ is a short root of $\Lambda$ then we let $x=r\w$ and $R$
be the $(-\w)$-reflection in $r$. Then $\Lambda_0$ contains
$Rx-x=r(-\w)\w-r\w=r(-\wbar-\w)=r$. 

\rom2
If $r$ is a long root of $\Lambda$ then we let $x=-r$, let $R$ be
the biflection in $r$, and observe $Rx-x=2r$. To
prove the second claim, suppose $x\in \Lambda$ has inner product
$-1$ with $r$ and take $R$ to be the biflection in $r$. Then
$Rx-x\in \Lambda_0$ is proportional to $r$ and has inner product
$2$ with $r$, so it coincides with $r$.

\rom3 Follows immediately from
Eq.~\eqtag{eq-commutator-of-translations} by taking commutators of 
translations of $\reflec L$.

\rom4
The null vectors  of height $1$ in $L$ are exactly those vectors
$(\lambda;1,z)$ with $\lambda\in \Lambda$, $z\in\ring$ and $\re
z= -\lambda^2/2$, and the translations in $\aut L$ permute them
transitively. Since the simultaneous stabilizer of $\rho$ and
one of these, say $(0;1,0)$, is the finite group $\aut
\Lambda$,  it 
suffices to prove that the group of translations in $\reflec L$
has finite index in the group of those in $\aut L$. This follows
from \rom1--\rom3:  $\Lambda_0$ has finite index in $\Lambda$
and $\cals$ has finite index in $\im\ring$. 
\endproof

Now we will exhibit a large number of reflections of $L$.
It is straightforward to
enumerate the roots of $L$ of any given height $h$; 
For $h=1$ one finds that these are the vectors 
$$\eqalign{
\hbox{Norm 2:}\qquad& (\lambda;1,z),\qquad  \re z=(2-\lambda^2)/2\cr
\hbox{Norm 1:}\qquad& (\lambda;1,z),\qquad  \re z=(1-\lambda^2)/2\;,
\cr}$$
with $\lambda\in \Lambda$ and $z\in\ring$.
If $\ring=\eisen\hbox{ or }\Hurwitz$, then 
height $1$ 
roots of both norms lie over each $\lambda\in \Lambda$, and  the 
translations of $L$ act simply  transitively on each set. If
$\ring=\gauss$ then height one  roots lie over each
$\lambda\in \Lambda$: long roots over $\lambda$ of even norm and
short roots over $\lambda$ of odd norm.
Again, the translations  act simply transitively on each
set. The differing behavior of the different rings
is another manifestation of the fact that
$\Eisenstein$ and $\Hurwitz$ have elements
with half-integer real part, while $\Gaussian$ does not.
One may also enumerate roots of larger heights\emdash for example, if $\Lambda$
is an $\eisen$-lattice, then there are short roots of $L$
of height $\theta$ 
over each $\lambda\theta^{-1}\in \Lambda\theta^{-1}$ with 
$\lambda^2\congruent 1$ modulo $3$.

Now we will discuss the second idea of this section:
the effects of reflections in roots of small
height $h$. This will occupy the rest of the section.

\beginproclaim Lemma
{\Tag{thm-gaussian-long-roots-height-reduction}}.  Suppose
$\Lambda$ is a positive-definite $\Gaussian$-lattice and
$L=\Lambda\oplus\II1,1$. Suppose $r$ is a long root of $L$ of
height $1$ lying over $\lambda\in\Lambda$, and that $v$ is an
isotropic vector of $L\tensor\R$ of height $1$ that lies over
$\ell\in \Lambda\tensor\R$. Suppose that
$(\ell-\lambda)^2<\sqrt3$. Then there is another long root
$r'\in L$ of height $1$, also lying over $\lambda$, such that
the biflection in $r'$ reduces the height of $v$.
\endproclaim 

\beginproof{Proof:}
Since $v$ has height $1$ and norm $0$ and lies over $\ell$, we
know that for some $w\in\im\K$ we have
$v=(\ell;1,w-\ell^2/2)$. Similarly, we deduce that
$$
r=\(\lambda;1,z_0+{2-\lambda^2\over2}\)
\eqno\eqTag{eq-coordinates-of-r}
$$ 
for some $z_0\in\im\K$. Every
other height $1$ long root of $L$ lying over $\lambda$ has the
form $r'=r+(0;0,z)$ for some $z\in\im\Gaussian$. We will obtain
the theorem by choosing $z$ appropriately.

We have
$$\eqalignno{
\ip{r'}{v}
&=\ip{\lambda}{\ell}+\(w-\ell^2/2\)+\(\bar{z}_0+\bar{z}+1-\lambda^2/2\)
	&\eqTag{eq-bar-1}\cr
&=1-{1\over2}\(\ell^2-2\ip{\lambda}{\ell}+\lambda^2\)+w+\bar{z}_0+\bar{z}
	&\eqTag{eq-bar-2}\cr
&=1-{1\over2}\(\ell^2-\ip{\lambda}{\ell}-\ip{\ell}{\lambda}+\lambda^2\)
	+{1\over2}\(\ip{\lambda}{\ell}-\ip{\ell}{\lambda}\)
	+w+\bar{z}_0+\bar{z}
	&\eqTag{eq-bar-3}\cr
&=\[1+{1\over2}(\ell-\lambda)^2\]
	+\[\im\ip{\lambda}{\ell}+w+\bar{z}_0+\bar{z}\]
	&\eqTag{eq-bar-4}\cr
&=a+B\;,&\eqTag{eq-bar-5}\cr
}
$$
where $a$ is the first bracketed expression and $B$ is the
second. The important thing to observe here is that $a$ depends
on $(\ell-\lambda)^2$, which is bounded by hypothesis, and $B$
depends on $z$, over which we have some control. Let $v'$ be the
image of $v$ under biflection in $r'$. Since
$v'=v-r'\ip{r'}{v}$, we have
$$\eqalign{
\ip{\rho}{v'}
&=\ip{\rho}{v}-\ip{\rho}{r'}\ip{r'}{v}\cr
&=1-(1)(a+B)\cr
&=(\ell-\lambda)^2/2-B.
\cr}$$
Since the covering radius of $\im\Gaussian$ is $1/2$, we may choose $z$
so that $|B|\leq1/2$. Then
$$ 
\left|\ip{\rho}{v'}\right|^2
=\left|(\ell-\lambda)^2\over2\right|+|B|^2
<\left|\sqrt3\over2\right|^2+\left|1\over2\right|^2
=1\;,
$$
so that $\height v'<\height v$.
\endproof

\beginproclaim Lemma
{\Tag{thm-short-root-height-reduction}}. 
Suppose $\Lambda$ is a positive-definite $\ring$-lattice and
$L=\Lambda\oplus\II1,1$. Let $h=1$ if $\ring=\Gaussian$, $h=1$
or $\theta$ if $\ring=\Eisenstein$, and $h=1$ or $1+\i$ if
$\ring=\Hurwitz$. Suppose $r$ is a short root of $L$ of height
$h$ lying over $\lambda h^{-1}$, with $\lambda\in \Lambda$. Let
$v\in L\tensor\R$ be isotropic, have height $1$, and lie over
$\ell\in \Lambda\tensor\R$. Set $D^2=(\ell-\lambda h^{-1})^2$
and suppose $D^2\leq 1/|h|^2$. Then there exists a short root
$r'$ of $L$, also of height $h$ and lying over $\lambda h^{-1}$,
such that one of the following holds:
\item{\rom1}
some reflection in $r'$ carries $v$ to a vector of smaller
height than $v$;
\item{\rom2}
$D^2=1/|h|^2$ and $\ip{r'}{v}=0$; or
\item{\rom3}
$\ring=\Hurwitz$, $h=1+\i$, $D^2=1/|h|^2=1/2$ and
$\ip{r'}{v}=(1+\i)/2$. 
\par
\endproclaim 

\beginproof{Proof:}
From the given norms and heights of $v$ and $r$, together with
the fact that they lie over $\ell$ and $\lambda h^{-1}$, we
deduce
$$ 
v=\(\ell;1,w-\ell^2/2\)
\qquad\hbox{and}\qquad
r=\(\lambda;h,z_0+{1-\lambda^2\over2|h|^2}h\)
\eqno\eqTag{eq-coordinates-of-r-2}
$$
for some $w\in\im\K$ and $z_0\in\K$ such that
$\re(\bar{h}z_0)=0$. The other height $h$ short roots of $L$ lying over
$\lambda h^{-1}$ have the form $r'=r+(0;0,z)$ for $z\in\ring$
such that $\re(\bar{h}z)=0$. The basic idea is similar to that
of Lemma~\tag{thm-gaussian-long-roots-height-reduction}: we will try to choose $z$, together with a
unit $\xi$ of $\ring$, such that the $\xi$-reflection in $r'$
carries $v$ to a vector of smaller height. It may happen that no
such choice is possible, which leads to the cases \rom2 and
\rom3 of the theorem.

A calculation similar to Eqs.~\eqtag{eq-bar-1}--\eqtag{eq-bar-5} reveals that
$$\eqalignno{
\ip{r'}{v}
=&\ip{\lambda}{\ell}+\bar h\(w-{\ell^2\over2}\)
	+\overline{\(z_0+z+{(1-\lambda^2)h\over2|h|^2}\)}\cr
=&h^{-1}\[\({1\over2}-{|h|^2\over2}D^2\) 
	+\(|h|^2\im\ip{\lambda h^{-1}}{\ell}+|h|^2w
	+h{\bar z}_0+h\bar z\)\]\cr
=&h^{-1}[a+B]&\eqTag{eq-ip-with-r'-2}
\cr}$$
where $a=(1-|h|^2D^2)/2$ is the real part of the term in
brackets and $B$ is the imaginary part. The slight difference
between the terms $a$ in Eqs.~\eqtag{eq-bar-5} and
\eqtag{eq-ip-with-r'-2} is due to the 
replacement of  $2-\lambda^2$ in Eq.~\eqtag{eq-coordinates-of-r} by
$1-\lambda^2$ in Eq.~\eqtag{eq-coordinates-of-r-2}, which 
is due to the fact that $r$ is now a short root.

We take $v'$ to be the image of $v$ under $\xi$-reflection in
$r'$ (we will choose $\xi$ later). Since
$v'=v-r'(1-\xi)\ip{r'}{v}$, we have
$$\eqalignno{
%\height v'=
\ip{\rho}{v'}
&=\ip{\rho}{v}-\ip{\rho}{r'}(1-\xi)\ip{r'}{v}\cr
&=1+{h(\xi-1)\bar h\over|h|^2}[a+B]\;.&\eqTag{eq-height-of-v'-2}
\cr}$$
By hypothesis, $D^2\leq1/|h|^2$, so $a\in[0,1/2]$.
We may change the value of $B$ by $h\bar{z}$, where $z$ may be
any element of $\ring\cap\im(\bar{h}\K)$. That is, we may change
$B$ by any element of
$$\eqalign{
h\cdot\overline{\ring\cap(h\cdot\im\K)}
&=h\cdot(\ring\cap(\im\K)\cdot\bar{h})\cr
&=(h\ring)\cap h\cdot(\im\K)\cdot\bar{h}\cr
&=(h\ring)\cap \im\K\cr
&=\im(h\ring)\;.
\cr}$$
We will try to choose $\xi$ and $z$ so that
\eqtag{eq-height-of-v'-2} has absolute value less than
$1=\ip{\rho}{v}$. This requires treating the different
possibilities for $\ring$ and $h$ separately. We will treat only
the case $\ring=\Hurwitz$, $h=1+\i$, which is more
involved than the other four cases.

We  write $B$ as $b\i +c\j +d\k $
with $b,c,d\in\R$. 
We first carry out a computation that will allow us
to use the $24$ units of $\Hurwitz$ effectively: we claim that
there is a unit $\xi'$ of $\Hurwitz$ with $\re\xi'=-1/2$ such
that 
$$
|1+\xi'(a+B)|^2=(a-1/2)^2+(|b|-1/2)^2+(|c|-1/2)^2+(|d|-1/2)^2\;.
\eqno\eqTag{eq-foo3}
$$
For any unit $\xi'$, the left side is just
the square of the distance between $a+B$ and $-\bar{\xi}'$
(proof: left-multiply by $1=|-\bar{\xi}'|^2$). Setting
$-\bar{\xi}'=(1\pm \i \pm \j \pm \k )/2$, with each of its $\i $, 
$\j $ and $\k $ components having the same sign as the corresponding
component of $B$ (or a random sign if that component 
vanishes), the right hand side becomes another expression for
this squared distance, proving the claim.

Next, we investigate our freedom to choose $B$. By choice of $z$
we may vary $B$ by any element of $\im(h\Hurwitz)$. It is easy
to check that
$$ 
\im((1+\i )\Hurwitz)=\set{b\i +c\j +d\k }{b,c,d\in\Z,\ 
b+c+d\congruent0\pmod2}\;.  
$$
That is, $\im(h\Hurwitz)$ is spanned by $\j +\k $, $\k +\i $ and
$\i -\k $, so by choice of $z$ we may take $b\in(-1,1]$ and
$c,d\in[0,1)$.

Suppose for a moment that there is a unit
$\xi$ of $\Hurwitz$ such that 
$$
\xi'={h(\xi-1)\bar h\over|h|^2}\;,
\eqno\eqTag{eq-relation-between-xi's}
$$ 
where $\xi'$ is as in Eq.~\tag{eq-foo3}.
Then by
Eqs.~\eqtag{eq-height-of-v'-2}, \eqtag{eq-relation-between-xi's} and
\eqtag{eq-foo3}, 
$$\eqalignno{
|\height v'|^2
%&=|\ip{\rho}{v'}|^2\cr
&=|1+\xi'(a+B)|^2\cr
&=(a-1/2)^2+(|b|-1/2)^2+(|c|-1/2)^2+(|d|-1/2)^2\;.
&\eqTag{eq-foo4}
\cr}$$
We have already shown that $a\in[0,1/2]$. By this and
the constraints on $b$, $c$ and $d$ obtained above, we see that
the right hand side of Eq.~\eqtag{eq-foo4} is less than
$1=|\height v|^2$ (so that conclusion \rom1 applies) unless
$a=0$, $b\in\{0,1\}$ and $c=d=0$. In this exceptional case,
$a=0$ implies that $D^2=1/|h|^2$, and $\ip{r'}{v}$ can be read
from Eq.~\eqtag{eq-ip-with-r'-2}. If $b=0$ we have $r'\bot v$
and conclusion \rom2 applies, and if $b=1$ then conclusion \rom3
applies.

It remains only to show that given a unit $\xi'$ of $\Hurwitz$
with $\re\xi'=-1/2$, there is another unit $\xi$ of $\Hurwitz$
satisfying Eq.~\eqtag{eq-relation-between-xi's}.  We simply
solve for $\xi$:
$$
\xi=|h|^2\cdot h^{-1}\xi'\bar{h}^{-1}+1 = 
{(1-\i )\xi'(1+\i )\over\sqrt2\phantom{\i (\xi')q}\sqrt2}+1\;.
\eqno\eqTag{eq-relation-between-xi's-2}
$$
The most straightforward way to show that $\xi$ is a unit of
$\Hurwitz$ 
is to simply evaluate the right hand side of
Eq.~\eqtag{eq-relation-between-xi's-2} for each of
the eight possibilities $\xi'=(-1\pm \i \pm \j \pm \k )/2$.
(What is really going on here is that the units of $\Hurwitz$
together with $(1+\i)/\sqrt2$ generate the binary octahedral
group, which normalizes the binary tetrahedral group consisting
of the units of $\Hurwitz$.)
\endproof

The idea used in the proofs of the last two lemmas can also be
applied for some other values of $h$. The cases stated above are
the ones that will be used later, but for completeness we
summarize in table~\tag{tab-table-of-heights} all the cases we
have been able to treat with this method. The table should be
read as follows. Suppose $\Lambda$ is an $\ring$-lattice,
$L=\Lambda\oplus\II1,1$, and $r$ is a short root of $L$ whose
height $h$ appears in the table, lying over $\lambda h^{-1}$,
with $\lambda\in \Lambda$. Suppose $v$ is a primitive null
vector of $L\tensor\R$ of height $1$, lying over $\ell\in
\Lambda\tensor\R$, and that $D^2=(\ell-\lambda h^{-1})^2$
satisfies $D^2\leq R^2$, where $R^2$ is given by the table. Then
there is another root $r'$ of $L$, of the same height and length
as $r$ and also lying over $\lambda h^{-1}$, such that either
some reflection in $r'$ preserves $L$ and reduces the height of
$v$, or else $D^2=R^2$ and $\ip{r'}{v}$ takes one of the values
given in the table. Note that $R^2=1/|h|^2$ in all cases except
that of long roots of height 1 in Gaussian lattices.

\midinsert
\leavevmode\hfil\vbox{%
\halign{
\hfil#\hfil&\qquad\hfil#\hfil&\qquad\hfil$#$\hfil&\quad\hfil$#$\hfil&\qquad\hfil#\hfil\cr
\bf
The~ring~$\ring$&\bf root~length&\hbox{\bf height~$h$}&R^2&$\ip{r'}{v}$\cr  
\noalign{\medskip}
$\Gaussian$&long&1&\sqrt3&$1-{\sqrt3\over2}+{\i\over2}$\cr
&short&1&1&$0$\cr
&&1+\i&1/2&$0$ or $h^{-1}\i$\cr
&&2&1/4&$0$ or $h^{-1}\i$\cr
\noalign{\medskip}
$\Eisenstein$&long&1&1&$-\wbar$\cr
&&\theta&1/3&$-h^{-1}\wbar$\cr
&short&1&1&$0$\cr
&&\theta&1/3&$0$\cr
&&2&1/4&$0$ or $h^{-1}\theta$\cr
&&2\theta&1/12&$0$ or $h^{-1}\theta$\cr
\noalign{\medskip}
$\Hurwitz$&long&1&1&$-\wbar$\cr
&short&1&1&$0$\cr
&&1+\i&1/2&$0$ or ${1\over2}(1+\i)$\cr
&&2&1/4&${1\over2}(1+a\i+b\j+c\k)$ for $a,b,c\in\{0,1\}$\cr
}}\hfil
\medskip
\centerline{Table~\Tag{tab-table-of-heights}. Summary of
Lemmas~\tag{thm-gaussian-long-roots-height-reduction}
and~\tag{thm-short-root-height-reduction}, and generalizations
thereof.}
\endinsert

\section{\Tag{sec-reflection-groups}. The reflection groups}

This section is the heart of the paper: we will apply the
results of Section~\tag{sec-reflections} to find Lorentzian
lattices that are reflective.  We begin by providing a general
criterion for a lattice to be reflective, and give a number of
examples
(Thms.~\tag{thm-cov-radius-implies-reflective}--\tag{thm-examples-of-reflective-lattices}). Then
we will study in much greater detail the lattices
$\Eisenstein^{n,1}$ and $\Hurwitz^{n,1}$ for small $n$
(Lemma~\tag{thm-basics-for-eisenstein-lattices}--Thm.~\tag{thm-qI2,1-qI3,1}),
and also two high-dimensional examples, acting on $\ch7$
and $\qh5$
(Lemma~\tag{thm-linear-independence-descends}--Thm.~\tag{thm-qI5,1}). At
the end of the section we return to low dimensions, discussing
the lattices $\rI1,1$ and $\gII1,1$.  We begin with the most basic of our
results:

\beginproclaim Theorem
\Tag{thm-cov-radius-implies-reflective}.
Suppose $\Lambda$ is a positive-definite $\ring$-lattice which is
spanned up to finite index by its roots and has covering radius $\leq1$. Then 
$L=\Lambda\oplus\II1,1$ is reflective.
Furthermore, if the covering radius is~$<1$ then any two primitive
isotropic vectors of $L$ are equivalent (up to a scalar) under
$\reflec L$.
\endproclaim 

\beginproof{Proof:}
According to Lemma~\tag{thm-stabilizer-of-rho}, the stabilizer
of $\rho$ in $\reflec L$ has finite index in the stabilizer in
$\aut L$. Now we study the $\reflec L$ orbits of primitive null
vectors in $L$. Suppose $v$ is such a vector, that it is not a
multiple of $\rho$, and that it has minimal height in its
$\reflec L$ orbit. Let $\ell$ be the element of
$\Lambda\tensor\Q$ over which it lies, let $\lambda$ be an
element of $\Lambda$ nearest $\ell$, and let $r$ be a short root
of $L$ of height 1 lying over $\lambda$ (or a long root if
$\ring=\Gaussian$ and $\lambda^2$ is even). We must have
$(\ell-\lambda)^2\geq1$, for else
Lemma~\tag{thm-short-root-height-reduction} (or
Lemma~\tag{thm-gaussian-long-roots-height-reduction} if
$\ring=\Gaussian$ and $r$ is long) assures us that $v$ is not of
minimal height in its $\reflec L$ orbit. In particular, if
$\Lambda$ has covering radius~$<1$ then $v$ cannot exist and we
have proven that every null vector of $L$ is equivalent under
$\reflec L$ to a multiple of $\rho$. This is the second
part of the theorem.

In  case the covering radius of $\Lambda$ is exactly
$1$, we can still deduce that there are only finitely many
$\reflec L$ orbits of primitive null vectors in $L$. For if one
cannot reduce the height of $v$ by a reflection, then $\ell$ is
a deep hole of $\Lambda$, and if $\lambda$ is any vertex of the
hole then there is a short root $r$ of $L$ of height 1 that lies
over $\lambda$ and is orthogonal to $v$. Now, $v$ is determined
up to a unit scalar by the point $\ell$ of $\Lambda\tensor\Q$
and the root $r$ (lying over $\lambda$) to which it is
orthogonal. Since the stabilizer of $\rho$ in $\reflec L$
contains a finite-index subgroup of the translations of $L$, we
may take $r$ to lie in some fixed finite set of roots. Then $\ell$ is a
deep hole nearest $\lambda$, for which there are only finitely
many possibilities. That is, there are only finitely many
$\reflec L$ orbits of primitive null vectors in $L$. The fact
that $\reflec L$ has finite index in $\aut L$ follows from this,
together with the fact that for one particular primitive null
vector, namely $\rho$, its stabilizer in $\reflec L$ has finite
index in its stabilizer in $\aut L$.
\endproof

Essentially the same argument, using
Lemma~\tag{thm-gaussian-long-roots-height-reduction} in place of
Lemma~\tag{thm-short-root-height-reduction}, proves the
following theorem.

\beginproclaim Theorem 
{\Tag{thm-gaussian-lattices-reflective}}. 
Suppose $\Lambda$ is an even positive-definite
$\Gaussian$-lattice which is spanned up to finite index by its
roots and has covering radius~$<\root4\of3$. Then
$L=\Lambda\oplus\II1,1$ is reflective and any two primitive null
vectors of $L$ are equivalent (up to a scalar) under $\reflec
L$. 
\QED
\endproclaim 

\beginproclaim Corollary 
{\Tag{thm-examples-of-reflective-lattices}}. 
Let $\Lambda$ be any of the $\ring$-lattices 
$$
\vbox{\halign{
#\hfil&\qquad#\hfil&\qquad#\hfil\cr
$\Gaussian$, $2^{1/2}\Gaussian$, $D_4^\Gaussian$,
$D_6^\Gaussian$ or $\geeight$&&if $\ring=\Gaussian$,\cr
$\Eisenstein$, $\Eisenstein^2$,&$\Eisenstein^3$ or
$D_4(\theta)$&if $\ring=\Eisenstein$, or\cr
$\Hurwitz$,&$2^{1/2}\Hurwitz$, $\Hurwitz^2$ or $\qeeight$&if
$\ring=\Hurwitz$.\cr
}}$$
Then $L$ is reflective. Furthermore, if $\Lambda$ appears in the
first column of the list then any two primitive null vectors of
$L$ are equivalent (up to a scalar) under $\reflec L$.
\endproclaim 

\remark{Remark:}
The lattices appearing
here are all described in Section~\tag{sec-lattice-reference}.
Thms.~\tag{thm-eI2,1-eI3,1-eI4,1} and
\tag{thm-qI2,1-qI3,1} give much more precise information about 
$\reflec L$ for $\Lambda=\Eisenstein$, $\Eisenstein^2$,
$\Eisenstein^3$, $\Hurwitz$ or $\Hurwitz^2$.

\beginproof{Proof:}
All these lattices are spanned by their roots.
The covering radii of the Gaussian lattices are $1/\sqrt2$, $1$,
$1$, $\sqrt{3/2}$ and $1$, and all but the first are even. The
covering radii of the Eisenstein lattices are $1/\sqrt3$,
$\sqrt{2/3}$, $1$ and $1$, and those of the Hurwitz lattices are
$1/\sqrt2$, $1$, $1$ and $1$. The result follows from
Thms.~\tag{thm-cov-radius-implies-reflective} and~\tag{thm-gaussian-lattices-reflective}.
\endproof

We will now study in more detail the reflection groups of some
low-dimensional selfdual Lorentzian
lattices over $\Eisenstein$ and $\Hurwitz$. If $L$ is any
lattice we will write $\shr L$ for the subgroup of $\reflec L$
generated by the reflections in the short roots of $L$.

\beginproclaim Lemma 
{\Tag{thm-basics-for-eisenstein-lattices}}. 
Suppose $\ring=\Eisenstein$, $\Lambda=\Eisenstein^n$ ($n>0$),
and $L=\Lambda\oplus\II1,1$. Then
\item{\rom1}
$\shr L$ contains all the translations of $L$.
\item{\rom2}
$\shr L$ contains a transformation acting trivially on $\Lambda$
and as $\w$ on $\II1,1$.
\item{\rom3}
The stabilizers of $\spanof{\rho}$ in $G$ and $\aut L$ coincide,
where $G$ is the group generated by $\shr L$ and the central
involution $-I$ of $L$. Furthermore, $G\sset\reflec L$.
\par
\endproclaim

\beginproof{Proof:}
\rom1
By Lemma~\tag{thm-stabilizer-of-rho}\rom1, $\reflec L$ contains a translation
$T_{x,z}$ for each $x\in\Eisenstein^n$. The proof shows that
these translations actually lie in $\shr L$. Taking commutators
as in Lemma~\tag{thm-stabilizer-of-rho}\rom3 shows that $\shr L$ contains all the
reflections of $L$.

\rom2
We have $T_{0,-\theta}\in\shr L$
by \rom1. Let $F$ be the transformation composed of
$T_{0,-\theta}$ followed by $(-\w)$-reflection in the short root
$(0;1,-\w)$. It is obvious that $F$ acts trivially on $\Lambda$
and computation reveals that it acts on $\rII1,1$ by left
multiplication by the matrix
$$ 
\pmatrix{0&\wbar\cr\wbar&0\cr}\;.
$$
The square of this matrix is the scalar $\w$ of $\II1,1$, which
proves the claim.

\rom3
Since $\shr L$ contains the central involution of $\Lambda$, $G$
contains the central involution $J$ of $\II1,1$. The biflection $B$
in $b=(0,\ldots,0;1,1)$ acts trivially on $\Lambda$ and on
$\II1,1$ as $\smallmatrix{0}{-1}{-1}{0}$. One can check that
$J=F^3B$, where $F\in\shr L$ is as in \rom2. This proves that
$B\in G$ and also that $G=\spanof{\shr L,B}$, hence
$G\sset\reflec L$. We also note that since $G$ contains $J$ and
also $F^2$, it contains all the scalars of $\II1,1$, so  it
suffices to show that $G$ contains the full stabilizer in $\aut
L$ of $\rho$. In light of \rom1 it suffices to merely show that $G$
contains $\aut \Lambda$. 

If $n=1$ then $\aut \Lambda$ is generated by reflections in its
short roots, as desired.  If $n>1$ then it suffices to prove
that $G$ contains the coordinate permutations with respect to
the chosen basis of $\Lambda$. That is, we must show that $G$
contains the biflections in vectors like
$x=(1,-1,0,\ldots,0;0,0)$.  It suffices to show that $x$ and $b$
are equivalent under $G$. To see this, observe that
$T_{(1,0,\ldots,0),\theta/2}$ followed by $F$, followed by the
scalar $-\w$, followed by $T_{(-1,1,0,\ldots,0),0}$, carries $x$
to $b$.
\endproof

\remark{Remarks:}
The 
condition $n>0$ is necessary; one can show
that $\reflec \eII1,1$ contains no scalars except the
identity. 

\beginproclaim Lemma
\Tag{thm-braiding-reflections}.
If $r$ and $r'$ are short roots in a lattice  over
$\ring=\Eisenstein$ or $\Hurwitz$ and $|\ip{r}{r'}|=1$, then $r$
and $r'$ are equivalent under the group generated by the
reflections in them.\endproclaim

\beginproof{Proof:}
One checks that the $(-\w)$-reflections $R$ and $R'$ in $r$ and
$r'$ satisfy the braid relation $RR'R=R'RR'$. (Because the
Hermitian form is degenerate on the span of $r$ and $r'$, one
must check that this relation holds by using
Eq.~\eqtag{eq-def-of-reflection}, not by just multiplying matrices for
the actions of $R$ 
and $R'$ on the span of $r$ and $r'$.) Rewriting this as
$R'^{-1}RR'=RR'R^{-1}$ we see that $R$ and $R'$ are conjugate in
the group they generate, which implies the lemma.
\endproof

\remark{Remark:}
The proof suggests connections between the braid groups and
complex reflection groups. This connection was first observed by Coxeter
\cite{coxeter:quotients-of-braid-groups}, and the braid groups
play a central role in the work of Deligne and Mostow
\cite{deligne:monodromy-of-hypergeometric-functions}, Mostow
\cite{mostow:picard-lattices-from-half-inegral-conditions} and
Thurston \cite{thurston:shapes-of-polyhedra}. They are also
important in work of the author, J. Carlson and D. Toledo
\cite{allcock:ch4-cubic-moduli}, 
\cite{allcock:ch4-cubic-moduli-announcement} on moduli of cubic surfaces.

\beginproclaim Theorem
\Tag{thm-eI2,1-eI3,1-eI4,1}.
Let $\ring=\Eisenstein$, $\Lambda=\Eisenstein^n$ and
$L=\Lambda\oplus\eII1,1$.
\item{\rom1} 
If $n=1$ then $\shr L$ acts with exactly $2$ orbits
of primitive null vectors, represented by $\pm\rho$. If $n=2$ or
$3$ then $\shr L$ acts transitively on the primitive null vectors of $L$.
\item{\rom2} 
If $n=1$, $2$ or $3$ then 
$\aut L=\reflec L=\shr L\times\{\pm I\}$.
\endproclaim 

\beginproof{Proof:}
First we show that $\shr L$ acts transitively on the
1-dimensional primitive null lattices in $L$.
For $n=1$ or $2$ this follows from
Thm.~\tag{thm-examples-of-reflective-lattices}. So suppose
$n=3$ and that $v\in L$ is a
primitive null vector not proportional to $\rho$ and of smallest
height in its orbit under $\shr L$. Since the covering radius of
$\Eisenstein^3$ is $1$,
Lemma~\tag{thm-short-root-height-reduction}\rom2 implies that $v$ is 
orthogonal to a  short root of height $1$. By applying a
translation (Lemma~\tag{thm-basics-for-eisenstein-lattices}\rom1) we may suppose that this root is
$r_1=(0,0,0;1,-\w)$. Taking $r_2=(0,0,1;0,1)$ and
$r_3=(0,0,1;0,0)$ we see that
$\ip{r_1}{r_2}=\ip{r_2}{r_3}=1$, so 
Lemma~\tag{thm-braiding-reflections} shows that $r_1$ 
is equivalent to 
$r_3$ under $\shr L$. Thus $v$ is equivalent to an element of
$r_3^\perp$, which is a copy of
$\Eisenstein^2\oplus\eII1,1$. Applying the $n=2$ case, we see
that $\shr L$ acts transitively on the 
primitive null sublattices of $L$.

It follows from
Lemma~\tag{thm-basics-for-eisenstein-lattices}\rom3 that $\aut
L$ is generated by $\shr L$ and $\{\pm I\}$. Now we will show
that $-I\notin\shr L$, which will establish the equality $\aut
L=\shr L\times\{\pm I\}$. Since $-I\in\reflec L$ by
Lemma~\tag{thm-basics-for-eisenstein-lattices}\rom3 we will have
proven \rom2. To prove $-I\notin\shr L$ we must consider the
finite vector space $V=L/L\theta$ over
$\F_3=\Eisenstein/\theta\Eisenstein$; we write $q$ for both
natural maps $L\to V$ and $\Eisenstein\to\F_3$. The Hermitian
form on $L$ gives rise to a symmetric bilinear form on $V$,
given by $\ip{q(v)}{q(w)}=q(\ip{v}{w})$. Since
$\eII1,1\isomorphism\eI1,1$, $L$ and hence $V$ admits an
orthogonal basis with $n+1$ vectors of norm $1$ and one of norm
$-1$. There is a homomorphism called the spinor norm from $\aut
V$ to the set $\{\pm1\}$ of nonzero square classes in
$\F_3$. This is characterized by the property that the
reflection in a vector of $V$ of norm $\pm1$ has spinor norm
$\pm1$. It is clear that $-I$ acts on $V$ with spinor norm
$-1$. Since a reflection of $L$ in a short root $r$ acts on $V$
either trivially or as the reflection in $q(r)$, of norm
$q(r^2)=1$, every element of $\shr L$ acts on $V$ with spinor
norm $+1$. Hence $-I\notin\shr L$, as desired. This also
characterizes $\shr L$ as the subgroup of $\aut L$ whose
elements act on $V$ with spinor norm $+1$.

Now we will establish \rom1. By the first part of the proof it
suffices to determine which multiples of $\rho$ are equivalent
to each other. By
Lemma~\tag{thm-basics-for-eisenstein-lattices}\rom2 it suffices
to determine whether $\pm\rho$ are equivalent. If $n=1$ then
they are inequivalent, because the stabilizers of $\rho$ in
$\shr L$ and $\aut L$ coincide. If $\rho$ and $-\rho$ were
equivalent then we would have transitivity on primitive null
vectors and  $\shr L=\aut L$ would follow. Since this is not
true, $\rho$ and $-\rho$ are inequivalent. On the other hand, if
$n=2$ or $3$, then $\pm\rho$ are equivalent. To see this, apply
the product of $-I$ and the biflection in any long root of $L$
that is orthogonal to $\rho$. The product exchanges $\pm\rho$,
and by spinor norm considerations it lies in $\shr L$.
\endproof

\remark{Remark:}
The case $n=3$ arises in algebraic geometry: the quotient of
$\ch4$ by $\shr \eI4,1$ may be identified with the moduli space
of stable cubic surfaces in $\cp3$.  One can also construct the
moduli space of marked stable cubic surfaces by
taking the quotient of $\ch4$ by the congruence subgroup of
$\shr
\eI4,1$ associated to the prime $\theta\in\Eisenstein$. The
quotient of $\shr \eI4,1$ by this normal subgroup is the $E_6$
Weyl group, also known as ``the group of the 27 lines on a cubic
surface''.  See \cite{allcock:ch4-cubic-moduli} for details.

\beginproclaim Lemma 
{\Tag{thm-basics-for-Hurwitz-lattices}}. 
Suppose $\ring=\Hurwitz$, $\Lambda=\Hurwitz^n$ $(n>0)$ and
$L=\Lambda\oplus\II1,1$. Then
\item{\rom1}
$\shr L$ contains a translation $T_{x,z}$ for each $x\in
\Lambda$, and also the central translations $T_{0,a\i+b\j+c\k}$
with $a\congruent b\congruent c$ $(\mod~2)$. In particular, coset
representatives for the translations of $\shr L$ in those of
$\aut L$ may be taken from
$\{T_{0,0},T_{0,\i},T_{0,\j},T_{0,\k}\}$. 
\item{\rom2}
$\shr L$ contains transformations acting trivially on $\Lambda$
and on $\II1,1$ by left scalar multiplication by any given
unit of $\Hurwitz$
\item{\rom3}
The stabilizer of $\spanof{\rho}$ in $\shr L$ has index~$\leq4$
in the stabilizer in $\aut L$; coset representatives may be
taken from the set given in \rom1.
\par
\endproclaim

\beginproof{Proof:}
\rom1
The first part follows immediately from
Lemma~\tag{thm-stabilizer-of-rho}\rom1. The second part may be
obtained by taking commutators: if $\lambda,\lambda'\in
\Lambda$ and $z,z'\in\im\K$ are such that
$T_{\lambda,z},T_{\lambda',z'}\in\shr L$, then
$T_{0,\pm2\im\ip{\lambda}{\lambda'}}\in\shr L$ by
Eq.~\eqtag{eq-commutator-of-translations}. Since $\Lambda$ contains
vectors $\lambda$ and $\lambda'$ with
$\ip{\lambda}{\lambda'}=\a$ for any given unit $\a$ of
$\Hurwitz$, we see that $\shr L$ contains $T_{0,2\i}$,
$T_{0,2\j}$, $T_{0,2\k}$ and $T_{0,\i+\j+\k}$. These generate
the group of central translations given in the statement of the
lemma.

\rom2
The argument of Lemma~\tag{thm-basics-for-eisenstein-lattices}\rom2 shows that $\shr L$ contains
an element acting trivially on $\Lambda$ and on $\II1,1$ as
left-multiplication by $\w$. Taking conjugates of this by the
group $\aut\Hurwitz$ acting on $\II1,1$, which normalizes $\shr
L$ even though it doesn't act $\Hurwitz$-linearly, we see that
$\shr L$ contains elements acting on $\II1,1$ as
left-multiplication by any of the units $(-1\pm\i\pm\j\pm\k)/2$ of
$\Hurwitz$. These generate the group of all units of $\Hurwitz$,
proving \rom2.

\rom3
Follows immediately from \rom1 and \rom2 and the arguments given
for Lemma~\tag{thm-basics-for-eisenstein-lattices}\rom3. Note
the curious fact that $\shr L$ contains the biflection $B$,
which it did not in the Eisenstein case.
\endproof

\beginproclaim Theorem
\Tag{thm-qI2,1-qI3,1}.
Let $\ring=\Hurwitz$, $\Lambda=\Hurwitz^n$ for  $n=1$ or $2$, and
$L=\Lambda\oplus\qII1,1$. Then $\shr L$ acts transitively on the
primitive null vectors of $L$ and has index at most $4$ in
$\aut L$; coset representative may be taken from 
$\{T_{0,0},T_{0,\i},T_{0,\j},T_{0,\k}\}$.
\endproclaim 

\beginproof{Proof:}
We first claim that $\shr L$ acts transitively on primitive null
lattices in $L$. For $n=1$ this follows from
Thm.~\tag{thm-examples-of-reflective-lattices}.  For $n=2$ it
follows from an argument similar to the $n=3$ case of
Thm.~\tag{thm-eI2,1-eI3,1-eI4,1}.  That is, the covering radius
of $\Lambda=\Hurwitz^2$ is $1$, so if $v\in L$ is primitive,
isotropic and of smallest height in its orbit under $\shr L$
then by Lemma~\tag{thm-short-root-height-reduction}\rom2 we see
that $v$ is either proportional to $\rho$ or orthogonal to a
short root $r_1$ of height $1$. In the latter case, after
applying a translation of $\shr L$, courtesy of
Thm.~\tag{thm-basics-for-Hurwitz-lattices}\rom1, we may take
$r=(0,0;1,x-\w)$, where $x$ is one of $0$, $\i$, $\j$ and
$\k$. In any of these cases, upon taking $r_2=(0,1;0,1)$ and
$r_3=(0,1;0,0)$ we have $\ip{r_1}{r_2}=\ip{r_2}{r_3}=1$. By
Lemma~\tag{thm-braiding-reflections}, $v$ is equivalent under
$\shr L$ to an element of $r_3^\perp$. Since $r_3^\perp$ is a
copy of $\Hurwitz^1\oplus\qII1,1$, the transitivity follows from
the case $n=1$.

The transitivity on primitive null vectors follows from
Thm.~\tag{thm-basics-for-Hurwitz-lattices}\rom2. The rest of the
theorem follows from
Lemma~\tag{thm-basics-for-Hurwitz-lattices}\rom3.
\endproof

Now we move on to higher-dimensional examples\emdash we will 
construct a group acting on $\ch7$ and 
another acting on $\qh5$. These arise  from our basic
construction by taking $\Lambda=\cct$ or
$\qbwl$. 

\beginproclaim Lemma
\Tag{thm-linear-independence-descends}.
Suppose $v,r_1,\ldots,r_m\in\K^n\oplus\K^{1,1}$ lie
over $\ell,\lambda_1,\ldots,\lambda_m\in\K^n$ respectively. Suppose
$v^2=0$, that $\ip{r_i}{v}=0$ for all $i$, and that the vectors
$\lambda_i-\ell$ are linearly independent in $\K^n$. Then the
images of the $r_i$ in $v^\perp/\spanof{v}$ are linearly
independent. \endproclaim

\beginproof{Proof:}
We may obviously replace $v$ and the $r_i$ by any scalar
multiples of themselves and so suppose that they have height
$1$. Thus $v=(\ell;1,?)$ and $r_i=(\lambda_i;1,?)$ where the
question marks denote irrelevant (and possibly distinct) elements
of $\K$. Let $T$ be the translation carrying $v$ to $(0;1,0)$,
so $T(r_i)=(\lambda_i-\ell;1,0)$. (The last coordinate vanishes
because $\ip{T(r_i)}{Tv}=0$.) Since the image of $T(r_i)$ in
$(Tv)^\perp/\spanof{Tv}$ may be identified with its first
coordinate, namely $\lambda_i-\ell$, the images of the $T(r_i)$ in
$(Tv)^\perp/\spanof{Tv}$ are linearly independent. The lemma
follows immediately.
\endproof

\beginproclaim Lemma
\Tag{thm-coxeter-todd-well-covered}.
$\cct\tensor\R$ is covered by the closed balls of radius $1$ centered
at points of $\cct$, together with those of radius
$1/\sqrt3$ centered at points $\lambda\theta^{-1}$ with
$\lambda\in\cct$ and $\lambda^2\congruent1(\mod~3)$. 
\endproclaim

\beginproof{Proof:}
Section~7 of \cite{jhc:Coxeter-Todd} defines a linear ``gluing map'' 
$g:{\textstyle{1\over\theta}}\cct/\cct\to{\textstyle{1\over\theta}}\cct/\cct$
with the property that the Leech lattice $\leech$, scaled
down by $2^{1/2}$, is the real form of the lattice of  vectors 
$(x_1,x_2)\in\({1\over\theta}\cct\)^2$ satisfying 
$g(x_1+\cct)=x_2+\cct$. Identifying $\cct$ with the set of such
$(x_1,x_2)$ with 
$x_2=0$, we see that the only points of $2^{-1/2}\leech$ 
at distance $<1$ from $\cct\tensor\R$ are those in $\cct$ and
those of the form $(x_1\theta^{-1},x_2\theta^{-1})$ with $x_2$ a
minimal vector of $\cct$ (a long root). 
The definition of $g$ (see \cite{jhc:Coxeter-Todd}) shows that
$x_1^2\congruent1(\mod~3)$ if and only if there is a long root $x_2$
of $\cct$ such that 
$(x_1\theta^{-1},x_2\theta^{-1})\in2^{-1/2}\leech$.
By \cite{jhc:leech-radius}, the covering radius of
$2^{-1/2}\leech$ is $1$. Therefore the balls
of radius $1$ centered at the points of $\cct$ and at the points
$(x_1\theta^{-1},x_2\theta^{-1})$ with $x_1^2\congruent1(\mod~3)$ and
$x_2^2=2$ cover 
$\cct\tensor\R$. Computing the radius of the intersection of a
ball of the second family with $\cct\tensor\R$ yields the lemma.
\endproof

\beginproclaim Theorem
\Tag{thm-eI7,1}.
Let $\Lambda=\cct$ and $L=\Lambda\oplus\eII1,1\isomorphism\eI7,1$. 
\item{\rom1} If $v\in L$ is primitive, isotropic and not
equivalent to a multiple of $\rho$ under $\reflec_0 L$, then
$v^\perp/\spanof{v}\isomorphism\Eisenstein^6$.
\item{\rom2}  $\aut L=\reflec L$. In particular, $L$ is
reflective. \endproclaim 

\beginproof{Proof:}
\rom1 
Suppose that $v$ is a primitive isotropic vector of smallest
height in its orbit under $\reflec_0 L$. Suppose $v$ is not a
multiple of $\rho$, so that it lies over some $\ell\in
\Lambda\tensor \R$. By
Lemma~\tag{thm-short-root-height-reduction} and the minimality
of the height of $v$, $\ell$ lies at distance $\geq1$ from
each lattice point $\lambda\in \Lambda$ and at distance $\geq
3^{-1/2}$ from each $\lambda\theta^{-1}$ with $\lambda\in
\Lambda$ and $\lambda^2\congruent1(\mod~3)$. By
Lemma~\tag{thm-coxeter-todd-well-covered} the set $S$ of such
points in $\Lambda\tensor\R$ is discrete. Let
$\mu_1,\ldots,\mu_n$ be the elements of $\Lambda$ with
$(\ell-\mu_i)^2=1$ and let $\nu_1,\ldots,\nu_m$ be those vectors
of the form $\lambda\theta^{-1}$ with $\lambda\in \Lambda$ and
$\lambda^2\congruent1(\mod~3)$ such that
$(\ell-\nu_i)^2=1/3$. Over each $\mu_i$ (resp. $\nu_i$) there is
a short root of $L$, say $r_i$ (resp. $s_i$), of height $1$
(resp. $\theta$). By
Lemma~\tag{thm-short-root-height-reduction}, we may suppose that
the $r_i$ and $s_i$ are orthogonal to $v$. Because $S$ is
discrete the vectors $\mu_i-\ell$ and $\nu_i-\ell$ span
$\Lambda\tensor\R$, and therefore there are 6 among them that
are linearly independent over $\C$. By
Lemma~\tag{thm-linear-independence-descends} this implies that
among the images in $v^\perp/\spanof{v}$ of the vectors $r_i$
and $s_i$ are $6$ short roots that are linearly independent over
$\Eisenstein$. Since $v^\perp/\spanof{v}$ is positive-definite,
it follows that $v^\perp/\spanof{v}\isomorphism
\Eisenstein^6$. 

\rom2
We have seem the transitivity of $\shr L$ on the primitive null
sublattices that, like $\spanof{\rho}$, are orthogonal to no short
roots. Since
$L\isomorphism\eI7,1\isomorphism\Eisenstein^6\oplus\II1,1$,
Lemma~\tag{thm-basics-for-eisenstein-lattices}\rom3 implies that $\spanof{\shr L, -I}$ lies in
$\reflec L$ and contains the scalars of $L$. Therefore $\reflec
L$ acts transitively on the primitive null vectors that, like
$\rho$, are orthogonal to no short roots. Now it suffices to show
that $\reflec L$ contains the full stabilizer of $\rho$. Since
$\Lambda$ is selfdual and spanned by its long roots,
Lemma~\tag{thm-stabilizer-of-rho}\rom2 and~\rom3 imply that $\reflec L$ contains all
the translations of $L$. Since $\aut \Lambda$ is a reflection
group, the proof is complete.
\endproof

\remark{Remark:}
it would be nice to understand the groups $\shr\eI7,1$ and
$\shr\qI5,1$. 

We now study the quaternionic lattice $\qI5,1$. The analysis is
surprisingly similar to our study of $\eI7,1$.  In particular,
Lemma~\tag{thm-barnes-wall-well-covered} below is very similar
to Lemma~\tag{thm-coxeter-todd-well-covered}.  Section~5 of
\cite{jhc:laminated-lattices} describes an embedding of the real
form $\rbwl$ of $2^{1/2}\qbwl$ into the Leech lattice $\leech$.
When we refer to concepts involving $\leech$ while discussing
$\qbwl$, we implicitly refer to this embedding.  (Up to isometry
of $\leech$, there is only one embedding.)

\beginproclaim Lemma
\Tag{thm-barnes-wall-well-covered}.
$\qbwl\tensor\R$ is covered by the closed balls of radius $1$
centered at points of $\qbwl$, together with those of radius
$1/\sqrt2$ centered at points $\lambda(1+\i )^{-1}$ with
$\lambda\in\qbwl$ of odd norm.  Any point of $\qbwl\tensor\R$
not in the interior of one of these balls is a deep hole of
$2^{-1/2}\leech$.\endproclaim

\beginproof{Proof:}
The orthogonal complement of $\qbwl$ in $2^{-1/2}\leech$ is
a copy of the $E_8$ lattice. Properties of the embedding are
described in \cite{jhc:laminated-lattices} and include the following.
If 
$(x,y)\in(\qbwl\tensor\R)\times(E_8\tensor\R)$ lies in
$2^{-1/2}\leech$, then $y\in{1\over2}E_8$ and hence has
norm $n/2$ for some nonnegative integer $n$. 
We write $B(x,y)$ for the ball of radius 1 with center
$(x,y)\in2^{-1/2}\leech$. Only if the norm of $y$ is $0$ or
$1/2$ does the interior of $B(x,y)$ meet $\qbwl\tensor\R$. Those
$x$ for which $(x,y)\in2^{-1/2}\leech$ for some $y$ of norm
$1/2$ are exactly the deep
holes of $\qbwl$. By Thm.~\tag{thm-barnes-wall-properties},
the set of such $x$ 
coincides exactly with 
$\set{\lambda(1+\i )^{-1}}{\lambda\in\qbwl,\lambda^2\congruent1(\mod~2)}$.
For such $(x,y)$, the ball $B(x,y)$ meets $\qbwl\tensor\R$ in a
ball of radius $2^{-1/2}$. The theorem follows from the fact
that the covering radius of
$2^{-1/2}\leech$ is 1 (see \cite{jhc:leech-radius}). 
\endproof

\beginproclaim Lemma 
{\Tag{thm-leech-holes-in-BW}}. 
Any deep hole of $\leech$ that lies in $\rbwl\tensor\R$ has a
vertex in $\rbwl$.
\endproclaim

\beginproof{Proof:}
\def\D{\Delta}
The natural language for discussing the deep holes of $\leech$
is that of affine Coxeter-Dynkin diagrams, using the slightly
nonstandard conventions of \cite{jhc:leech-radius}. If $h$ is a deep hole of
$\leech$ then its vertices $v_i$ lie at distance $\sqrt2$ from
$h$, and we define the diagram $\D$ of $h$ to be the graph whose
vertices are the $v_i$, with $v_i$ and $v_j$ unjoined, singly
joined or doubly joined according to whether $(v_i-v_j)^2$ is
$4$, $6$ or $8$. Each component of $\D$ is an affine diagram of
type $A_n$, $D_n$ or $E_n$. For the rest of the proof we will
take $h$ as the origin. The definition of $\D$ and the fact that
$(v_i-h)^2=2$ for all $i$ means that the inner product of $v_i$
and $v_j$ is $0$, $-1$ or $-2$ according to whether the
corresponding vertices of $\D$ are unjoined, singly joined, or
doubly joined. It follows that the subspaces spanned by
different components of $\D$ are orthogonal, and that the
vertices corresponding to each component form a system of simple
roots of the corresponding type, together with the lowest root,
which corresponds to the extending node in the diagram. Now,
$\rbwl$ is the fixed-point set of an involution $\phi$ of
$\leech$. Since $\phi$ fixes $h$, it acts on $\D$. We will show
that $\phi$ preserves a vertex $v$ of $\D$, which forces $v$ to
lie in $\rbwl$, which proves the lemma.

For each component $C$ of $\D$ we write $S_C$ for the
corresponding subspace of $\R^{24}$. If $\phi$ preserves $C$
then we write $F_C$ for the subspace of $S_C$ fixed pointwise by
$\phi$. We write $F$ for the subspace of $\R^{24}$ fixed
pointwise by $\phi$. It is easy to see that
$$
\dimension F=\sum_{\phi(C)=C}\dimension F_C
+\sum_{\phi(C)\neq C}{\dimension S_C\over2}\;,
\eqno\eqTag{eq-dimension-counting}
$$
where the sums are over the components $C$ of $\D$ that are
(resp. are not) preserved by $\phi$. If $\phi$ preserves $C$
then its action on $C$ determines $F_C$. The explicit
description in terms of root systems allows one to deduce that
$\dimension F_C$ equals the number of vertices of $C$ preserved
by $\phi$, plus half the number not preserved, minus 1. It
follows that if $\phi$ permutes the $v_i$ freely then each term
in each sum in Eq.~\eqtag{eq-dimension-counting} is bounded by
${1\over2}\dimension S_C$. Since $\sum_C \dimension S_C=24$ we
would obtain $\dimension F\leq 12$, which is impossible since
$\rbwl\sset F$.
\endproof

\remark{Remark:}
A more involved analysis shows that any deep hole of $\leech$
lying in $\rbwl\tensor\R$ has at least nine vertices in $\rbwl$,
and that this bound cannot be improved.

\beginproclaim Theorem
\Tag{thm-qI5,1}.
Let $\Lambda=\qbwl$ and $L=\Lambda\oplus\qII1,1\isomorphism\qI5,1$. 
\item{\rom1} 
If $v\in L$ is a primitive isotropic vector not
equivalent under  $\shr L$ to a multiple of $\rho$ then $v^\perp/\spanof{v}$
contains a short root.
\item{\rom2} 
The index of  $\reflec L$ in $\aut L$ is at most $4$, so that
$L$ is reflective. More precisely, coset representatives for
$\reflec L$ in $\aut L$ may be taken from
$\{T_{0,0},T_{0,\i},T_{0,\j},T_{0,\k}\}$. 
\par
\endproclaim

\beginproof{Proof:}
\rom1
This is very similar to the proof of
Thm.~\tag{thm-eI7,1}\rom1. Suppose 
that $v$ is a primitive isotropic vector of $L$ of smallest
height in its orbit under $\reflec L$. Suppose $v$ is not a multiple
of $\rho$, so that $v$ lies over some $\ell\in
\Lambda\tensor\R$. By
Lemma~\tag{thm-short-root-height-reduction}, $\ell$ lies at 
distance $\geq1$ from each lattice point $\lambda\in \Lambda$
and at distance $\geq 2^{-1/2}$ from each point
$\lambda(1+\i )^{-1}$ with $\lambda\in \Lambda$ and
$\lambda^2\congruent1(\mod~2)$. By Lemma~\tag{thm-barnes-wall-well-covered},
$\ell$ must be a deep 
hole of $2^{-1/2}\leech$. By
Lemma~\tag{thm-leech-holes-in-BW} there 
is a vertex $\lambda\in\qbwl$ of the hole with
$(\ell-\lambda)^2=1$. There is a short root of $L$ lying
over $\lambda$, and by Lemma~\tag{thm-short-root-height-reduction} there is also one
orthogonal to $v$.

\rom2
Since $\Lambda$ is selfdual and spanned by its long roots,
Lemma~\tag{thm-stabilizer-of-rho}\rom2 shows that $\reflec L$ contains a translation
$T_{x,z}$ for each $x\in \Lambda$. Taking commutators as in
Lemma~\tag{thm-basics-for-Hurwitz-lattices}\rom1 shows that $\reflec L$ contains the central
translations $T_{0,a\i+b\j+c\k}$ with $a\congruent b\congruent
c$ (mod~$2$). Then the proof of Lemma~\tag{thm-basics-for-Hurwitz-lattices}\rom2 shows that
$\reflec L$ contains elements acting on $\II1,1$ by
left-multiplication by the units of $\Hurwitz$. Together
with \rom1 this proves the transitivity of $\reflec L$ on
primitive null vectors that, like $\rho$. are orthogonal to no
short roots. Then \rom2 follows from the facts that $\aut
\Lambda$ is a reflection group and $\reflec L$ contains the
translations just discussed.
\endproof

We close this section by returning to low dimensions, studying
the 2-dimensional selfdual Lorentzian lattices. If
$\ring=\Eisenstein$ or $\Gaussian$ then one can obtain very
explicit descriptions of the groups by drawing pictures in
$\ch1\sset\cp1$. In particular, if we represent a point
$(a,b)\in\II1,1$ by $a/b\in\cp1$ then the hyperbolic space
becomes the right half-plane and $\rho$ the point at
infinity. It is easy to find the points of $\cp1$ corresponding
to the roots of $L$ of small height, and then one can work out the
group $\reflec L$. For example, one can check that
$\reflec\eII1,1$ acts
as the triangle group $\trianglegroup{2}{6}{\infty}$. 
One can also show that $\aut\gII1,1$ acts on $\ch1$
as $\trianglegroup{2}{3}{\infty}$ and its subgroup
of index 2 consisting of elements with determinant $+1$ is conjugate in
$\GL_2(\gauss)$ to $\SL_2\Z$. 
The group  $\reflec\gII1,1$ 
%has index 6 in $\aut\gII1,1$ and
is generated by 3 biflections, which act by rotations by $\pi$
around the three finite corners of a quadrilateral in $\ch1$
with corner angles $\pi/2$, $\pi/2$, $\pi/2$ and $\pi/\infty$.
For completeness we mention that $\aut\gI1,1$ acts on $\ch1$ as
$\trianglegroup{2}{4}{\infty}$, and its reflection subgroup acts
as $\trianglegroup{4}{4}{\infty}$.
See \cite{cox:gens-and-rlns} for descriptions of the
groups $\trianglegroup{p}{q}{r}$ and other information.

One can also treat the quaternionic case: an adaptation of the
argument of \ecite{allcock:oh2}{Thm.~5.3\rom1} shows that
$\reflec\qII1,1$ acts on $\qh1\isomorphism\rh4$ as the rotation
subgroup of the real hyperbolic reflection group with the
Coxeter diagram below.
Note that the 6 outer nodes generate an affine reflection
group, so this graph is a special case of
the usual procedure of ``hyperbolizing'' an affine reflection
group by adjoining an extra node.
$$
%
%  a Coxeter Diagram
%
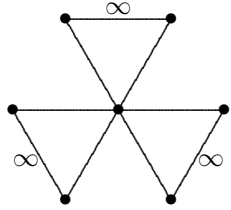
\beginpicture
\def\magfactor{2}
\newdimen\vstep
\vstep=17.13 pt
\newdimen\hstep
\hstep=20 pt
\multiply\vstep by\magfactor
\multiply\hstep by\magfactor
\setcoordinatesystem units <\hstep,\vstep>
\put {$\bullet$} at  .5 -1 
\put {$\bullet$} at  -.5 -1 
\put {$\bullet$} at  1 0
\put {$\bullet$} at  0 0  
\put {$\bullet$} at  -1 0 
\put {$\bullet$} at  .5 1  
\put {$\bullet$} at  -.5 1   
\plot 0 0 .5 1 -.5 1 0 0 -1 0 -.5 -1 0 0 .5 -1 1 0 0 0 /
\put {$\infty$} [b]  <0pt,2pt> at 0 1 
\put {$\infty$} [tr] <0pt,0pt> at -.75 -.5
\put {$\infty$} [tl] <0pt,0pt> at  .75 -.5
\endpicture
$$

\section{\Tag{sec-selfdual-classification}. Enumeration of
selfdual lattices}

As we explain below, the orbits of primitive isotropic lattices
in the Lorentzian 
lattice $\rI n+1,1$ are in natural $1$-$1$ correspondence with
the equivalence classes of positive-definite selfdual lattices of
dimension $n$ over $\ring$. This means that one may classify
such lattices by studying $\aut\rI n+1,1$. Since we have made
such a study in the previous section, in terms of the geometry
of various special lattices, we can now 
classify selfdual lattices in low dimensions. We begin with
an analogue of a  result well-known for lattices over $\Z$.

\beginproclaim Theorem
\Tag{thm-indefinite-selfdual-classification}.
An indefinite selfdual lattice $L$ over $\ring=\Eisenstein$ or
$\Hurwitz$ is characterized up to isometry by its dimension and
signature. An indefinite selfdual lattice $L$ over
$\ring=\Gaussian$ is characterized up to isometry by its
dimension, signature, and whether it is even; if $L$ is even
with signature $(n,m)$ then $n-m$ is divisible by $4$.
\endproclaim

\beginproof{Proof:}
First we show that any indefinite selfdual $\ring$-lattice $L$
contains an isotropic vector. If $\ring=\Hurwitz$, or if
$\dimension L>2$, then the real form of $L\tensor\Q$ is an
indefinite rational bilinear form of rank $>4$, so Meyer's
theorem
\ecite{milnor:bilinear-forms}{Chap.~2} asserts the existence of an
isotropic vector. If $\dimension L=2$ and $\ring=\gauss$ or
$\eisen$, then we consider the $2\times2$ matrix of inner
products of the elements of a basis for $L$. This may be
diagonalized by row and column operations
over $\ring\tensor\Q$ to a diagonal matrix
$[a,-a^{-1}]$ with $a\in\Q$. (Each term is real because the
matrix is Hermitian, and each determines the other because the
determinant is $-1$.)  Then  the vector $(1,a)$ is
isotropic. Having obtained an isotropic vector in
$L\tensor \Q$, we may multiply by a scalar to obtain one in $L$. 

If $L$ is odd then the proof of  Thm.~4.3 in
\ecite{milnor:bilinear-forms}{Chap.~2} applies, and $L\isomorphism
\rI n,m$ for some $n$ and $m$. This completes the proof of the
first claim, since any selfdual lattice over $\Eisenstein$ or
$\Hurwitz$ is odd: if $v,w\in L$ satisfy $\ip{v}{w}=\w$ then
$v^2$, $w^2$ and $(v+w)^2$ cannot all have the same
parity. This also proves that an odd indefinite
selfdual Gaussian lattice is characterized by its dimension
and signature.

One may construct lattices $N$ from an odd Gaussian lattice $M$ by
considering the sublattice $M^e$ consisting of the elements  of
$M$ of even norm, and considering the 3 lattices $N$ such that
$M_e'\sset N\sset M_e$. When $M$ is $\gI1,1$, then
$N$ may be chosen to be $\gII1,1$.
Now consider an indefinite 
even selfdual $\gauss$-lattice $L$. We know that $L$ contains
an isotropic vector, and as in \cite{milnor:bilinear-forms} there
is a  decomposition $L=\Lambda\oplus\gII1,1$. We see that $L$
arises by applying the construction above to the odd selfdual lattice
$\Lambda\oplus\gI1,1$. Since  $\Lambda\oplus\gI1,1$
is isomorphic to $\gI n,m$, it is clear
that all possible $L$ can be constructed by applying our
construction to the various $\gI n,m$. No even lattices arise unless
$n-m\congruent0(\mod~4)$, when two isometric ones do.
\endproof

Special cases of
Thm.~\tag{thm-indefinite-selfdual-classification}  are
$\eI7,1\isomorphism\cct\oplus\eII1,1$ and
$\qI5,1\isomorphism\qbwl\oplus\qII1,1$, which are the lattices
studied in
Thms.~\tag{thm-eI7,1} and
\tag{thm-qI5,1}. 
Thm.~\tag{thm-indefinite-selfdual-classification} also
provides the correspondence mentioned above: if $V$ is a
primitive isotropic lattice in $\rI n+1,1$ then it is easy to
check that $V^\perp/V$ is an $n$-dimensional positive-definite
selfdual lattice, and that this establishes a one-to-one
correspondence between orbits of primitive isotropic lattices of
$\rI n+1,1$ and isometry classes of selfdual positive-definite
lattices in dimension $n$. Similarly, the orbits of primitive
isotropic lattices of $\gII n+1,1$ correspond to the classes of
positive-definite even selfdual Gaussian lattices of dimension
$n$.

\beginproclaim Theorem 
{\Tag{thm-classification-of-selfdual-lattices}}. 
The positive-definite selfdual $\Eisenstein$-lattices in
dimensions~$\leq6$ are $\Eisenstein^n$ and $\cct$. The
positive-definite selfdual $\Hurwitz$-lattices in
dimensions~$\leq4$ are $\Hurwitz^n$ and $\qbwl$. The 
positive-definite even selfdual $\Gaussian$-lattices in
dimensions~$\leq4$ are $\{0\}$ and $\geeight$.
\endproclaim 

\beginproof{Proof:}
By Thm.~\tag{thm-eI7,1}\rom1, any primitive null vector $v$ of
$\eI7,1$ satisfies $v^\perp/\spanof{v}\isomorphism\cct$ or
$v^\perp/\spanof{v}\isomorphism\Eisenstein^6$; the first claim
follows immediately. To see the last claim, suppose that
$\Lambda$ is 
an even selfdual $\Gaussian$-lattice of dimension~$\leq4$. By
the signature condition, the dimension is either $0$ or $4$. 
In the latter case the isomorphism $\Lambda\isomorphism\geeight$
follows from the equivalence of any two primitive null vectors
in $\gII5,1=\geeight\oplus\II1,1$ (Thm.~\tag{thm-examples-of-reflective-lattices}).

We will only sketch the quaternionic case. By
Thm.~\tag{thm-qI5,1}\rom1, any 4-dimensional positive-definite
selfdual $\Hurwitz$-lattice $\Lambda$ besides $\qbwl$ has a
short root. We claim that in fact $\Lambda$ has a pair of
independent (hence orthogonal) short roots. This follows from
the remark after Lemma~\tag{thm-leech-holes-in-BW}: if $\ell$ is a
deep hole of $2^{-1/2}\leech$ lying in $\qbwl\tensor\R$ then it
has 9 vertices $v_i$ in $\qbwl$. By considering the hole diagram
of $\ell$ one can show that the $v_i-\ell$ span a space of real
dimension~$\geq5$, hence of dimension~$\geq2$ over $\H$. Then
the argument of Thm.~\tag{thm-qI5,1}\rom1 establishes the
claim. Therefore $\Lambda$ is the direct sum of $\Hurwitz^2$ and
a two-dimension selfdual $\Hurwitz$-lattice. The second summand
must also be $\Hurwitz^2$, by the treatment of $\qI3,1$ in
Thm.~\tag{thm-qI2,1-qI3,1}.
\endproof

These results have been obtained before, by very different
means. Feit \cite{feit:unimodular-Eisenstein} found examples of
many positive-definite selfdual $\Eisenstein$-lattices. He
derived a version of the mass formula to verify that his list
was complete for dimensions $n\leq12$. Conway and Sloane
\ecite{jhc:Coxeter-Todd}{Thm.~3} provide a nice proof of this
classification in dimensions $n\leq6$ based on theta series and
modular forms. (Their proof does not apply for $6<n<12$: in the
second-to-last sentence of the proof, ``12'' should be replaced
by ``7''.)  Although selfdual $\Gaussian$-lattices have not been
tabulated, it would be easy (and boring) to enumerate them
through dimension $12$ by using the fact that the real form of a
selfdual $\Gaussian$-lattice is selfdual over $\Z$.  An
enumeration of positive-definite selfdual $\Hurwitz$-lattices
for dimensions $n\leq7$ has recently been completed by Bachoc
\cite{bachoc:selfdual-Hurwitz-lattices} and for $n=8$ by Bachoc
and Nebe \cite{bachoc:8-dimensional-Hurwitz-lattices}. 
These enumerations are based on a generalization of Kneser's notion
of ``neighboring'' lattices, together with a suitable version of
the mass formula.

\section{\Tag{sec-comparison}. Comparison with the groups of
Deligne and Mostow}

In this section we justify the word ``new'' in our title, by
showing that our largest three reflection groups do not
appear on the lists of Mostow
\cite{mostow:monodromy-groups-on-the-complex-n-ball} and
Thurston 
\cite{thurston:shapes-of-polyhedra}. Deligne and Mostow
\cite{deligne:monodromy-of-hypergeometric-functions} and Mostow
\cite{mostow:picard-lattices-from-half-inegral-conditions} 
constructed 94 reflection groups acting on $\chn$ for various
$n=2,\ldots,9$ by considering the monodromy of hypergeometric
functions. Thurston \cite{thurston:shapes-of-polyhedra}
constructed the same set of groups 
in terms of moduli of flat metrics (with specified
sorts of singularities) on the sphere $S^2$. We will generally refer to
these groups as the DM groups. We show here (Thm.~\tag{thm-my-groups-are-new}) that none
of the groups $\reflec\eI n,1$ ($n\geq4$) or $\reflec\gII
4n+1,1$ ($n\geq1$) appear on their lists. In particular, our
groups $\reflec\eI7,1$, $\reflec\eI4,1$ and $\reflec\gII5,1$ are
new.  We will also identify 
$\reflec\eI3,1$ with one of the DM groups. We leave open
the question of whether our other groups appear on their
lists and also the question of commensurability.

We will distinguish our groups from the DM groups by considering
the orders of the reflections in the groups. We begin by showing
that the only reflections of the selfdual lattices are the
obvious ones, a result well-known for lattices over $\Z$.

\beginproclaim Lemma
\Tag{thm-classify-honest-reflections}.
Any reflection $R$ of a selfdual lattice $M$ over
$\ring=\Eisenstein$ or $\Gaussian$ is either a reflection in a
lattice vector of norm $\pm1$ or a biflection in a lattice
vector of norm $\pm2$.\endproclaim

\beginproof{Proof:}
By considering the determinant of $R$ we discover that its only
nontrivial eigenvalue is a unit of $\ring$, so $M$ contains an
element of the corresponding eigenspace, so $R$ is the
$\a$-reflection in some lattice vector $v$, where $\a$ is a unit
of $\ring$. Taking $v$ to be
primitive, every lattice vector in the complex span of $v$ 
lies in the $\ring$-span of $v$. (This uses the fact that
$\ring$ is a principal ideal domain.) Furthermore, by
the selfduality of $M$, there exists $w\in M$ satisfying
$\ip{v}{w}=1$. Then $R(w)=w-v(1-\a)/v^2$ and so 
$w-R(w)=v(1-\a)/v^2$ lies in $M$. Therefore
$(1-\a)/v^2\in\ring$. Unless $\a=-1$ this requires $v^2=\pm1$
and if $\a=-1$ then it requires that $v^2$ divide $2$.
\endproof

In order to compare our groups to the DM groups we will also
need to consider the transformations of projective space that
arise from linear reflections, which we call projective
reflections. If $L$ is a Lorentzian lattice then $\paut L$ may
contain projective reflections that are not represented by any
reflection of $L$. For an example, consider
$\aut\gII1,1$.  The subgroup of elements of determinant one is
conjugate to $\SL_2\Z$ and hence contains an element acting on
$\ch1$ as a triflection. This happens despite the fact
(Lemma~\tag{thm-classify-honest-reflections}) that the only
reflections of $\gII1,1$ are biflections. The following lemma
assures us that this is merely a low-dimensional phenomenon.

\beginproclaim Lemma
\Tag{thm-projective-reflections-are-honest}.
Suppose $M$ is an $n$-dimensional lattice over
$\ring=\Eisenstein$ or $\Gaussian$ and that $R$ is a projective
reflection in $\paut M$, of order $m<n$. Then $R$ is represented
by a reflection of $M$.
\endproclaim

\beginproof{Proof:}
We will also write $R$ for any element of $\aut M$ representing
$R$. 
Since $R$ acts on $\cpnm$ as a projective reflection, it has two distinct
eigenvalues $\lambda$ and $\lambda'$, with one (say $\lambda$)
having multiplicity $n-1$. Furthermore, since $R^m$ preserves
$M$ and acts trivially on $\cpnm$, we see that there is a unit
$\a$ of $\ring$ such that $\lambda^m=\lambda'^m=\a$. The
characteristic polynomial of $R$ is
$(\lambda-x)^{n-1}(\lambda'-x)$, and since $R\in \GL_n\ring$ the
coefficients  must all lie in
$\ring$. We write $y$ and $z$ for the coefficients of $x^{n-1}$ and
$x^{n-m-1}$, and compute
$$\eqalign{
{n-1\choose 1}\lambda^{\phantom{m+1}} + 
	{n-1\choose 0} \lambda'\phantom{\lambda^m}&=
	(-1)^{n-1}y\cr
{n-1\choose m+1}\lambda^{m+1} + 
	{n-1\choose m}\lambda'\lambda^{m} &= 
	(-1)^{n-m-1}z\rlap{$\;$.}\cr 
}$$
Because $\lambda^m=\a\in\ring$, the
second equation reduces to a linear equation in $\lambda$ and
$\lambda'$. For $n>m$ this is a nonsingular system of equations,
so $\lambda,\lambda'\in\ring\tensor\Q$. Since $\lambda,\lambda'$
are roots of unity they must actually lie in $\ring$. Then
$\lambda^{-1}R\in\aut M$ has eigenvalues $1$ (with multiplicity
$n-1$) and $\lambda^{-1}\lambda'$, completing the proof.
\endproof

Now we will consider the DM groups. If $\Gamma$ is a group
acting on $\chn$ then a projective reflection in $\Gamma$ is
called primitive if it is not a power of a projective reflection
in $\Gamma$ of larger order. The construction of the DM groups
allows one to find primitive projective reflections in them.
This requires a sketch of the construction, for which we use
Thurston's approach. Let $n\geq4$ and let
$\a=(\a_1,\ldots,\a_n)$ be an $n$-tuple of numbers in the
interval $(0,2\pi)$ that sum to $4\pi$. Let $P(\a)$ be the
moduli space of pairs $(p,g)$ where $p$ is an injective map from
$\{1,\ldots,n\}$ to an oriented sphere $S^2$ and $g$ is a
singular Riemannian metric on $S^2$ which is flat except on the
image of $p$, with $p(i)$ being a cone point of curvature
$\a_i$. We denote $p(i)$ also by $p_i$. Two such pairs are
considered equivalent if they differ by an
orientation-preserving similarity that identifies the
corresponding points $p_i$ with each other. This moduli space is
a manifold of $2(n-3)$ real dimensions and admits a metric which
is locally isometric to $\ch{n-3}$. Let $H$ be the group of
elements $\sigma$ of the symmetric group $S_n$ satisfying
$\a_{\sigma(i)}=\a_i$ for all $i=1,\ldots,n$. Then $H$ acts by
isometries of $P(\a)$, by permuting the points $p_i$. We denote
the quotient orbifold by $C(\a)$. The fundamental group of
$P(\a)$ is the pure (spherical) braid group on $n$ strands, and
the orbifold fundamental group of $C(\a)$ is the subgroup of the
full (spherical) braid group that maps to $H$ under the usual
map from the braid group to the symmetric group.

If the $\a_i$ satisfy certain conditions then the metric
completion $\Cbar(\a)$ of $C(\a)$ turns
out to be the quotient of $\ch{n-3}$ by a reflection group
$\G(\a)$. There are only 94 choices for $\a$ (with $n\geq5$)
satisfying these conditions, and the corresponding $\G(\a)$ are
the DM groups. The points of $\Cbar(\a)\setminus C(\a)$ are the
images of the mirrors of certain reflections of $\G(\a)$. One
can figure out the orders of the primitive reflections
associated to these mirrors by finding the ``cone angle'' at
each generic point of $\Cbar(\a)\setminus C(\a)$: if the cone
angle is $2\pi/m$ then the corresponding primitive projective
reflections have order $m$. (This cone angle should not be
confused with the cone angles at the points $p_i\in S^2$.) The
generic points of $\Cbar(\a)\setminus C(\a)$ are associated to
``collisions'' between pairs of points $p_i$ and $p_j$ on $S^2$
for which $\a_i+\a_j<2\pi$. We quote Thurston's Proposition~3.5,
which provides a way to compute the cone angles at these points
of $\Cbar(\a)\setminus C(\a)$.

\beginproclaim Proposition
\Tag{thm-cone-angles}. 
Let $S$ be the stratum of $\Cbar(\a_1,\ldots,\a_n)$ where two
cone points of $S^2$ of curvature $\a_i$ and $\a_j$ collide. If $\a_i=\a_j$
then the cone angle around $S$ is $\pi-\a_i$; otherwise it is
$2\pi-\a_i-\a_j$. 
\QED
\endproclaim 

For example, take $\a$ to be the 10-tuple
$({2\pi\over3},{2\pi\over3},{\pi\over3},{\pi\over3},{\pi\over3},{\pi\over3},{\pi\over3},{\pi\over3},{\pi\over3},{\pi\over3})$,
which is number~13 on Thurston's list and number~4 on Mostow's.
Then at the singular strata of $\Cbar(\a)$ where two cone points
of curvature ${2\pi\over3}$ (resp. two of curvature $\pi\over3$,
resp. one of each curvature) collide, the cone angle is
$\pi-{2\pi\over3}={\pi\over3}$
(resp. $\pi-{\pi\over3}={2\pi\over3}$, resp.
$2\pi-{2\pi\over3}-{\pi\over3}=\pi$). We deduce that $\G(\a)$
contains primitive projective reflections of orders~6, 3 and 2.

\beginproclaim Theorem
\Tag{thm-my-groups-are-new}.
If $L$ is $\eI n,1$ ($n\geq4$) or $\gII4n+1,1$ ($n\geq1$) then
$\reflec L$ does not appear among the Deligne-Mostow groups.
\endproclaim 

\beginproof{Proof:}
By Lemmas~\tag{thm-classify-honest-reflections} and~\tag{thm-projective-reflections-are-honest}, $\paut L$ contains no primitive
projective reflections of order~$3$ or~4. Also, $\aut L$ is not
cocompact because $L$ contains isotropic vectors. Turning to the
DM groups,  Prop.~\tag{thm-cone-angles} and the list of $n$-tuples
$\a$ provided in
\cite{mostow:monodromy-groups-on-the-complex-n-ball} or
\cite{thurston:shapes-of-polyhedra} make it easy 
to compute the cone angles at all the generic points of
$\Cbar(\a)\setminus C(\a)$ for each $n$-tuple $\a$ with $n\ge7$. The
author wrote a short computer program to do this, and also
performed the computation by hand. The only one 
for which none of the cone angles are $2\pi/4$ or $2\pi/3$ is
number~50 on Thurston's list (number~21 on Mostow's). According
to Thurston's table, $\G(\a)$ is cocompact for this choice of
$\a$.  Therefore each DM groups acting on $\chn$ for $n\geq4$ is
either cocompact or contains a primitive projective reflection
of order~3 or~4. The theorem follows.
\endproof

We close by sketching a proof that $\reflec\eI3,1$
is one of the DM groups\emdash it is the group $\G(\a)$
with
$\a=({2\pi\over3},{2\pi\over3},{2\pi\over3},{2\pi\over3},{2\pi\over3},{2\pi\over3})$,
which is number~1 on Thurston's list and number~23 on Mostow's.
Because all the $\a_i$ are equal, the orbifold fundamental group
of $C(\a)$ is the spherical braid group $B_6$ on six strands. A
standard generator for $B_6$, braiding two points $p_i$ and
$p_{i+1}$, corresponds to a loop in $C(\a)$ encircling the
singular stratum $S$ of $\Cbar(\a)$ associated to a collision
between $p_i$ and $p_{i+1}$. Since the cone angle at $S$ is
$\pi/3$ we find that the standard generators map to 6-fold
reflections. This fact, together with the braid relations and
the fact that the image of $B_6$ is not finite, specifies the
representation uniquely up to complex conjugation. The five
standard generators may be taken to map to $(-\w)$-reflections
in short roots of $\eI3,1$, which are orthogonal if the
corresponding braid generators commute and have inner product
$+1$ otherwise. One may then use the techniques of
Sections~\tag{sec-reflections} and
\tag{sec-reflection-groups} to show that 
the image of $B_6$ is $\shr\eI3,1$. The
arguments we have sketched here concerning the braid group
representation are carried out in detail in
\cite{allcock:ch4-cubic-moduli}. 

\section{References}
%

% Do not tamper with the following line
% ---begin bibentries---

\bibitem{allcock:ch13}
D.~Allcock.
 The {L}eech lattice and complex hyperbolic reflections.
 Submitted, 1999.

\bibitem{allcock:ch4-cubic-moduli}
D.~Allcock, J.~Carlson, and D.~Toledo.
 The complex hyperbolic geometry of the moduli space of cubic
  surfaces.
 Preprint, 1998.

\bibitem{allcock:ch4-cubic-moduli-announcement}
D.~Allcock, J.~Carlson, and D.~Toledo.
 A complex hyperbolic structure for moduli of cubic surfaces.
 {\it CRAS}, 326:49--54, 1998.

\bibitem{allcock:oh2}
D.~J. Allcock.
 Reflection groups on the octave hyperbolic plane.
 {\it J. Alg.}, to appear.

\bibitem{bachoc:selfdual-Hurwitz-lattices}
C.~Bachoc.
 Voisinage au sens de {K}neser pour les r\'{e}seaux quaternioniens.
 {\it Comment. Math. Helvetici}, 70:350--374, 1995.

\bibitem{bachoc:8-dimensional-Hurwitz-lattices}
C.~Bachoc and G.~Nebe.
 Classification of two genera of 32-dimensional lattices of rank 8
  over the {H}urwitz order.
 {\it Experimental Mathematics}, to appear.

\bibitem{reb:lzl}
R.~E. Borcherds.
 Automorphism groups of {L}orentzian lattices.
 {\it Journal of Algebra}, 111:133--53, 1987.

\bibitem{reb:like-leech}
R.~E. Borcherds.
 Lattices like the {L}eech lattice.
 {\it Journal of Algebra}, 130:219--34, 1990.

\bibitem{borel:arithmetic-groups}
A.~Borel and Harish-Chandra.
 Arithmetic subgroups of algebraic groups.
 {\it Annals of Mathematics}, 75:485--535, 1962.

\bibitem{jhc:26dim}
J.~H. Conway.
 The automorphism group of the 26-dimensional even unimodular
  {L}orentzian lattice.
 {\it Journal of Algebra}, 80:159--163, 1983.
 Reprinted in \cite{splag}.

\bibitem{ATLAS}
J.~H. Conway, R.~T. Curtis, S.~P. Norton, R.~A. Parker, and R.~A. Wilson.
 {\it {A}{T}{L}{A}{S} of Finite Groups}.
 Oxford, 1985.

\bibitem{jhc:leech-radius}
J.~H. Conway, R.~A. Parker, and N.~J.~A. Sloane.
 The covering radius of the {L}eech lattice.
 {\it Proceedings of the Royal Society of London}, A380:261--90, 1982.
 Reprinted in \cite{splag}.

\bibitem{jhc:laminated-lattices}
J.~H. Conway and N.~J.~A. Sloane.
 Laminated lattices.
 {\it Annals of Mathematics}, 116:593--620, 1982.
 Reprinted in \cite{splag}.

\bibitem{jhc:Coxeter-Todd}
J.~H. Conway and N.~J.~A. Sloane.
 the {C}oxeter-{T}odd lattice, the {M}itchell group, and related
  sphere packings.
 {\it Proceedings of the Cambridge Philosophical Society}, 93:421--40,
  1983.

\bibitem{splag}
J.~H. Conway and N.~J.~A. Sloane.
 {\it Sphere Packings, Lattices, and Groups}.
 Springer-Verlag, 1988.

\bibitem{coxeter:quotients-of-braid-groups}
H.~S.~M. Coxeter.
 Factor groups of the braid groups.
 In {\it Proceedings of the Fourth Canadian Mathematical Congress,
  Banff, 1957}, pages 95--122. Toronto University Press, 1959.

\bibitem{cox:gens-and-rlns}
H.~S.~M. Coxeter and W.~O.~J. Moser.
 {\it Generators and Relations for Discrete Groups}.
 Springer-Verlag, 1984.

\bibitem{deligne:monodromy-of-hypergeometric-functions}
P.~Deligne and G.~D. Mostow.
 Monodromy of hypergeometric functions and non-lattice integral
  monodromy.
 {\it Pub. I.H.E.S.}, 63:5--90, 1986.

\bibitem{feit:unimodular-Eisenstein}
W.~Feit.
 Some lattices over ${Q}(\sqrt{-3})$.
 {\it Journal of Algebra}, 52:248--263, 1978.

\bibitem{milnor:bilinear-forms}
J.~Milnor and D.~Husemoller.
 {\it Symmetric Bilinear Forms}.
 Springer-Verlag, 1973.

\bibitem{mostow:remarkable-polyhedra}
G.~D. Mostow.
 A remarkable class of polyhedra in complex hyperbolic space.
 {\it Pacific Journal of Mathematics}, 86:171--276, 1980.

\bibitem{mostow:picard-lattices-from-half-inegral-conditions}
G.~D. Mostow.
 Generalized {P}icard lattices arising from half-integral conditions.
 {\it Pub. I.H.E.S.}, 63:91--106, 1986.

\bibitem{mostow:monodromy-groups-on-the-complex-n-ball}
G.~D. Mostow.
 On discontinuous action of monodromy groups on the complex $n$-ball.
 {\it Jour. A.M.S.}, 1:555--586, 1988.

\bibitem{thurston:shapes-of-polyhedra}
W.~P. Thurston.
 Shapes of polyhedra.
 Geometry Supercomputer Project, University of Minnesota, research
  report GCG 7, 1987.

\bibitem{vin:groups-of-quad-forms}
E.~B. Vinberg.
 On the groups of units of certain quadratic forms.
 {\it Mathematics of the USSR---Sbornik}, 16(1):17--35, 1972.

\bibitem{vin:unimodular-quad-forms}
E.~B. Vinberg.
 On unimodular integral quadratic forms.
 {\it Funct. Anal. Applic.}, 6:105--11, 1972.

\bibitem{vin:I18.1-and-I19.1}
E.~B. Vinberg and I.~M. Kaplinskaja.
 On the groups ${O}_{18,1}({Z})$ and ${O}_{19,1}({Z})$.
 {\it Soviet Mathematics, Doklady}, 19(1):194--197, 1978.

% ---end bibentries---
% Do not tamper with the previous line.
%

\email{allcock@math.harvard.edu}

\homepage{http://www.math.harvard.edu/$\sim$allcock}

\address{Department of Mathematics,
Harvard University, One Oxford St, Cambridge, MA 02138}

\bye